\documentclass[11pt]{article}
\RequirePackage[OT1]{fontenc} 
\RequirePackage{amsthm,amsmath}
\RequirePackage[colorlinks,citecolor=blue,urlcolor=blue]{hyperref}

\usepackage{graphics,epsfig,epsf,psfrag}
\usepackage{url}
\usepackage{color, appendix}
\usepackage{graphicx}
\usepackage{amsmath, amsfonts, amsthm,amssymb}
\usepackage{epstopdf}
\usepackage{hyperref}
\usepackage{enumerate, framed}
\usepackage{graphicx,lscape,rotate}
\usepackage{algorithmic, algorithm}
\usepackage{epsfig}
\usepackage{mathrsfs}
\usepackage{amsfonts}
\usepackage{mathtools}
\usepackage{tasks}
\usepackage{ifthen}
\usepackage{bbm}
\graphicspath{{eps/}}
\usepackage{amsthm}
\usepackage[utf8]{inputenc}
\usepackage[english]{babel}
\usepackage{bm}
\let\BLS=\baselinestretch
\makeatletter
\newcommand{\singlespacing}{\let\CS=\@currsize\renewcommand{\baselinestretch}{1}\small\CS}
\newcommand{\doublespacing}{\let\CS=\@currsize\renewcommand{\baselinestretch}{1.5}\small\CS}
\newcommand{\normalspacing}{\let\CS=\@currsize\renewcommand{\baselinestretch}{\BLS}\small\CS}
\makeatother
\newcommand{\sK}{\mathscr{K}}

\numberwithin{equation}{section}
\theoremstyle{plain}
\newtheorem{theorem}{Theorem}[section]
\newtheorem{remark}{Remark}
\newtheorem{proposition}{Proposition}[section]
\newtheorem{corollary}{Corollary}[section]
\newtheorem{lemma}{Lemma}[section]

\newtheorem{example}[theorem]{Example}
\theoremstyle{remark}

\newcommand{\s}{\mathcal{L}}
\newcommand{\J}{ \mathscr{J}}
\newcommand{\eq}{\begin{equation}}
\newcommand{\eeq}{\end{equation}}

\newcount\colveccount
\newcommand*\colvec[1]{
        \global\colveccount#1
        \begin{pmatrix}
        \colvecnext
}
\def\colvecnext#1{
        #1
        \global\advance\colveccount-1
        \ifnum\colveccount>0
                \\
                \expandafter\colvecnext
        \else
                \end{pmatrix}
        \fi
}

\makeatletter
\newcommand{\Spvek}[2][r]{%
  \gdef\@VORNE{1}
  \left(\hskip-\arraycolsep%
    \begin{array}{#1}\vekSp@lten{#2}\end{array}%
  \hskip-\arraycolsep\right)}

\def\vekSp@lten#1{\xvekSp@lten#1;vekL@stLine;}
\def\vekL@stLine{vekL@stLine}
\def\xvekSp@lten#1;{\def\temp{#1}%
  \ifx\temp\vekL@stLine
  \else
    \ifnum\@VORNE=1\gdef\@VORNE{0}
    \else\@arraycr\fi%
    #1%
    \expandafter\xvekSp@lten
  \fi}
\makeatother

\begin{document}


\title{$\s^p$ Boundedness of the Scattering Wave Operators of Schr\"odinger Dynamics with Time-dependent Potentials and Applications -Part \uppercase\expandafter{\romannumeral1}}

\author{\textbf{Avy Soffer}$^1$\footnote{Email: \texttt{soffer@math.rutgers.edu} (Avy Soffer)},~~~~\textbf{Xiaoxu Wu}$^1$\footnote{Email: \texttt{xw292@math.rutgers.edu} (Xiaoxu Wu)}\\
1.Department of Mathematics, Piscataway, NJ 08854, USA \\
\\
}

\date{}  

\maketitle

\begin{abstract}
This paper establishes the $\s^p$ boundedness of wave operators localized at high-frequency for linear Schr\"odinger equations in $\mathbb{R}^3$ with time-dependent potentials. The approach to the proof is based on new cancellation lemmas. As a typical application based on this method, combined with Strichartz estimates is the existence and scattering for nonlinear dispersive equations. For example, we prove global existence
and uniform boundedness in $\s^{\infty}$, for a class of  Hartree nonlinear Schr\"odinger equations in $\s^2(\mathbb{R}^3),$ allowing the presence of solitons. We also prove the existence of free channel wave operators in $\s^p(\mathbb{R}^3)$ for $p>6$.
\end{abstract}

\section{Introduction}
In this paper, we let $H_0=-\Delta_x$, where $\Delta_x=(\partial/\partial x_1)^2+\cdots+(\partial/\partial x_n)^2$ is the Laplacian in $\s^2(\mathbb{R}^n), n\geq1$. The paper is devoted to the study of $\s^p$ boundedness of the wave operator $\Omega_\pm$, associated with a pair $H_0, H$ of  self-adjoint operators, and its conjugate $\Omega_\pm^*$ for all $1\leq p\leq \infty$:
\begin{equation}
\Omega_{\pm}=s\text{-}\lim\limits_{t\to\pm\infty}U(0,t)e^{-iH_0t}, \quad \text{ on }\s^p(\mathbb R^n)\cap \s^2(\mathbb R^n)
\end{equation}
\begin{equation}
\Omega_{\pm}^*=s\text{-}\lim\limits_{t\to\pm\infty}e^{itH_0}U(t,0)P_c, \quad \text{ on }\s^p(\mathbb R^n)\cap \s^2(\mathbb R^n)
\end{equation}
for the time-dependent problem
$$
i\partial_t\psi(t)=H(t)\psi(x),
$$
where the time-dependent Hamiltonian $H(t)$ is self-adjoint on the domain of $H(t)$ given by
$$
H(t)=-\Delta_x+V(x,t).
$$
Here $U(t,0)$ denotes the dynamical group of the Schr\"odinger equation generated by $H(t)$ and $P_c$ denotes the projection on the space of the scattering states of $H(t)$, the range of $\Omega_\pm$. (For example, when $H=-\Delta_x+W(x)$, $P_c$ is equal to the projection on the continuous spectrum of $H$). Here, it is worth noting that the definition of $\Omega_\pm^*$ relies on $P_c$ since without $P_c$, the strong limit of $e^{itH_0}U(t,0)$ may not exist. For example, when $H(t)$ is time-independent having a bound state $\psi_0$, $ e^{itH_0}U(t,0)\psi_0$ goes to $0$ only in the weak $\s^2$-sense. If $V(x,t)$ has sufficient decay in $x,$ and is bounded uniformly in $t,$ the existence of $\Omega_\pm$ and $\Omega_\pm^*$ follows by Cook's method and local decay estimates, and followed by density of $\s^1\cap \s^2$ in $\s^2$ (see \cite{S2018}). Based on the existence of $\Omega_\pm,$ $\Omega_\pm^*$ on a dense set of $\s^p$, we prove their $\s^p$ boundedness with a high-frequency using basic continuity in Banach Spaces. to extend the domain to the full $\s^p$ space by continuity for all $1\leq p\leq \infty$.

Throughout this paper, we focus on the case \( T \to \infty \) and work in dimension \( n = 3 \). For \( n \geq 3 \), we believe that the results can be extended using a similar argument. We adopt the notations \( \Omega \equiv \Omega_+ \) and \( \mathcal{L}^p \equiv \s^p(\mathbb{R}^3) \) for simplicity. The Fourier transform of $f(x)$ in $x$ variable in $n$-dimension is defined by
 \begin{equation}
\hat{f}(k,t):=\mathscr{F}_x[f(x)](k,t)=\frac{1}{(2\pi)^{\frac{n}{2}}}\int_{\mathbb{R}^n}e^{-ik\cdot x}f(x,t)d^nx,
\end{equation}
and
\begin{equation}
f(x,t)=\mathscr{F}^{-1}_k[\hat{f}(k,t)](x,t)=\frac{1}{(2\pi)^{\frac{n}{2}}}\int_{\mathbb{R}^n}e^{ik\cdot x}\hat{f}(k,t)d^nk.
\end{equation}

\subsection{Background and previous method}

The first general approach to the proof of these estimates  was developed by Journ\'e, Soffer, and Sogge \cite{JSS1991}. They proved decay estimates for time-independent potentials, by using a time-dependent method which combined spectral and scattering theory with harmonic analysis. Their method involved splitting solutions into high- and low-energy parts, and using Kato's smoothing and the local energy decay on the corresponding pieces. Both parts relied on cancellation lemma (CL): 
\begin{equation}
\sK_t(V(x)):=e^{iH_0t}V(x)e^{-iH_0t}: L^p \to L^p, \text{ is bounded for }1\leq p\leq \infty.
\end{equation}
Throughout this paper, we refer to $\sK_t(V(x,t))$ as the \emph{time translated}($tT$) potential for general time-dependent potentials $V(x,t)$. They also assumed that zero is neither an eigenvalue, nor a resonance, and, roughly $|V(x)|\leq C|x|^{-4-n}$, $\hat{V}\in \s^1$, where $\hat V$ denotes the Fourier transform of $V$ in $x$ variable. Recall that a resonance is a distributional solution of $H\psi=0$ so that $\psi\notin \s^2$ but $(1+|x|^2)^{-\frac{\delta}{2}}\psi(x)\in \s^2$ for any $\delta>1/2$ but not for $\delta=0$, see \cite{JK1979}.

 Their work was preceded by related results of Rauch \cite{R1978}, Jensen, Kato \cite{JK1979}, and Jensen \cite{J11980}, \cite{J21984}, who established decay estimates on weighted $\s^2$ space
\begin{equation}
\|\langle x\rangle^{-\delta}e^{-itH}f\|_{\s^2(\mathbb{R}^n)}\leq Ct^{-n/2}\|\langle x\rangle^{\delta^\prime}f \|_{\s^2(\mathbb{R}^n)}
\end{equation}
for some sufficiently large $\delta$ and $\delta^\prime$, and developed the small energy asymptotic expansions of the resolvent which are used in \cite{JSS1991} to deal with low energy estimates. Here $\langle x\rangle=\sqrt{|x|^2+1}.$

After the work of \cite{JSS1991}, many works followed. $\s^p$ estimates for wave operators were first introduced by Yajima \cite{Y21995}. He used a stationary method to prove the $\s^p$ boundedness of the wave operators, either when the Fourier transform of $\langle x\rangle^{\delta}V$ is small in some norm, or when $\partial^\alpha V/\partial x^\alpha$ decays rapidly for $|\alpha|\leq N,$ some $N\in \mathbb{N}^+.$ These assumptions on the potential are weaker than those in \cite{JSS1991}. His theorem implies the dispersive bounds by using intertwining property of the wave operators. In fact, in time-independent situation, the intertwining property holds between $H$ and $H_0$. It implies that $\Omega$ and $\Omega^*$ intertwine the part $H_c$ of $H$, the continuous spectral subspace $\s^2_c(H)$ and $H_0$: $H_c=\Omega H_0\Omega^*$ on $\s^2_c(H)$. Hence, the $\s^p$ boundedness of $\Omega$ implies that the functions $f(H_0)$ and $f(H)P_c(H)$, $P_c(H)$ being the orthogonal projection onto $\s^2_c(H)$, have equivalent operator norms from $\s^p(\mathbb{R}^n)$ to $L^{p^\prime}(\mathbb{R}^n)$ for all $1\leq p\leq 2$. However, the intertwining property holds for time-independent potentials but may fail for time-dependent systems since $U(t+s,t)$ will not generally have a nice limit as $t\to\infty$. But for potentials periodic in time with a period $T$, the intertwining property does hold since $U(t,t+T)=U(t+kT,t+(k+1)T)$, $k\in \mathbb{Z}$, see \cite{RS1980}. See also Weder \cite{W2000} for results of time-independent case in one dimension, $n = 1$, and Yajima \cite{Y31999} for $n = 2$.

    For time-dependent potentials, the analogue of Kato's scattering result was proved by Howland \cite{H11980}. When $V(x,t)$  decays in time (in the sense  of integrability), wave operators were constructed by Howland \cite{H21974} and Davies \cite{D1974}.

     For potentials periodic in $t$, Soffer and Weinstein \cite{SW1998} introduced a theory of resonances for a class of nonautonomous Hamiltonians to treat the problem related to time-periodic potentials. A further consequence of the $\s^p$ decay estimates is the Strichartz estimate \cite{JSS1991}. The non-endpoint Strichartz estimates (when $q\neq 2$ in~\eqref{str}) were addressed in \cite{GV1992}, \cite{Y11987} and of course the original work of Strichartz \cite{Str1977}. The more delicate endpoint cases are established by Keel and Tao \cite{KT1998}.

Closely related to the boundedness of the wave operator on $\s^p,$ are $\s^p$ decay estimates for the free Schr\"odinger equation $(H(t)=H_0)$ on $\mathbb{R}^n$ given by
\begin{equation}
\| e^{itH_0}f\|_{\s^p}\leq C_p |t|^{-n(\frac{1}{2}-\frac{1}{p})}\| f\|_{\s^{p^\prime}}, \quad p\geq 2, \frac{1}{p}+\frac{1}{p^\prime}=1.
\end{equation}
They imply the Strichartz estimates
\begin{equation}
\| e^{itH_0}f \|_{\s^q_t\s^r_x}\leq C_q \| f\|_{\s^2} , \quad 2\leq r, q\leq \infty, \frac{n}{r}+\frac{2}{q}=\frac{n}{2}, \text{ and }(q,r,n)\neq (2,\infty,2).\label{str}
\end{equation}
The proof for the non-endpoint case can be found in \cite{T2006}.

Such decay estimates play a fundamental role in the theory of nonlinear dispersive equations, among other things. The extension of such theories to inhomogeneous problems (either due to curvature, local potentials, or coherent structure such as solitons, vortices, etc.) then motivated the efforts to establish the $\s^p$ decay estimates for more general Hamiltonians. \par

     Rodnianski and Schlag \cite{RS2004} proved decay estimates for small time-dependent potentials which  satisfy the following condition
\begin{equation}
\sup\limits_{t}\| V(t,\cdot)\|_{\s^{3/2}(\mathbb{R}^3)}+\sup\limits_{y\in \mathbb{R}^3}\int_{\mathbb{R}^3}\int \frac{|V(\hat{\tau},x)|}{|x-y|}d\tau dx<c_0,\text{ for some constant }c_0>0.
\end{equation}
Their proof uses the representation of $U(t,0)$ as an infinite series of oscillatory integrals; they also established non-endpoint Strichartz estimates for large time-independent potentials with $\langle x\rangle^{-2-\epsilon}$ decay.

 Goldberg proved, in \cite{G12006} dispersive estimates for almost-critical potentials and, in \cite{G22009}, Strichartz estimates for $\s^{n/2}$ and thus scaling-critical potentials. Later, Beceanu \cite{B2011} proved Strichartz estimates for time-dependent potentials by using new forms of Wiener's theorem.  \par
Now we go back to the wave operator.
The construction of wave operators, and in particular the use of the intertwining property has a long history, going back at least to Friedrich. But the application to the case where the potential perturbation is time-dependent is largely unknown.

 \subsection{Improved cancellation lemma}
We introduce an {\bf improved cancellation lemma} (ICL): 
\eq
I\sK(V):=\int_0^\infty dt e^{itH_0}V(x,t)e^{-itH_0}: \s^p_x\to \s^p_x, \text{  is bounded for }1\leq p\leq\infty. \label{Dec.14.1}
\eeq
We refer to $I\sK(V)$ as the integrated $tT$ potential. Throughout the paper, we write $I\sK$ to represent $I\sK(V)$ for convenience. We also refer to the $\s^p$ boundedness of $I\sK(V)$ as the advanced cancellation lemma.

We prove $\s^p$ boundedness of wave operator $\Omega$ on high frequency part of $\s^p$ space using ICL. To be precise, let $\beta(\lambda)$ be a smooth cut-off function satisfying $\beta(\lambda)=0$ for $-\infty<\lambda<1/2$ and $\beta(\lambda)=1$ for $\lambda\geq 1$ and write $\beta(t>M):=\beta(\frac{t}{M})$ for $M>0$. The wave operator on high-frequency cut-off $\s^p$ space denoted by $\Omega\beta(H_0>M)$ is understood in the sense of Abelian limits:
\begin{equation}
\Omega\beta(H_0>M)=s\text{-}\lim\limits_{\epsilon\downarrow0}\Omega_\epsilon\beta(H_0>M), \text{ on }\s^p,\label{def: omMe}
\end{equation}
where $\Omega_\epsilon, \epsilon>0,$ is given by  
\begin{equation}
\Omega_\epsilon=1+i\int_0^\infty dt e^{-\epsilon t}\Omega(t) e^{itH_0}V(x,t)e^{-itH_0}, \quad \Omega(t)=:U(0,t)e^{-itH_0}.
\end{equation}
By applying Duhamel's principle and iterating it infinitely many times on the right-hand side of~\eqref{def: omMe}, we obtain  
\begin{equation}  \label{expan: Omegae}
    \Omega_\epsilon = 1 + \sum_{j=1}^\infty i^j I_\epsilon^{(j)},  
\end{equation}  
where \( I_\epsilon \equiv I_\epsilon^{(1)} \) is defined as  
\begin{equation}  
    I_\epsilon := \int_0^\infty e^{-\epsilon t} \sK_t(V(x, t)) \, dt,  
\end{equation}  
and for \( j \geq 2 \), \( I_\epsilon^{(j)} \) is given by  
\begin{equation}  
    I_\epsilon^{(j)} := \int_0^\infty \int_0^{t_1} \cdots \int_0^{t_{j-1}} e^{-\epsilon t_j} \sK_{t_j} \cdots \sK_{t_1} \, dt_j \cdots dt_1,  
\end{equation}  
where we use the shorthand notation \( \sK_t \equiv \sK_t(V(x, t)) \). It is worth noting that $ I\mathscr{K}=s\text{-}\lim\limits_{\epsilon \downarrow 0} I_\epsilon$ on $\s^p$, for $1\leq p\leq \infty$, provided that the strong limit makes sense. Additionally, thanks to factor $e^{-\epsilon t}$, by the CL introduced in \cite{JSS1991} $I_\epsilon$ is bounded on $\s^p$ for all $\epsilon>0, 1\leq p\leq \infty$. We prove the $\s^p$ boundedness of $I_\epsilon$ for all $1\leq p\leq \infty$ uniformly in $\epsilon\in [0,1)$ and this establishes the ICL defined in~\eqref{Dec.14.1}. Additionally, we find that the ICL and the high-frequency cutoff \( \beta(H_0 > M) \) yield the following estimates:  
\begin{equation}  
    \| I_\epsilon^{(j)} \beta(H_0 > M) \|_{\s^p \to \s^p} \leq \frac{C^j}{M^{\frac{j-1}{2}}}, \qquad j = 2, 3, \dots,  \epsilon\in [0,1)\label{est: Iej}
\end{equation}  
for some constant \( C = C(V) > 0 \), provided \( V \) satisfies the assumed conditions introduced in subsection~\ref{sec: main}. The estimates in~\eqref{est: Iej}, combined with the ICL and~\eqref{expan: Omegae}, establish the \(\s^p\)-boundedness of \(\Omega_\epsilon \beta(H_0 > M)\) for all \(1 \leq p \leq \infty\) and \(\epsilon \in [0, 1]\), provided \(M > C^2\), where \(C\) is the constant from~\eqref{est: Iej}. Consequently, \(\Omega \beta(H_0 > M) \equiv \Omega_0 \beta(H_0 > M)\) is bounded on \(\s^p\) for all \(1 \leq p \leq \infty\). By duality, we obtain the $\s^p$ boundedness of $\beta(H_0>M)\Omega^*$.

Unfortunately, this scheme fails for the low-frequency component of the data. In fact, for the low-frequency part, additional challenges arise due to the possible presence of resonances. The proof of the \(\s^p\)-boundedness of the wave operator for the low-frequency part will be addressed in future work.

\subsection{Main results and applications to NLS equations}\label{sec: main}

We outline several key cases of time-dependent potentials that are the focus of our study. Let $\{e_1,e_2,e_3\}$ be an orthogonal basis in $\mathbb R^3$. We first consider a class of Mikhlin-type potentials (in the $t$ variable) $V(x,t)$. We assume $V(x,t)$ satisfies the condition that there exists $\hat{V}_0(\xi)\in \s^1_\xi (\mathbb{R}^3)\cap \s^\infty_\xi(\mathbb{R}^3)$ such that
\begin{equation}
\sup\limits_{t\in \mathbb{R}}\frac{(1+|t|)^a}{a!}\sum\limits_{l,j=0}^2\sum\limits_{m,r=1}^3|\frac{\partial^a}{\partial t^a}\left[\partial_{\xi\cdot e_r}^l\partial_{\xi\cdot e_m}^j\hat{V}(\xi,t)\right]|\leq c^a |\hat{V}_0(\xi)|, \text{ for all }a\in \mathbb{N}, \text{ some }c\geq 1.\label{A.19.1.1}
\end{equation}

\begin{theorem}\label{main3}If \( V(x, t) \) satisfies condition~\eqref{pp1}, then there exists \( M = M(V(x, t)) > 0 \) such that both \(\Omega \beta(H_0 > M)\), as defined in~\eqref{def: omMe}, and its adjoint exist on \(\s^p\) for all \(1 \leq p \leq \infty\) and are bounded.

\end{theorem}
The proof of Theorem~\ref{main3} is presented in section \ref{section 3}. Typical examples are
\begin{equation}
V(x,t)=V_0(x)+\frac{\sin(\ln(1+|t|))}{(1+|t|)^\delta}V_1(x), \text{ for }\delta\geq 0,
\end{equation}
and
\begin{equation}\label{exp 2}
V(x,t)=V_0(x)+V_1(x-\frac{\sin(\omega \ln (|t|+1))}{(1+|t|)^\delta}v),\text{ for }\delta\geq 0,
\end{equation}
where $V_j(x), j=0,1,$ are real-valued functions satisfying the condition~\eqref{A.19.1.1}. Please refer to Corollaries \ref{ext1} and \ref{timesum1} for a detailed discussion.
\begin{remark}
The first example above involves a potential that decays arbitrarily slowly in time (\(\delta > 0\)) or has no decay at all (\(\delta = 0\)), approaching a time-independent potential. Since the decay in time is not in \(\s^1\), this case is not addressed by existing results, even in \(\s^2\). See, for example, \cite{SW1999}. The second example is more involved: it corresponds to a charge transfer type Hamiltonian, where the moving potential is along a non-linear path in time. Previous works required the path to be linear up to a fast decaying term. The case of a general path, which however converges to an end point, was considered in \cite{BS2012};    the path was allowed to be a rough function of time. However, this method did not apply to the charge transfer case defined in~\eqref{exp 2}, as a time-independent part $V_0$ was not allowed. All previous works focused on proving time decay estimates of the dynamics, but not $\s^p$ boundedness of the wave operator. See e.g. \cite{RSS2005}, \cite{B2011}.
\end{remark}

\begin{remark}When \(\partial_t[V](t, \xi) \in \s^1_t(0, \infty)\), it indicates that the asymptotic energy exists and is bounded. This may suggest that \(\delta > 0\) is optimal. Interestingly, our method is capable of addressing the case when \(\delta = 0\). However, in this scenario, it remains an open question whether the frequency support of the solution generally stays bounded.
\end{remark}

\bigskip

We also consider the case of self similar potentials:
\begin{equation}
V(x,t)=V_1(g(t)x,t)+\frac{1}{(2\pi)^{3/2}}\sum\limits_{j=1}^\infty f_j(t)e^{i g_j(t)x\cdot a_j}\label{type21}
\end{equation}
where $g_j(t)$ are real-valued functions on $[0,\infty)$ and $V_1$ and $f_j, j=1,\cdots,$ satisfy the condition 
\begin{equation}
\int d^n\xi |\hat{V}_1(\xi,t)|+\sum\limits_{j=1}^\infty |f_j(t)|\in \s^1_t[0,\infty).\label{710.3}
\end{equation}
Let 
\begin{equation}
    h(t):=\int d^n\xi |\hat{V}_1(\xi,t)|+\sum\limits_{j=1}^\infty |f_j(t)|.
\end{equation}
\begin{theorem}\label{dilationg} Let $V(x,t)$ be as in~\eqref{type21}. If $g_j(t)$ are real-valued functions on $[0,\infty)$ and condition~\ref{710.3} holds true, then
\begin{equation}
\lim\limits_{t\to \pm \infty}\| U(0,t)e^{-itH_0}-\Omega \|_{\s^p\to L^p}= 0 
\end{equation}
and
\begin{equation}
    \|\Omega\|_{\s^p\to L^p}\leq \exp\left(\frac{\| h(t)\|_{\s^1_t(0,\infty)}}{(2\pi)^{\frac{3}{2}}}\right).
\end{equation}
\end{theorem}
A typical example is  
\begin{equation}  
V(x, t) = \frac{\chi(|t| \geq c)\sin(\omega t)}{t^{3/2}}V_1\left(\frac{x}{t}\right), \quad \text{for some } c > 0, \, \omega \in \mathbb{R},  
\end{equation}  
where \(\hat{V}_1(\xi)\) is a finite measure. This formulation can be applied to study self-similar solutions of certain nonlinear Schrödinger equations (NLS). For a detailed proof, refer to Section~\ref{section 4}.\par
Theorems~\ref{main3} and~\ref{dilationg} imply the corresponding $\s^p$ decay estimates and Strichartz estimates for the high frequency part in  $\s^p$ space. To be precise, let
\begin{equation}
\Omega(0,t):=U(0,t)e^{-iH_0t}, \qquad t\geq0.\label{1.23}
\end{equation}
We observe that the \(\s^p\)-boundedness of \(\Omega \beta(|H_0| > M^2)\) implies 
\begin{center}
the \(\s^p\)-boundedness of \(P_c \Omega(0, t) \beta(|H_0| > M^2)\), uniformly in \(t \in \mathbb{R}^+\), 
\end{center}
provided that \(P_c\) is bounded on \(\s^p\) for all \(1 \leq p \leq \infty\). See Corollaries~\ref{3Tu} and~\ref{4Tu}. This further leads to \(\s^{p'}\)-decay estimates via the inequality  
\begin{equation}  
\| P_c U(0, T) \beta(|H_0| > M) \|_{\s^{p} \to \s^{p^\prime}} \leq \sup_{T \in \mathbb{R}} \| P_c \Omega(0, T) \beta(H_0 > M^2) \|_{\s^{p^\prime} \to \s^{p^\prime}} \| e^{itH_0} \|_{\s^p \to \s^{p^\prime}},  
\end{equation}  
for all \(p \in [1, 2]\).\par
\bigskip

We also establish decay estimates directly using the same strategy for the charge transfer case, where the potential is \(V(x - \sqrt{1 + |t|}v)\). Here, the real-valued function \(V(x)\) satisfies the condition  
\begin{equation}  
||| V(x)|||_p := \sum_{l=0}^2 \sum_{m=1}^3 \| (|\xi| + 1)^3 |\partial_{\xi \cdot e_m}^l \hat{V}(\xi)| \|_{\s^1_\xi} + \|V(x)\|_{\s^1_x \cap \s^2_x} < \infty. \label{sqrtc}  
\end{equation}  

\begin{theorem}If $V(x)$ satisfies the condition~\ref{sqrtc}, then for a sufficiently large $M>0$,
\begin{equation}
\sup\limits_{T\in \mathbb{R}} |T|^{3/2}\|P_cU(0,T)\beta(|P|>M)\|_{\s^1_x\to \s^\infty_x}<\infty.
\end{equation}
\end{theorem}

Such type of potentials may be observed in the motion of vortices for example.
See section \ref{section 4} for more details. \par

\bigskip

The methods developed here may be applied  to NLS Equations.  Let
\eq
\mathcal{F}\s^1_x:=\left\{ f(x): \hat{f}(\xi)\in \s^1_\xi(\mathbb{R}^3)\right\}.
\eeq
We use advanced CL to deal with Hartree-type NLS equations
\begin{equation}
i\partial_t\psi(t)=(H_0+V(x,t))\psi(t)+\mathcal{N}(|\psi(t)|)\psi(t), \quad \psi(0)=\psi_0\in \s^2(\mathbb{R}^3) \label{generalNLS}
\end{equation}
with $\mathcal{N}(\cdot): \s^2_x\cap \s^p_x\to \s^2_x\cap \mathcal{F}\s^1_x $ for some $2\leq p<6$, satisfying the following \textbf{advanced cancellation criterion(ACC)} and the condition: for some $(q,r)$ satisfying
\eq
\frac{2}{q}+\frac{3}{r}=\frac{3}{2}, \quad 2\leq q\leq \infty, 2\leq r<6\label{strr},
\eeq
\begin{enumerate}[(A)]
    \item \label{Ncondition1}(\textbf{ACC1}): For $\psi(t)\in \s^q_t([-T,T])\s^r_x\cap C_t\s^2_x$, all $1\leq p\leq\infty$, some $k_1> 1$,
 \eq
\|(\mathcal{N}(|\psi(t)|))\|_{\s^{k_1}_t([-T,T])\mathcal{F}\s^1_x}\lesssim  C(\| \psi(t)\|_{\s^q_t([-T,T])\s^r_x\cap C_t\s^2_x}).\label{N1}
\eeq

    \item \label{Ncondition2} (\textbf{ACC2}): For $\psi(t), \phi(t)\in \s^q_t([-T,T])\s^r_x\cap C_t\s^2_x$, all $1\leq p\leq\infty$,
    \begin{multline}
   \int_{-T}^T dt\|\mathcal{N}(|\psi(t)|)-\mathcal{N}(|\phi(t)|)\|_{\mathcal{F}\s^1_x} \lesssim C(T)\| \psi(t)-\phi(t)\|_{\s^q_t([-T,T])\s^r_x\cap C_t\s^2_x}  \times \\
    C(\| \psi(t)\|_{\s^q_t([-T,T])\s^r_x\cap C_t\s^2_x},\|\phi(t)\|_{\s^q_t([-T,T])\s^r_x\cap C_t\s^2_x})\label{N2}
    \end{multline}
    with some constant $C(T)$ satisfying 
    \eq
C(T)\to 0, \text{ as }T\to 0.
\eeq
    \item\label{Ncondition3} (\textbf{Condition}):
    \eq
    \| \mathcal{N}(|f(x)|)f(x)\|_{\s^1_x}\lesssim \|f(x)\|_{\s^r_x}^{q_0}
    \eeq
    with
    \eq
    0\leq q_0\leq q.
    \eeq
\end{enumerate}
Here the potential $V(x,t)$ satisfies following \textbf{advanced cancellation criterion} and some condition:
\begin{enumerate}
\item \label{Vcondition1}(\textbf{ACC3}): For all $1\leq p\leq\infty$, any $a\in \mathbb{R}$, some $k_2> 1$,
\eq
\| V(x,t+a)\|_{\s^{k_2}_t([-T,T])\mathcal{F}\s^1_x}\lesssim_T 1.
\eeq
    \item\label{Vcondition2} (\textbf{Condition}): for any $a,T\in \mathbb{R},$
    \eq
    \| V(x,t)\|_{ \s^{q_1}_t([a,a+T])\s^{r'}_x}\lesssim_T 1 
    \eeq
    with
    \eq
    \frac{1}{r'}+\frac{1}{r}=1, \quad  \frac{1}{q'}+\frac{1}{q}=1, \quad q_1\geq q'.
    \eeq
\end{enumerate}
\begin{theorem}\label{HartreeNLS}If $V(x,t)$ satisfies the conditions \ref{Vcondition1} and \ref{Vcondition2} and if $\mathcal{N}$ satisfies the conditions \eqref{Ncondition1}-\eqref{Ncondition2}, then system \eqref{generalNLS} has global wellposedness in $\s^2_x$ and in addition, if $\psi_0\in \s^1_x\cap \s^2_x$ and $\mathcal{N}$ also satisfies the condition~\eqref{Ncondition3}, then for any $c>0$,
\eq
\sup\limits_{|t|\geq c}\| \psi(t)\|_{\s^\infty_x}\lesssim_{\|\psi_0\|_{\s_x^1\cap \s^2_x},c} 1.
\eeq
\end{theorem}

\begin{remark}Here for global wellposedness, we believe that $k_1$ in \eqref{Ncondition1} can be equal to $1$.
\end{remark}
The proof for Theorem \ref{HartreeNLS} relies on advanced CL by using advanced cancellation criterion. Based on such advanced CL for $\mathcal{N}(|e^{-itH_0}\psi_0|)$, an adapted iteration scheme and standard contraction mapping argument, we get local wellposedness in $\s^2_x$ and local Strichartz estimate for solution $\psi(t)$. Based on such result, we are able to build advanced CL for $\mathcal{N}(|\psi(t)|)$, which helps to establish the $\s^\infty_x$ boundedness for $\psi(t)$ when $|t|\geq 1$. Such upper bound is independent on $t\in (-\infty, -c]\cup [c,\infty)$ with given $c>0$. Typical examples are
\eq
\mathcal{N}(|\psi(t)|)=\pm \lambda[\frac{1}{|x|^{3/2-\delta}}*|\psi(t)|^2](x), \text{ for }\delta\in (0,\frac{3}{2}), \lambda>0\label{ex11}
\eeq
and
\eq
\mathcal{N}(|\psi(t)|)=\pm \lambda[\frac{e^{-c|x|}}{|x|^{3/2-\delta}}*|\psi(t)|^2](x), \text{ for }\delta\in (0,\frac{3}{2}), \lambda>0, c>0.\label{ex12}
\eeq
For \eqref{ex11}, we establish global well-posedness, and for \eqref{ex12}, we prove both global well-posedness and $\s^\infty$ boundedness for $|t|\geq c$ for any $c>0$. To illustrate the theory, we also prove Theorem \ref{HartreeNLS} by showing that how the method works in an example:
\begin{theorem}\label{2020LinftyNLS}System 
\begin{equation}
\begin{cases}
i\partial_t\psi(t)=H_0\psi(t)+[f*|\psi(t)|^2](x)\psi(t), \\
\quad \psi(0)=\psi_0\in \s^2(\mathbb{R}^3)
\end{cases}, \quad \text{ with }f(x,t)\in C_t\s^2_x \label{generalNLS00}
\end{equation}
has global wellposedness in $\s^2_x$ and in addition, if $\psi_0\in \s^1_x\cap \s^2_x$, then for any $c>0$,
\eq
\sup\limits_{|t|\geq c}\| \psi(t)\|_{\s^\infty_x}\lesssim_{\|\psi_0\|_{\s_x^1\cap \s^2_x},c} 1.
\eeq
\end{theorem}
When it comes to $\mathcal{H}^1_x$, we consider the following NLS
\eq
\begin{cases}
i\partial_t\psi(x,t)=(-\Delta_x+\mathcal{N}(|\psi(x,t)|))\psi(x,t)\\
\psi(x,0)=\psi_0(x)\in \mathcal{H}^1_x(\mathbb{R}^3)\cap \s^1_x(\mathbb{R}^3)
\end{cases}, \quad \text{ in }3 \text{ dimensions}\label{NLS}
\eeq
where $\mathcal{H}^1_x$ denotes the Sobolev space with integer $1$. We show the $\s^{p}$ boundedness of $e^{itH_0}U(t,0)-1$(including $\Omega_\pm^*-1$) on $\s^{p_0}_x\cap \mathcal{H}^1_x$ for any $p_0\in (6,\infty], p\in [2,\infty]$ if $\psi_0\in \mathcal{H}^1_x$ leads to a global solution with $\mathcal{H}^1_x$ uniformly bounded in $t$ and if $\mathcal{N}$ satisfies the condition 
\eq
\begin{cases}
\mathcal{N}(\cdot): \mathcal{H}^1_x\to \s^2_x, \text{ is bounded}\\
\mathcal{N}_1(\cdot): \mathcal{H}^1_x\to \s^2_x, \text{ is bounded}\\
\mathcal{N}'(\cdot): \mathcal{H}^1_x\to \s^3_x, \text{ is bounded}
\end{cases}\label{condition}
\eeq
where
\eq
\mathcal{N}'(k):=\frac{d}{dk}[\mathcal{N}(k)], \quad \mathcal{N}_1(k)=\frac{\mathcal{N}(k)}{|k|}:
\eeq
\begin{theorem}[Existence of free channel wave operator in $\s^p_x$]\label{LinftyNLS}For any $p\in [2,\infty ], p_0\in (6,\infty]$, if $\mathcal{N}$ satisfies \eqref{condition} and if 
\eq
\sup\limits_{t\in \mathbb{R}}\|\psi(t)\|_{\mathcal{H}^1_x}\leq C(\|\psi_0\|_{\mathcal{H}^1_x}),\label{Dec.24.1}
\eeq
then 
\eq
\| (e^{itH_0}U(t,0)-1)\psi_0\|_{\s^p_x}\leq C(\|\psi_0\|_{\mathcal{H}^1_x\cap \s^{p_0}_x}, \sup\limits_{t\in \mathbb{R}}\|\psi(t)\|_{\mathcal{H}^1_x} ).
\eeq
Furthermore, if we also have
\eq
\|\mathcal{N}(|f(x)|)f(x)\|_{\s^{p'}_x}\lesssim_{\|f(x)\|_{H^1_x}} 1,\quad\text{ for some }p\in(6,\infty]
\eeq
then for $\psi_0\in \mathcal{H}^1_x\cap \s^{p}_x$ satisfying \eqref{Dec.24.1}, for $p>6$,
\eq
\Omega_\pm^*\psi_0:= \lim\limits_{t\to \pm \infty} e^{itH_0}U(t,0)\psi_0\text{ exists in }\s^p_x
\eeq
and
\eq
\|\Omega_\pm^*\psi_0 \|_{\s^p_x}\leq C(\|\psi_0\|_{\mathcal{H}^1_x\cap \s^p_x}, \sup\limits_{t\in \mathbb{R}}\|\psi(t)\|_{\mathcal{H}^1_x} ).
\eeq

\end{theorem}
{
\begin{remark}Here $p>6$ makes $e^{itH_0}: \s^{p'}_x\to \s^p_x$, bounded with a bound $|t|^{-3(\frac{1}{2}-\frac{1}{p'})}$ integrable on $\mathbb{R}\setminus(-1,1)$. We will give a proof for the case when $p=\infty$ and the result for other $p\in (6,\infty)$ will follow in a similar way.
\end{remark}

In addition, if we only have $\psi_0\in \mathcal{H}^1_x$, we can have $\s^p$ boundedness of $e^{itH_0}U(t,0)-1$ for $2\leq p<\infty$:
\begin{theorem}\label{LinftyNLSp}For any $p\in [2,\infty )$, if $\mathcal{N}$ satisfies \eqref{condition} and if 
\eq
\sup\limits_{t\in \mathbb{R}}\|\psi(t)\|_{\mathcal{H}^1_x}\lesssim_{\|\psi_0\|_{\mathcal{H}^1_x}}1,
\eeq
then $e^{itH_0}U(t,0)-1: \mathcal{H}^1_x\to \s^p_x$, is bounded uniformly in $t\in (-\infty,-1]\cup[1,\infty)$. Moreover, if $\psi_0\in \s^\infty_x\cap \mathcal{H}^1_x$, then $e^{itH_0}U(t,0)\psi_0\in \s^\infty_x$.
\end{theorem}
The proof for Theorem \ref{LinftyNLS} mainly relies on $\s^\infty$ boundedness of $ e^{itH_0}U(t,0)-1$ on $ \mathcal{H}^1_x\cap \s^{p_0}_x$ since $ e^{itH_0}U(t,0)-1$ is already bounded on $\mathcal{H}^1_x$ and since the $\s^p$ result can follow via interpolation inequality. And the $\s^\infty $ boundedness relies on the method of $ItT$ potential(advanced CL). The proof for Theorem \ref{LinftyNLSp} mainly relies on the statement that if $\psi_0\in \mathcal{H}^1_x$, then $\psi(t)\in \s^\infty_x+\mathcal{F}\s^{1+\epsilon}_x$ for any $\epsilon \in (0,1)$.

Here are some examples: 
\eq
\mathcal{N}(f):=|f|^3,\text{ or } -|f|^2+|f|^3,
\eeq
the assumption \eqref{condition} is satisfied:
When $\mathcal{N}(f)=|f|^3,$ that $\psi_0\in \mathcal{H}^1_x$ implies global wellposedness in $\mathcal{H}^1_x$ due to energy conservation 
\eq
E(\psi(t)):=\int d^nx (\frac{1}{2}|\nabla_x\psi(t)|^2+\frac{1}{2}|\psi(t)|^4).
\eeq
When $\mathcal{N}(f)= -|f|^2+|f|^3$, we have following lemma:
\begin{lemma}[\cite{TVZ2007}]\label{TVZ}If $ \psi_0(x)\in \mathcal{H}^1_x$, then with $\mathcal{N}(f)= -|f|^2+|f|^3$, 
\eq
\|\psi(t)\|_{\dot{S}(I\times \mathbb{R}^3)}\lesssim C(|I|, \|\psi_0 \|_{\mathcal{H}^1_x}).
\eeq
\end{lemma}

Our method also has some other applications, e.g. the ionization problem for more general potentials \cite{SW1999}. Decay estimates of $\beta(|P|\leq M)P_cU(t,0)$ with rough potentials will be treated in a future publication.
In Theorem \ref{LinftyNLS}, the solution is not always dispersive, due to the possible presence of solitons or other bound states.

\subsection{Other Results}
In this subsection, we discuss the result when \(V\) is time-independent. Let \(L_{\eta, l, j}(k, \hat{\xi}, \epsilon)\) denote the Fourier transform of \(\chi(|\xi| \geq 0) |\xi| \partial_{\xi \cdot e_l}^j [\hat{V}(\xi - \eta)] e^{- \frac{\epsilon}{|\xi|}}\) in the \(|\xi|\) variable, where \(l = 1, 2, 3\), \(j = 0, 1, 2\), \(\eta \in \mathbb{R}^3\), and  
\begin{equation}  
K_1(V(x), \eta) := \max_{l=1,2,3, \, j=0,1,2} \int_{S^2} d\sigma(\xi) \int_{-\infty}^\infty dk \sup_{\epsilon \geq 0} |L_{\eta, l, j}(k, \hat{\xi}, \epsilon)|.  
\end{equation}  

\begin{theorem}\label{intime}If
\begin{equation}
\hat{V}_a(\xi):=\sum\limits_{j,l=0}^2\sum\limits_{r,m=1}^3|\partial_{\xi\cdot e_r}^j\partial_{\xi\cdot e_m}^l\hat{V}(\xi)|\in \s^1_\xi\text{ and }K_m(V(x)):=\sup\limits_{\eta\in \mathbb{R}^3}|K_1(V(x),\eta)|<\infty, \label{intimec}
\end{equation}
for $l=0,1,2,$ then there exists $M=M(V(x))>0$ such that both $\Omega\beta(H_0>M)$ and its adjoint exist on $\s^p$ and are bounded.
\end{theorem}
 Here the assumption $K_m(V(x))<\infty$  is satisfied if $\langle|P_\xi|\rangle^2[\hat{V}](\xi)\in \s_\xi^\infty$, where $(P_\xi)_j:=-i\partial_{\xi_j}$:
\begin{proposition}If  $ \langle|P_\xi|\rangle^2[\hat{V}](\xi)\in \s_\xi^\infty $ and
\begin{equation}
\| \langle |P_\xi|\rangle^4[\hat{V}](\xi)\|_{\mathcal{K}(\mathbb{R}^3)}<\infty
\end{equation}
where $\|\cdot \|_{\mathcal{K}(\mathbb{R}^3)}$ denotes Kato norm, then $ K_m(V)<\infty$.
\end{proposition}
\begin{proof} We begin with the case where \(l = 1\), \(j = 0\), and \(\epsilon = 0\):  
\[
K^{1,0} := \int_{S^2} d\sigma(\xi) \int_{-\infty}^\infty dk \left| \int_0^\infty d|\xi| \, |\xi| \hat{V}(\xi - \eta) e^{ik|\xi|} \right|.  
\]  
For \(|k| \leq 1\), we proceed directly without additional steps. For \(|k| > 1\), we take integration by parts twice with respect to the \(|\xi|\) variable:  
\begin{align*}  
|K^{1,0}| \leq & \int_{S^2} d\sigma(\xi) \int_{-1}^1 dk \int_0^\infty d|\xi| \, |\xi| |\hat{V}(\xi - \eta)| + \\  
& \int_{S^2} d\sigma(\xi) \int_{|k| > 1} \frac{dk}{k^2} \int_0^\infty d|\xi| \left| \partial_{|\xi|}^2 \left[ |\xi| \hat{V}(\xi - \eta) \right] \right| + \\  
& \int_{S^2} d\sigma(\xi) \int_{|k| > 1} \frac{dk}{k^2} |\hat{V}(-\eta)| \quad \text{(Boundary term)}.  
\end{align*}  

Using Fubini's theorem and changing from spherical to standard Euclidean coordinates, we obtain:  
\begin{align*}  
|K^{1,0}| \leq & \int_{-1}^1 dk \int d^3\xi \frac{|\hat{V}(\xi - \eta)|}{|\xi|} + \int_{|k| > 1} \frac{dk}{k^2} \int d^3\xi \frac{\left| \partial_{|\xi|}^2 \left[ |\xi| \hat{V}(\xi - \eta) \right] \right|}{|\xi|^2} + \\  
& \int_{S^2} d\sigma(\xi) \int_{|k| > 1} \frac{dk}{k^2} |\hat{V}(-\eta)| \\  
\leq & 2 \| \langle |P_\xi| \rangle^4 [\hat{V}](\xi) \|_{\mathcal{K}(\mathbb{R}^3)} + (2 \| \langle |P_\xi| \rangle^4 [\hat{V}](\xi) \|_{\mathcal{K}(\mathbb{R}^3)} + 8\pi \| \hat{V}(\xi) \|_{\s^\infty_\xi}) + 8\pi \| \hat{V}(\xi) \|_{\s^\infty_\xi}.  
\end{align*}  

Similarly, we find:  
\begin{equation}  
|K_1(V(x), \eta)| \leq 4 \| \langle |P_\xi| \rangle^4 [\hat{V}](\xi) \|_{\mathcal{K}(\mathbb{R}^3)} + 16\pi \| \langle |P_\xi| \rangle^2 \hat{V}(\xi) \|_{\s^\infty_\xi},  
\end{equation}  
and consequently:  
\begin{equation}  
K_m(V) = \sup_{\eta \in \mathbb{R}^3} |K_1(V(x), \eta)| \leq 4 \| \langle |P_\xi| \rangle^4 [\hat{V}](\xi) \|_{\mathcal{K}(\mathbb{R}^3)} + 16\pi \| \langle |P_\xi| \rangle^2 \hat{V}(\xi) \|_{\s^\infty_\xi}.  
\end{equation}  

\end{proof}

\section{Improved Cancellation Lemma (CL)}\label{section C}
We begin by introducing additional notation that will be used throughout this paper, followed by the definitions of the CL and the improved CL.

\subsection{Notation}
In this paper, even though in most situation, $n=3$, $n$ will always denote the dimension of the ambient physical space, the configuration space. If $x=(x_1,\cdots,x_n)$ and $\xi=(\xi_1,\cdots,\xi_n)$ lie in $\mathbb{R}^n$, we use $x\cdot\xi$ to denote the dot production $x\cdot \xi:=x_1\xi_1+\cdots +x_n\xi_n$, and $|x|$ to denote the magnitude $|x|:=(x_1^2+\cdots+x_n^2)^{1/2}$. We also use $\langle x\rangle$ to denote the inhomogeneous magnitude ( Japanese $x$) $\langle x\rangle:=(1+|x|^2)^{1/2}$ of $x$. The derivatives will either be interpreted in the classical sense or the distributional sense. \par
If $X$ and $Y$ are two quantities, we use $X\lesssim Y$ to denote the statement that $X\leq CY$ for some absolute constant $C>0$. More generally, given some parameters $a_1,\cdots, a_k$, we use $X\lesssim_{a_1,\cdots,a_k} Y$ to denote the statement that $X\leq C_{a_1,\cdots,a_k} Y$ for some constant $C_{a_1,\cdots,a_k}>0$ which can depend on the parameters $a_1,\cdots, a_k$. \par
Throughout the paper, $P_j:=-i\partial_{x_j}$ and $Q_j$ is multiplication by $x_j$. Sometimes we use $x_j$ denote the operator of multiplication by $x_j$. The commutator $ i[P_j,Q_k]=\delta_{jk}$ and $P^2=P_jP_j=-\Delta_x$ where $\delta_{jk}$ is the Kronecker delta. $\{e_1,\cdots, e_n\}$ denotes a basis in $\mathbb{R}^n$. $\tau$ denotes the operator of dilation $(\tau_\delta f)(x)=f(\delta x)$.\par
We also assume $\beta(t\leq 1)=1-\beta(t>1)$ and
\begin{equation}
\sup\limits_{n=0,1,2,3,4} \|\beta^{(n)}(t)\|_{\s^\infty_t}\leq C_\beta,
\end{equation}
with 
\eq
\beta^{(n)}(t):=\frac{d^n}{dt^n}[\beta(t)].
\eeq
\subsection{CL and Improved CL}
We start with the introduction of the \emph {time translated (tT) Potential}, the translation being the flow under the free hamiltonian, the Laplacian:
\begin{equation}
\sK_t(V(x,s)):=e^{itH_0}V(Q,s)e^{-itH_0}\label{skt}.
\end{equation}
Since
\begin{equation}
d/dt(e^{itf(P)}g(Q)e^{-itf(P)})=e^{itf(P)}i[f(P),g(Q)]e^{-itf(P)},\label{wuhan 4}
\end{equation}
we have
\[
e^{itH_0}Qe^{-itH_0}=e^{itP^2}Qe^{-itP^2}=Q+\int_0^t (e^{itP^2}(2P)e^{-itP^2})dt=Q+2tP,
\]
which implies
\begin{equation}
e^{itH_0}e^{i\xi\cdot Q}e^{-itH_0}=e^{i\xi \cdot (Q+2tP)}, ~ \text{for }\xi\in \mathbb{R}^n.\label{0.2}
\end{equation}
Based on $i[P_j,Q_j]=1$, we have
\begin{equation}
[i\xi \cdot Q, i t \xi \cdot P]=\sum\limits_{l,j}-t\xi_j\xi_l[Q_j,P_l]=\sum\limits_{l,j}-it\xi_j\xi_l\delta_{jl}=-it\xi^2.\label{1.3}
\end{equation}
Then since $[i\xi \cdot Q, i t \xi \cdot P]$ is a $c$-number, according to Baker-Campbell-Hausdorff formula, we have
\begin{equation}
e^{i\xi\cdot(Q+2tP)}=e^{i\xi\cdot Q}\cdot e^{2it\xi \cdot P}\cdot e^{-\frac{1}{2}[i\xi\cdot Q,2it\xi \cdot P]}=e^{i\xi\cdot Q}\cdot e^{2it\xi \cdot P}\cdot e^{it\xi^2}.\label{1.4}
\end{equation}
Based on equation \eqref{1.4}, the representation of the tT  potential operator follows
\begin{equation}
\sK_t(V(x,t))=\frac{1}{(2\pi)^{\frac{n}{2}}}\int d^n\xi \hat{V}(\xi,t)e^{i\xi\cdot Q}e^{2it\xi \cdot P}\cdot e^{it\xi^2}.\label{3.1.1.1}
\end{equation}
Hence,the tT potential satisfies:
\begin{equation}
\|\sK_t(V(x,t))\|_{\s_x^p\to \s_x^p}\leq \frac{1}{(2\pi)^{\frac{n}{2}}}\|\hat{V}(\xi,t)\|_{\s^1_\xi}.
\end{equation}
If $\hat{V}(\xi,t)$ happens to be a finite measure in $\xi$ and its total variation is denoted by $m(t)$, and if $\sup\limits_{t\in \mathbb{R}}m(t)<\infty$, we also have
\begin{equation}
\|\sK_t(V(x,t))\|_{\s_x^p\to \s_x^p}\leq \frac{1}{(2\pi)^{\frac{n}{2}}}\sup\limits_{t\in \mathbb{R}}m(t).\label{mt}
\end{equation}
Then a Cancellation Lemma follows from the above:
\begin{lemma}[CL]\label{cancell}If $\hat{V}(\xi,t)$ is a finite measure whose total variation is denoted by $m(t)$ and if $\sup\limits_{t\in \mathbb{R}}m(t)<\infty$, then
\begin{equation}
\sup\limits_{s\in \mathbb{R}}\|\sK_s(V(x,t))\|_{\s_x^p\to \s_x^p}\leq \frac{1}{(2\pi)^{\frac{n}{2}}} \sup\limits_{t\in \mathbb{R}}m(t).\label{5.20.1}
\end{equation}
\end{lemma}
\begin{proposition}\label{exp1} Recall the definition of $\Omega(0,t)$, see \eqref{1.23}.
\begin{equation}
\ln\left(\left\| \Omega(0,t)\right\|_{\s^p_x\to \s^p_x}\right)\leq \int_0^{|t|}du\|\sK_u(V(x,u))\|_{\s_x^p\to \s_x^p}.
\end{equation}
Therefore if $\hat{V}(\xi,t)$ is a finite measure whose total variation is denoted by $m(t)$ and if
\begin{equation}
c(t):=\int_0^tds|m(s)|\lesssim_t1,
\end{equation}
then for $1\leq p\leq \infty$,
\begin{equation}
\ln\left(\left\|\Omega(0,t)\right\|_{\s^p_x\to \s^p_x}\right) \leq  \frac{c(t)}{(2\pi)^{\frac{n}{2}}}
\end{equation}
or
\begin{equation}
\left\|\Omega(0,t)\right\|_{\s^p_x\to \s^p_x}\leq \exp\left( \frac{c(t)}{(2\pi)^{\frac{n}{2}}}\right).
\end{equation}
Similarly, we have
 \begin{equation}
\ln\left(\left\| \Omega(0,t)^*\right\|_{\s^p_x\to \s^p_x}\right) \leq  \frac{c(t)}{(2\pi)^{\frac{n}{2}}}
\end{equation}
or
\begin{equation}
\left\| \Omega(0,t)^*\right\|_{\s^p_x\to \s^p_x}\leq \exp\left( \frac{c(t)}{(2\pi)^{\frac{n}{2}}}\right).
\end{equation}
\end{proposition}
\begin{proof} Since in $n$ dimensions,
\begin{align}
\sK_t(V(x,t))=e^{itH_0}V(x,t)e^{-itH_0}=&\frac{1}{(2\pi)^{\frac{n}{2}}}\int d^n\xi \hat{V}(\xi,t)e^{itH_0}e^{ix\cdot \xi}e^{-itH_0}\\
=&\frac{1}{(2\pi)^{\frac{n}{2}}}\int d^n\xi \hat{V}(\xi,t)e^{iQ\cdot \xi}e^{2it\xi\cdot P}e^{it\xi^2}
\end{align}
where $Q$ denotes the operator of multiplication by $x$, we obtain
\begin{equation}
\| e^{itH_0}V(x,t)e^{-itH_0}\|_{\s^p\to \s^p}\leq \frac{|m(t)|}{(2\pi)^{\frac{n}{2}}}.\label{eq1}
\end{equation}
Now we prove boundedness of $\Omega(0,t) $. For $\Omega(0,t)$, we use Duhamel's formula and iterate it for infinitely many times
\begin{equation}
\Omega(0,t)=\sum\limits_{k=0}^\infty i^kI^{(k)}(t),
\end{equation}
where
\begin{equation}
I^{(k)}(t):=\int_0^t dt_1\int_{t_1}^t dt_2\cdots \int_{t_{k-1}}^t dt_k \sK_{t_1}(V(x,t_1))\cdots \sK_{t_k}(V(x,t_k)), \quad k=0,1,\cdots,\label{Sep.24.2}
\end{equation}
$I^{(0)}(t)$ denotes the identity. Since
\begin{align*}
\| I^{(k)}(t)\|_{\s^p\to \s^p}\leq& \int_0^{|t|} dt_1\int_{t_1}^{|t|}dt_2\cdots \int_{t_{k-1}}^{|t|} dt_k\| \sK_{t_1}( V(x,t_1)) \|_{\s_x^p\to \s_x^p}\cdots \| \sK_{t_k}( V(x,t_k))\|_{\s_x^p\to \s_x^p}\\=&\frac{1}{k!}\left(\int_0^{|t|}ds \| \sK_{s}( V(x,s)) \|_{\s_x^p\to \s_x^p}\right)^k,
\end{align*}
we have
\begin{equation}
\|\Omega(0,t) \|_{\s_x^p\to \s_x^p}\leq \exp\left(\int_0^{|t|} ds \|\sK_s( V(x,s))\|_{\s_x^p\to \s_x^p} \right).
\end{equation}
Using $\hat{V}(\xi,t)\in \s^\infty_t\s^1_\xi$, due to \eqref{eq1}, we get
\begin{equation}
\ln\left(\left\| \Omega(0,t)\right\|_{\s^p_x\to \s^p_x}\right) \leq  \frac{c(t)}{(2\pi)^{\frac{n}{2}}},
\end{equation}
that is,
\begin{equation}
\|\Omega(0,t)\|_{\s^p_x\to \s^p_x} \leq  \exp\left(\frac{c(t)}{(2\pi)^{\frac{n}{2}}}\right).\label{Sep.24.3}
\end{equation}
Similarly, since
\begin{equation}
\Omega(0,t)^*=\sum\limits_{k=0}^\infty i^k\left(I^{(k)}\right)^*(t),
\end{equation}
where
\begin{equation}
\left(I^{(k)}\right)^*(t):=\int_0^t dt_1\int_0^{t_1} dt_2\cdots \int_0^{t_{k-1}} dt_k \sK_{t_1}(V(x,t_1))\cdots \sK_{t_k}(V(x,t_k)), \quad k=0,1,\cdots,
\end{equation}
we have
\begin{align*}
\|\left( I^{(k)}\right)^*(t)\|_{\s^p\to \s^p}\leq& \int_0^{|t|} dt_1\int_0^{t_1}dt_2\cdots \int_0^{t_{k-1}} dt_k\| \sK_{t_1}( V(x,t_1)) \|_{\s_x^p\to \s_x^p}\cdots \| \sK_{t_k}( V(x,t_k)) \|_{\s_x^p\to \s_x^p}\\=&\frac{1}{k!}\left(\int_0^{|t|}ds \| \sK_{s}( V(x,s)) \|_{\s_x^p\to \s_x^p}\right)^k
\end{align*}
and therefore
\begin{equation}
\|\Omega(0,t)^* \|_{\s_x^p\to \s_x^p}\leq \exp\left(\int_0^{|t|} ds \|\sK_s( V(x,s))\|_{\s_x^p\to \s_x^p} \right).
\end{equation}
Using $\hat{V}(\xi,t)\in \s^\infty_t\s^1_\xi$, due to \eqref{eq1}, we get
\begin{equation}
\ln\left(\left\|\Omega(0,t)^*\right\|_{\s^p_x\to \s^p_x}\right) \leq  \frac{c(t)}{(2\pi)^{\frac{n}{2}}},
\end{equation}
that is,
\begin{equation}
\|\Omega(0,t)^*\|_{\s^p_x\to \s^p_x} \leq  \exp\left(\frac{c(t)}{(2\pi)^{\frac{n}{2}}}\right).
\end{equation}
\end{proof}
It implies immediately the global boundedness of $\Omega(0,t)$ for  Schr\"odinger equations with general potentials; for example, quasi-periodic in $x$, on $\s^\infty$ space in one dimension:
\begin{corollary}In one dimension, if $V(x)$ is quasi periodic, (in other word, if $V(x)$ is a finite sum of terms of the form $a\cos(bx)$ or $a\sin(bx)$) and if the initial data is $de^{icx}$ for some $c,d\in \mathbb{R}$, then $\Omega(0,t)\psi(0)$ where $\psi(t)$ is the solution to system
\begin{equation}
i\partial_t \psi(x,t)=(H_0+V(x))\psi(x,t),\label{Sep.24.1}
\end{equation}
exists in $\s^\infty$ and is a sum of sine and cosine terms only, and is bounded for all times.
\end{corollary}
\begin{proof}
Assume
\begin{equation}
V(x)=\sum\limits_{k=0}^N a_k\cos(b_k x)+c_k\sin(d_kx)
\end{equation}
The boundedness follows from \eqref{Sep.24.3} with
\begin{equation}
c(t)\leq t\sum\limits_{k=0}^N|a_k|+|c_k|.
\end{equation}
The solution is a sum of sine and cosine terms only since
\begin{equation}
\sK_{t}(e^{iax})\psi(0)=\sK_{t}(e^{iax})(de^{icx})=de^{ita^2}e^{iax}e^{ic(x+2ta)}=de^{i(ta^2+2tac)}e^{ix(a+c)}.
\end{equation}
\end{proof}
In particular, if both initial data $\psi(x,0)$ and the potential $V(x,t)$ are smooth in $x$, then so is the solution:
\begin{corollary}
If both initial data $\psi(x,0)$ and the potential $V(x,t)$ are smooth in $x$, then so is the solution of \eqref{Sep.24.1}.
\end{corollary}
\begin{proof}
If the initial data $\psi(x,0)$is smooth in $x$, then in \eqref{Sep.24.2}, take $n$th order derivative on both sides and on the right hand side, one can commute through the derivative; it hits the potential term. So if $V(x,t)$ is smooth in $x$, then so is the solution for all times.
\end{proof}
Now we would like to introduce the \emph {Integrated tT Potential operator}
\begin{equation}
I\sK:= \int_0^\infty dt \sK_t(V(x,t))
\end{equation}
which is relevant to the  $\s^p$ boundedness of the wave operator. Based on Cook's method, one can prove the existence of $I\sK: \s_x^1\cap \s_x^2\to \s_x^2$ when $V(x,t)\in \s^\infty_t\s^\infty_x\cap \s^\infty_t\s^2_x $.
\begin{lemma}\label{Iexist}When $V(x,t)\in  \s^\infty_t \s^\infty_x\cap \s^\infty_t\s^2_x(\mathbb{R}\times\mathbb{R}^n)$, $n\geq 3$, $I \sK : \s_x^1\cap \s_x^2\to \s_x^2$ exists and is bounded.
\end{lemma}
\begin{proof}Let $\psi\in \s_x^1\cap \s_x^2$. Since
\begin{equation}
\| e^{itH_0}V(Q,t)e^{-itH_0}\psi\|_{\s^2_x}\lesssim_n \frac{1}{\langle t\rangle^{n/2}}\| V(x,t)\|_{\s^\infty_t\s^2_x\cap\s^\infty_t \s^\infty_x }\|\psi(x)\|_{\s^2_x\cap \s^1_x}
\end{equation}
where we use $e^{itH_0}$ is unitary on $\s^2$ and the decay estimates of $e^{itH_0}$ on $\s^1$, we have
\begin{equation}
\| I\sK\|_{\s_x^2\cap \s_x^1\to \s_x^2}\lesssim_n \| V(x,t)\|_{\s^\infty_t\s^2_x\cap \s^\infty_t\s^\infty_x }\int_0^\infty \frac{dt}{\langle t\rangle^{n/2}} \lesssim_n \| V(x,t)\|_{\s^\infty_t\s^2_x\cap \s^\infty_t\s^\infty_x }.
\end{equation}

\end{proof}
Once we know the existence of $I\sK$ on $\s_x^1\cap \s_x^2$, we can redefine $I\sK$ in Abelian limit sense
\begin{equation}
I\sK=s\text{-}\lim\limits_{\epsilon \downarrow 0}I\sK_\epsilon, \text{ on }\s_x^1\cap \s_x^2
\end{equation}
where
\begin{equation}
I\sK_\epsilon:=\int_0^\infty dt e^{-\epsilon t}\sK_t(V(x,t)).
\end{equation}
There is no confusion about this limit taking in strong sense since due to the same argument in Lemma \ref{Iexist} we have that $I\sK_{\epsilon}: \s_x^1\cap \s_x^2\to \s_x^2$ is uniformly bounded in $\epsilon\in [0,1]$. Based on this definition of $I\sK$, when $V$ is time-independent, we get the following representation of $I\sK$:
\begin{lemma}\label{2re}If $\hat{V}(\xi)\in \s^1_\xi,$ then for $\epsilon>0$,
\begin{equation}
I\sK_\epsilon=\frac{1}{(2\pi)^{n/2}}\int d^3\xi \hat{V}(\xi)e^{ix\cdot\xi}\frac{-1}{i(\xi^2+2\xi\cdot P)-\epsilon}\label{ire}
\end{equation}
\end{lemma}
\begin{proof}It suffices to check on a dense set of $\s_x^1\cap \s_x^2$. Choose $\psi\in \s_x^\infty\cap \s_x^1$. According to the identity \eqref{1.4},
\begin{equation}
I\sK_\epsilon\psi(x)=\frac{1}{(2\pi)^{n/2}}\int_0^\infty dt \int d^n\xi \hat{V}(\xi)e^{ix\cdot \xi}e^{it\xi^2-\epsilon t}\psi(x+2t\xi)\label{I1}.
\end{equation}
That $\psi\in \s_x^\infty$, $e^{-\epsilon t}\in \s^1_t[0,\infty)$ and $\hat{V}(\xi)\in \s^1_\xi$ imply
\begin{equation}
\hat{V}(\xi)e^{ix\cdot \xi}e^{it\xi^2-\epsilon t}\psi(x+2t\xi)\in \s^1_t[0,\infty)\s^1_\xi .
\end{equation}
Then by Fubini's theorem, we change the order of the integral and then take the integral over $t$
\begin{equation}
I\sK_\epsilon\psi=\frac{1}{(2\pi)^{n/2}}\int d^n\xi \hat{V}(\xi)e^{ix\cdot\xi}\frac{-1}{i(\xi^2+2\xi\cdot P)-\epsilon}\psi.
\end{equation}

\end{proof}
$I\sK$ is regarded as the limit of $I\sK_{\epsilon}$ as $\epsilon\downarrow0$ in strong topology. Based on Lemma \ref{2re},
\begin{equation}
I\sK=\frac{1}{(2\pi)^{n/2}}\int d^n\xi \hat{V}(\xi)e^{ix\cdot\xi}\frac{-1}{i(\xi^2+2\xi\cdot P)-0}.
\end{equation}

For the construction of the wave operator, we have to introduce another representation formula for $I\sK_\epsilon$. Choose $V(x,t)\in \mathcal{S}_t\mathcal{S}_x$. For $\psi\in \s_x^\infty \cap \s_x^p$, in identity \eqref{I1}, we use Fubini's theorem to integrate over $t$ first, use spherical coordinates of $\xi$, then change variables from $t\to u=|\xi|t$ and then change the order of the integral over $|\xi|$ and $u$
\begin{equation}
I\sK_\epsilon\psi(x)=\frac{1}{(2\pi)^{\frac{3}{2}}} \int_{S^2}d\sigma(\xi)\int_0^\infty du\int_0^\infty d|\xi||\xi|\hat{V}(\xi,\frac{u}{|\xi|}) e^{-\frac{\epsilon u}{|\xi|}+i(x\cdot \xi +u|\xi|)}\psi(x+2u\hat{\xi}).
\end{equation}
Then for $\psi\in L^p$ and general $V(x,t)$, we have a representation
\begin{equation}
I\sK_\epsilon\psi(x)=\frac{1}{(2\pi)^{\frac{3}{2}}} \int_{S^2}d\sigma(\xi)\int_0^\infty du\int_0^\infty d|\xi||\xi|\hat{V}(\xi, \frac{u}{|\xi|}) e^{-\frac{\epsilon u}{|\xi|}+i(x\cdot \xi +u|\xi|)}\psi(x+2u\hat{\xi}).\label{Ig1}
\end{equation}
\subsection{Improved CL For Time Dependent Potentials}
For the tT Potentials in general, we cannot prove the improved cancellation lemma (ICL) without regularity assumptions in $x$, when the potentials are time-dependent.\emph {To be precise, if we just assume $V(x,t)\in C_t\s^1_x $, the improved cancellation lemma fails}.

  \par Let $B_{\infty,2}(T)(T>0)$ denote the space of bounded linear transformation from $C_t([-T,T])\s^2_x$ to $\s^p_t([-T,T]) \s^2_x$ with $(p>1),$ and its standard norm is denoted by $\| \cdot\|_{B_{\infty,2}(T)}$. Now we consider the following linear transformation
\begin{equation}
\mathcal{L}_T: \mathcal{D}_T\to B_{\infty,2}(T), \quad  V(x,t)\mapsto \sK_t(V(x,t)) \label{Vregular}
\end{equation}
where
\begin{equation}
\mathcal{D}_T:=\{V(x,t)\in C_t([-T,T])L^1_x: \| \mathcal{L}_T( V(x,t))\|_{B_{\infty,2}(T)}<\infty\}.
\end{equation}
The following lemma reveals the unbounded nature of $\mathcal{L}_T$:
\begin{lemma}\label{2iteration}For all $T>0$, $\mathcal{L}_T$ defined in \eqref{Vregular} is unbounded.
\end{lemma}
\begin{proof}Prove by contradiction. Assume there exists $T_0>0$ such that
\begin{equation}
L_{T_0}:=\| \mathcal{L}_{\frac{T_0}{2}}\|_{ \mathcal{D}_{T_0/2}\to B_{\infty,2}(T_0/2) }<\infty\label{assumption}
 \end{equation}
 and therefore $\mathcal{D}_{T_0}=C_t([-T_0,T_0])\s^1_x $. According to the definition of $L_T$, we have 
 \eq
 L_{t_1}\leq L_{t_2}, \text{ if }0<t_1<t_2,
 \eeq
 which implies $L_T<\infty$ if $T\leq T_0$. In the following, we are going to use this to get a contradiction. We consider a NLS system
\begin{equation}
i\partial_t \psi(t)=H_0\psi(t)+|\psi(t)|^{p-1}\psi(t), \text{ with }p=3, n=3.\label{NLSex}
\end{equation}
We  are going to show that if \eqref{assumption} holds, it implies the local wellposedness of this NLS in $L^2_x(\mathbb{R}^n).$ This violates the known result that well-posedness in $H^s_x(\mathbb{R}^n)$ holds, if and only if $s \geq \max(s_c, 0)$, where $s_c := \frac{d}{2}-\frac{2}{p-1}.$

For $\psi(0)= \psi_0\in  L^2_x(\mathbb{R}^n)$, let us consider the following iteration
\begin{equation}
\phi_k(t)=e^{-itH_0}\psi_0+(-i)\int_0^tds e^{-itH_0}\sK_s(|\phi_{k-1}(s)|^2)e^{isH_0}\phi_k(s)
\end{equation}
with $\phi_0=e^{-itH_0}\psi_0.$ Since due to the definition of $\mathcal{L}_T$ and H\"older's inequality,
\begin{equation}
\| \int_0^tds e^{-itH_0}\sK_s(|\phi_{k-1}(s)|^2)f(x,s)\|_{C_t([-T,T])\s^2_x}\leq T^{p'} \s_T\| |\phi_{k-1}(t)|^2\|_{C_t([-T,T])\s^1_x}\| f(x,t)\|_{C_t([-T,T])\s^2_x},\label{3.12.1}
\end{equation}
due to Corollary \ref{exp1},
\begin{equation}
\| \phi_k(t)\|_{C_t([-T,T])\s^2_x}\leq \|\phi_0\|_{\s^2_x}\exp \left(T^{p'} \| |\phi_{k-1}(t)|^2\|_{C_t([-T,T])\s^1_x}\s_T \right),
\end{equation}
if $\phi_{k-1}(t)\in C_t([-T,T])\s^2_x $. Since $\phi_0=e^{-itH_0}\psi_0\in C_t([-T,T])\s^2_x$, due to conservation law, we have
\begin{equation}
\| \phi_k(t)\|_{\s^2_x}=\| \psi_0\|_{\s^2_x}, \text{ for all }k=0,\cdots.
\end{equation}
Since
\begin{align*}
&\sK_t(|\phi_{k-1}|^2)e^{itH_0}\phi_k-\sK_t(|\phi_{k}|^2)e^{itH_0}\phi_{k+1}\\
=&\sK_t((\phi_{k-1}-\phi_k)^*\phi_{k-1}) e^{itH_0}\phi_k+\sK_t(\phi_k^*(\phi_{k-1}-\phi_k))e^{itH_0}\phi_k+\sK_t(|\phi_k|^2)e^{itH_0}(\phi_k-\phi_{k+1}),
\end{align*}
applying estimate \eqref{3.12.1}, we get
\begin{align*}
&\| \phi_k(t)-\phi_{k+1}(t)\|_{C_t([-T,T])\s^2_x}\\\leq &2T^{p'}\s_T\|\psi_0\|_{\s^2_x}^2\|\phi_k(t)-\phi_{k-1}(t) \|_{C_t([-T,T])\s^2_x}+T^{p'}\s_T\|\psi_0\|_{\s^2_x}^2\|\phi_k(t)-\phi_{k+1}(t) \|_{C_t([-T,T])\s^2_x},
\end{align*}
which implies
\begin{align}
&\| \phi_k(t)-\phi_{k+1}(t)\|_{C_t([-T,T])\s^2_x}\\
\leq& 2T^{p'}\s_T\|\psi_0\|_{\s^2_x}^2\|\phi_k(t)-\phi_{k-1}(t) \|_{C_t([-T,T])\s^2_x}+T^{p'}\s_T\|\psi_0\|_{\s^2_x}^2\|\phi_k(t)-\phi_{k+1}(t) \|_{C_t([-T,T])\s^2_x}.
\end{align}
Choose $T$ small enough such that $ T^{p'}\s_T\|\psi_0\|_{\s^2_x}^2\leq \frac{1}{8}$. Then
\begin{equation}
\| \phi_k(t)-\phi_{k+1}(t)\|_{C_t([-T,T])\s^2_x}\leq \frac{1}{2}\|\phi_k(t)-\phi_{k-1}(t) \|_{C_t([-T,T])\s^2_x}.\label{L2}
\end{equation}
By contraction mapping principle, we get local wellposedness in $\s_x^2$. Then based on the same argument, we get global existence of \eqref{NLSex}. Contradiction since in \cite{MRS2014}, Merle, Rapha\"el and Szeftel showed there is a solution $u\in C_t([0,T))H^{1/2}_x\subseteq \s^2_x$ which blows up in $\s^2_x$ at time $T$. Also, in \cite{CCT2003}, Christ, Colliander and Tao sketched the proof of the ill-posedness in $\s^2_x$.

\end{proof}
\begin{remark}Lemma \ref{2iteration} implies the failure of local smoothing property for some $C_t\s^1_x$ localization. In other word, for some $V(x,t)\in C_t\s^1_x$, any $A>0$, the map $\mathcal{C}: C_t([-A,A])\s^2_x\to \s^1_t([-A,A])\s^2_x, f\mapsto V(x,t)e^{-itH_0}f$, is unbounded.

\end{remark}

By applying a similar argument we get as an application, useful for decay estimates for rough potentials the following:
\begin{lemma}If $V(x,t)\in \s^\infty_t\s^q_x(\mathbb{R}^3)$ for $q\in (\frac{4}{3}, 2]$, then for $t\in (0,1], s\in[\frac{t}{2},t)$, $\sK_s e^{itH_0}:\s_x^1\cap \s_x^2(\mathbb{R}^3)\to \s_x^\infty$ is bounded with
\begin{equation}
\|\sK_s e^{itH_0} \|_{\s_x^1\cap \s_x^2\to \s_x^\infty}\lesssim \frac{1}{t^{3/2}}\times \frac{1}{(t-s)^{1-\epsilon}}
\end{equation}
for some $\epsilon=\epsilon(q)\in (0,1]$.
\end{lemma}

\begin{proof}Let $\psi\in \mathcal{S}$ and $\hat{V}(\xi,t)\in \s^\infty_t\s^1_\xi $.  According to the same computation above,
\begin{equation}
\sK_se^{itH_0}\psi=\frac{1}{(2\pi)^{3/2}}\int d^3\xi \hat{V}(\xi,s)e^{ix\cdot \xi}e^{is|\xi+P|^2}e^{i(t-s)P^2}\psi.
\end{equation}
Let $\psi_{t-s}:=e^{i(t-s)P^2}\psi$. Then $\psi_{t-s}\in \s_x^2\cap \s_x^\infty$ when $s<t$.
\begin{equation}
e^{is|\xi+P|^2}e^{i(t-s)P^2}\psi=\frac{1}{(2\pi is)^{3/2}}\int d^3k e^{-i\frac{k^2}{2s}}\psi_{t-s}(x-k)e^{-ix\cdot \xi}e^{i(x-k)\cdot \xi}.
\end{equation}
Hence,
\begin{equation}
\sK_se^{itH_0}\psi=\frac{1}{(2\pi)^{3/2}}\times\frac{1}{(2\pi is)^{3/2}}\int d^3\xi d^3k\hat{V}(\xi,s)e^{i(x-k)\cdot \xi}e^{-i\frac{k^2}{2s}}\psi_{t-s}(x-k).\label{A.6.2}
\end{equation}
$\psi\in \mathcal{S}$ implies $\psi_{t-s}(x)\in \s^1_x $. Then $\hat{V}(\xi,s)\psi_{t-s}(x-k)\in \s^1_\xi \s^1_k$. By Fubini's theorem, we change the order of the integral and integrate over $\xi$ first
\begin{equation}
\sK_se^{itH_0}\psi=\frac{1}{(2\pi is)^{3/2}}\int d^3ke^{-i\frac{k^2}{2s}}\psi_{t-s}(x-k)V(x-k,s).
\end{equation}
Then when $V(x,t)\in \s^\infty_t\s^q_x$ for $q\in (\frac{4}{3}, 2]$, by H\"older's inequality,
\begin{equation}
\|\psi_{t-s}(x-k)V(x-k,s)\|_{\s^1_k}\leq \| \psi_{t-s}(x-k)\|_{\s^{q^\prime}_k}\| V(x-k,t) \|_{\s^\infty_t\s^q_k}\lesssim \frac{\| V(x-k,t) \|_{\s^\infty_t\s^q_k}\|\psi\|_{\s_x^1\cap \s_x^2}}{(t-s)^{3(2-q^\prime)/2}}.\label{A.6.1}
\end{equation}
$q\in (\frac{4}{3}, 2]$ implies $ 3(2-q^\prime)/2\in [0,1)$. Then we use Banach Space Continuity twice and get the same inequality \eqref{A.6.1} for $\psi\in \s^1_x\cap \s^2_x(\mathbb{R}^3)$, $V\in \s^\infty_t\s^q_x(\mathbb{R}^3)$. Combining this inequality with \eqref{A.6.2}, we complete the proof.
\end{proof}

For the construction of the wave operator, we also need to introduce the following operators
\begin{equation}
I_\epsilon^{(k)}:=\int_0^\infty dt_1\int_0^{t_1}dt_2\cdots \int_0^{t_{k-1}}dt_ke^{-\epsilon t_1}\sK_{t_k}(V(x,t_k))\cdots\sK_{t_1}(V(x,t_1)) ,\text{ for }k=1,2,\cdots.
\end{equation}

\section{Time-independent potentials in \texorpdfstring{$\mathbb{R}^3$}{R3} }\label{section 2}
In this section, we prove the $\s^p$ boundedness of the wave operator $\Omega$ on the high-frequencypart of the $\s^p$ space for time-independent potentials $V(x),$ on $\s^p$ space in $\mathbb{R}^3$.  We assume
 \begin{equation}
\begin{cases} K_m(V(x))= \sup\limits_{\eta\in \mathbb{R}^3}|K_1(V(x),\eta)|<\infty, \\
 \hat{V}_a(\xi)\in \s^1_\xi.
 \end{cases}
 \end{equation}
 Recall that $L_{\eta,l,j}(k,\hat{\xi},\epsilon)$ denotes the Fourier transform of $|\xi|\partial_{\xi\cdot e_l}^j[\hat{V}(\xi-\eta)]e^{- \frac{\epsilon}{|\xi|}}$ in $|\xi|$ variable for $l=1,2,3,$ $j=0,1,2$,
\begin{equation}
K_1(V(x),\eta)=\max\limits_{l=1,2,3, j=0,1,2}\int_{S^2}d\sigma(\xi)\int_{-\infty}^\infty dk \sup\limits_{\epsilon \geq0}|L_{\eta,l,j}(k,\hat{\xi},\epsilon)|
\end{equation}
 and
 \begin{equation}
 \hat{V}_a(\xi)=\sum\limits_{j,l=0}^2\sum\limits_{r,m=1}^3|\partial_{\xi\cdot e_r}^j\partial_{\xi\cdot e_m}^l\hat{V}(\xi)|, \text{ with a basis }\{e_1,e_2,e_3\}.
 \end{equation}
 We begin with some basic lemmas.
 \subsection{Some lemmas}\label{sectionbasic}
For the $\s^p$ estimates for $I\sK$ and wave operator in the following context, we need some lemmas:
\begin{lemma}\label{c1} Let $f(u)\in \s^1_u(\mathbb{R})$. The operator $\mathcal{T}_{\hat{\xi}}: \s^p(\mathbb{R}^n)\to \s^p(\mathbb{R}^n)$
\begin{equation}
\mathcal{T}_{\hat{\xi}}(\psi)(x):=\int_0^\infty dk f(x\cdot\hat{\xi}+k)\psi(x+2k\hat{\xi})
\end{equation}
is uniformly bounded in $\hat{\xi}\in S^{n-1}$ for $1\leq p\leq \infty$ with upper bound $\| f(k)\|_{\s^1_k}$.
\end{lemma}
\begin{proof}Write $x:=\sum\limits_{j=1}^n x_je_j=(x_1,\cdots,x_n)$ with $e_1:=\hat{\xi}$. We do a change of variables $k\to u=k+x\cdot \hat{\xi}$
\begin{equation}
\mathcal{T}_{\hat{\xi}}(\psi)(x)=\int^\infty_{x\cdot\hat{\xi}}du f(u)\psi(2u-x_1,x_2,\cdots,x_n).
\end{equation}
Then by Minkowski's integral inequality,
\begin{equation}
\|\mathcal{T}_{\hat{\xi}}(\psi)(x) \|_{\s^p_x}\leq \int  |f(u)| \| \psi(2u-x_1,x_2,\cdots,x_n )\|_{\s^p_x}du= \| f(u)\|_{\s^1_u}\|\psi(x)\|_{\s^p_x}.
\end{equation}
\end{proof}
\begin{lemma}\label{g}
For $d\in \{1,2,3,4\} $, $j\in \{0,1,2\}$, $M>1$, $\epsilon\in \mathbb{R}$, $1\leq p\leq \infty$, let
\begin{equation}
\mathscr{P}_{jd}(M,\epsilon):=\frac{\beta^{(j)}(| P|>2M)}{(P+i\epsilon)^d }: \s^p(\mathbb{R})\to \s^p(\mathbb{R}),
\end{equation}
be a Fourier multiplier. Then $\|\mathscr{P}_{jd}(M,\epsilon)\|_{\s^p\to \s^p}\lesssim \frac{1}{M^d}$. In addition, for $\psi\in \s^p$,
\begin{equation}
\|\sup\limits_{\epsilon\in [0,1]}| \mathscr{P}_{jd}(M,\epsilon)\psi(x)|\|_{\s^p_x}\lesssim \frac{1}{M^d}\|\psi(x)\|_{\s^p_x}.
\end{equation}
\end{lemma}
\begin{proof}When $d=1$, it suffices to show that it is the Fourier transform of some finite Borel measure $\mu_M$ whose total variation is less than $C/M$. Let
\begin{equation}
\mu(x):=\mathscr{F}^{-1}_q[\frac{\beta(|q|>2)}{q+i\epsilon/M}](x), \text{ and then }\mu_M(x)=[\mathscr{F}^{-1}\tau_{1/M}\mathscr{F}[\frac{\mu}{M}]](x)=[\tau_{M}\mu](x)=\mu(Mx)
\end{equation}
since $\mathscr{F}\sigma_\delta=|\delta|^{-1}\tau_{\delta^{-1}}\mathscr{F}$.
We are going to show $M\int |d\mu_M(x)|\lesssim 1$ for $d=1,$ and the other cases will follow by the same way. Since for $q$ large, $\frac{1}{q+i\epsilon/M}\sim \frac{1}{q}$, then for $|x|\leq 1$, $|d\mu(x)|\lesssim -\ln|x| dx$. For $|x|>1$, since $|d\mu(x)|\lesssim_N \frac{1}{|x|^N}dx$ for any $N\geq 1$, by the use of integration by parts, then $|\mu(x)|\lesssim \frac{1}{x^2}$. Hence,
\begin{equation}
\int |d\mu_M(x)|=\frac{1}{M}\int M|d\mu(Mx)|=\frac{1}{M}\int |d\mu(x)|\lesssim\frac{1}{M}.
\end{equation}

\end{proof}
In \cite{JSS1991}, Journ\'e, Soffer and Sogge proved that the high energy cutoff function $\gamma(H/M): \s^1(\mathbb{R}^n)\to \s^1(\mathbb{R}^n)$ is bounded for each $M>0$, when $\gamma\in C^\infty(\mathbb{R})$ satisfying $\gamma(\lambda)=1$ for $\lambda\geq 1,$ and $\beta(\lambda)=0$ for $-\infty <\lambda <1/2$; $H=H_0+V(x)$ for some nice $V(x)$ including the case when $H=H_0$.

 When $H=H_0$, this high energy function $\gamma(H_0>M)$ is Fourier multiplier, and it implies that $\beta(|P|>M)$ is also bounded on $\s^1$ by taking $\gamma(H_0/M^2)=\beta(\sqrt{H_0/M^2})$. By duality, we get the $ \s^p$ boundedness of $\beta(|P|>M)$ for all $1\leq p\leq\infty$. We will use the $\s^p$ boundedness of $\beta(|P|>M)$ throughout the following context. Let
\begin{equation}
E_{n,M}:=\max\left(\|\beta(|P|>M) \|_{\s_x^p(\mathbb{R}^n)\to \s_x^p(\mathbb{R}^n)}, \|\beta(|P|\leq M) \|_{\s_x^p(\mathbb{R}^n)\to \s_x^p(\mathbb{R}^n)}\right)\label{gg2}
\end{equation}
in dimension $n$.
\begin{lemma}\label{power}If $\mathscr{T}(\eta): \s^p(\mathbb{R}^n) \to \s^p(\mathbb{R}^n)$, is bounded with
\begin{equation}
A:=\sup\limits_{\eta\in \mathbb{R}^n} \|\mathscr{T}(\eta) \|_{\s^p(\mathbb{R}^n)\to \s^p(\mathbb{R}^n) },
\end{equation}
then for $f(\xi)\in \s^1_\xi (\mathbb{R}^n)$, we have
\begin{equation}
\left\|\int d^n\xi_1\cdots d^n\xi_n f(\xi_1)f(\xi_2-\xi_1)\cdots f(\xi_k-\xi_{k-1}) \mathscr{T}(\xi_k)\right\|_{\s^p\to \s^p}\leq A\| f(\xi)\|_{\s^1_\xi}^k.
\end{equation}
\end{lemma}
\begin{proof}
It follows from
\begin{align*}
&\left\|\int d^n\xi_1\cdots d^n\xi_n f(\xi_1)f(\xi_2-\xi_1)\cdots f(\xi_k-\xi_{k-1}) \mathscr{T}(\xi_k)\right\|_{\s^p\to \s^p}\\
\leq& \int d^n\xi_1\cdots d^n\xi_n |f(\xi_1)f(\xi_2-\xi_1)\cdots f(\xi_k-\xi_{k-1})| \sup\limits_{\eta\in \mathbb{R}^n}\|\mathscr{T}(\eta)\|_{\s^p\to \s^p}\\
=& A\| f(\xi)\|_{\s^1_\xi}^k.
\end{align*}
\end{proof}
\subsection{\texorpdfstring{$\s^p$}{\sp} boundedness for \texorpdfstring{$I^{(*)}$}{I(*)}}
Let
\begin{equation}
I^{(*)}\psi(x)=\sup\limits_{\epsilon\geq 0}| I_\epsilon\psi(x)|, \text{ for }\psi\in \s^p.
\end{equation}
\begin{theorem}\label{320I1}
If $K_1(V(x),0)<\infty$, then for $1\leq p\leq \infty$, $\psi\in \s^p$,
\begin{equation}
\| I^{(*)}\psi(x) \|_{\s^p_x}\lesssim K_1(V(x),0)\|\psi(x)\|_{\s^p_x}.
\end{equation}
\end{theorem}
\begin{proof}
According to equation \eqref{Ig1},
\begin{equation}
I_\epsilon\psi(x)=\frac{1}{(2\pi)^{\frac{3}{2}}} \int_{S^2}d\sigma(\xi)\int_0^\infty du\int_0^\infty d|\xi||\xi|\hat{V}(\xi) e^{-\frac{\epsilon u}{|\xi|}+i(x\cdot \xi +u|\xi|)}\psi(x+2u\hat{\xi}).
\end{equation}
Then
\begin{equation}
I^{(*)}\psi(x)\leq \int_{S^2}d\sigma(\xi)\int_{-\infty}^\infty du\left(\sup\limits_{\epsilon\geq 0}|L_{0,1,0}(x\cdot\hat{\xi}+u,\hat{\xi},\epsilon)|\right)|\psi(x+2u\hat{\xi})|
\end{equation}
where we use
\begin{equation}
L_{0,1,0}(k,\hat{\xi},u\epsilon)=\mathscr{F}_{|\xi|}(\chi(|\xi|\geq 0)|\xi|\hat{V}(\xi)e^{-\frac{u\epsilon}{|\xi|}})
\end{equation}
and
\begin{equation}
\sup\limits_{\epsilon \geq 0}| L_{0,1,0}(k,\hat{\xi},u\epsilon)|=\sup\limits_{\epsilon \geq 0}| L_{0,1,0}(k,\hat{\xi},\epsilon)|, \text{ for }u>0.
\end{equation}
Due to Lemma \ref{c1}, we have
\begin{equation}
\| I^{(*)}\psi(x) \|_{\s^p_x}\lesssim K_1(V(x),0)\|\psi(x)\|_{\s^p_x}.
\end{equation}
\end{proof}
Recall that
\begin{equation}
\sK_t(V(x,s))=e^{itH_0}V(Q,s)e^{-itH_0}.
\end{equation}
To Proceed, we need more general operators
\begin{equation}
T_\epsilon(\eta):=\int_0^\infty dt e^{-\epsilon t}\sK_t(V(x)e^{i\eta\cdot x}),\label{Tetae}
\end{equation}
and
\begin{equation}
\partial_{\eta\cdot e_j}^l[T_\epsilon(\eta)]:=\int_0^\infty dt e^{-\epsilon t}\sK_t((ix\cdot e_j)^l V(x)e^{i\eta\cdot x}), \text{ for }\epsilon\geq0.
\end{equation}
The corresponding maximal $T$ transform  is
\begin{equation}
[T_{j,l}(\eta)]^{(*)}\psi(x)=\sup\limits_{\epsilon\geq 0}|\partial_{\eta\cdot e_j}^l[T_\epsilon(\eta)]\psi(x)|.
\end{equation}
\begin{corollary}\label{T1}If $V(x)$ satisfies condition \eqref{intimec}, then
\begin{equation}
\|[T_{j,l}(\eta)]^{(*)}\psi(x)\|_{\s^p_x}\lesssim K_m\|\psi(x)\|_{\s^p_x},\quad j=1,2,3, ~l=0,1,2.
\end{equation}
\end{corollary}
\begin{proof} Repeating the proof of Theorem \ref{320I1} by replacing $\hat{V}(\xi)$ with $\partial_{\eta\cdot e_j}^l[\hat{V}(\xi-\eta)],$  and taking the supremum over $\eta\in \mathbb{R}^3,$  we will get the same an upper bound, with $K_m$ instead of $K_1$.
\end{proof}
\subsection{\texorpdfstring{$\s^p$}{\sp} boundedness of \texorpdfstring{$I_M^{(*,k)}$}{IM(*,k)}}
Let
\begin{equation}
I_M^{(*,k)}\psi(x):=\sup\limits_{\epsilon\geq 0}| I_\epsilon^{(k)}\beta(|P|>M)\psi(x) |, \text{ for }\psi\in \s_x^p.
\end{equation}
Before controlling the $\s_x^p$ norm of $I_M^{(*,k)}\psi(x)$, we introduce the following expression:
\begin{lemma}[Representation formula 1]\label{refor}For $\xi_i\in \mathbb{R}^n$, $i=1,\cdots, k$ $(k \in\mathbb{N}^+)$, $\psi(x)\in \s_x^p(\mathbb{R}^n)$, we have
\begin{equation}
\sK_{t_k}(e^{i(\xi_k-\xi_{k-1})\cdot x})\cdots \sK_{t_1}(e^{i(\xi_1-\xi_{0})\cdot x})\psi=\left[e^{i(Q\cdot \xi_k)}e^{it_k(\xi_k^2+2\xi_k\cdot P)}\Pi_{j=1}^{k-1} e^{i(t_j-t_{j+1})(\xi_j^2+2\xi_j\cdot P)}\right]\psi
\end{equation}
with $\xi_0=0\in \mathbb{R}^n$.
\end{lemma}
\begin{proof}
We prove by induction. When $k=1$,it follows from  equations \eqref{0.2} and \eqref{1.4}. Assume that when $k=m$, the representation formula holds. When $k=m+1$,
\begin{align*}
&\sK_{t_{m+1}}(e^{i(\xi_{m+1}-\xi_{m})\cdot x})\cdots \sK_{t_1}(e^{i(\xi_1-\xi_{0})\cdot x})\psi\\
=&\sK_{t_{m+1}}(e^{i(\xi_{m+1}-\xi_{m})\cdot x })\left[e^{i(Q\cdot \xi_m)}e^{it_m(\xi_m^2+2\xi_m\cdot P)}\Pi_{j=1}^{m-1} e^{i(t_j-t_{j+1})(\xi_j^2+2\xi_j\cdot P)}\right]\psi \\
=&e^{iQ\cdot(\xi_{m+1}-\xi_m)}e^{it_{m+1}[(\xi_{m+1}-\xi_m)^2+2(\xi_{m+1}-\xi_m)\cdot P]}\left[e^{i(Q\cdot \xi_m)}e^{it_m(\xi_m^2+2\xi_m\cdot P)}\Pi_{j=1}^{m-1} e^{i(t_j-t_{j+1})(\xi_j^2+2\xi_j\cdot P)}\right]\psi\\
=&\left[e^{i(Q\cdot \xi_{m+1})}e^{it_{m+1}(\xi_{m+1}^2+2\xi_{m+1}\cdot P)}\Pi_{j=1}^{m} e^{i(t_j-t_{j+1})(\xi_j^2+2\xi_j\cdot P)}\right]\psi.
\end{align*}

By induction, we finish the proof.
\end{proof}
Choose $V(x)\in \mathcal{S}_x$. For $\xi=(\xi_1,\cdots,\xi_k)\in \mathbb{R}^{3k}$, let
\begin{equation}
\mathcal{V}(\xi,k):=\frac{1}{(2\pi)^{\frac{3k}{2}}}\hat{V}(\xi_1)\hat{V}(\xi_2-\xi_1)\cdots\hat{V}(\xi_k-\xi_{k-1}).
\end{equation}
Writing $\sK_{t_j}(V(x))$ as
$$
\sK_{t_j}(V(x))=\frac{1}{(2\pi)^{3/2}}\int d^3\xi_j\hat{V}(\xi_j-\xi_{j-1}) \sK_{t_j}(e^{ix\cdot(\xi_j-\xi_{j-1})}), \text{ for }j=1,\cdots,k,
$$
and applying Lemma \ref{refor}, we have
\begin{align*}
I_\epsilon^{(k)}\psi(x)=&\int_0^\infty dt_k\int_{t_{k}}^\infty dt_{k-1}\cdots\int_{t_2}^\infty e^{-\epsilon t_1}dt_1\int
d^3\xi_1\cdots d^3\xi_k \mathcal{V}(\xi,k)\\
&\int d^3qe^{i(x\cdot (\xi_k+q)+t_k(\xi_k^2+2q\cdot \xi_k)+(t_{k-1}-t_k)(\xi_{k-1}^2+2q\cdot \xi_{k-1})+\cdots +(t_1-t_2)(\xi_1^2+2\xi_1\cdot q))}\frac{\hat{\psi}(q)}{(2\pi)^{\frac{3}{2}}}
\end{align*}
where $\frac{1}{(2\pi)^{\frac{3}{2}}}$ comes from the inverse Fourier transform in $q$. It is sufficient to work with $\psi\in \beta(|P|>32M)\mathcal{S}_x, V(x,t)\in \mathcal{S}_t\mathcal{S}_x$ to get concise representation of $I_\epsilon^{(k)}\psi(x).$ This can then be extended to all of $\s_x^p$ and general $V$ by the density continuation. For any $\epsilon>0$,
\begin{equation}
\int_0^\infty dt_k\cdots\int_{t_{2}}^\infty dt_1 \int d^3\xi_1\cdots d^3\xi_kd^3qe^{-\epsilon t_1} |\mathcal{V}(\xi,k)||\hat{\psi}(q)|\leq \frac{1}{(2\pi)^{3k/2}\epsilon^k}\|\hat{V}(\xi)\|_{\s^1_\xi}^k \|\hat{\psi}(q)\|_{\s^1_q}<\infty.\label{F2.26.1}
\end{equation}
Due to Fubini's theorem, we can change the order of the integral over $\xi_j, t_j$ and $q$ when  needed.  We change variables from $t_k$, to $t_k=s_k$, $t_j$, to $t_j=s_k+\cdots+s_j$, $j=1,\cdots, k-1$ with Jacobian $1$,
\begin{align*}
I_\epsilon^{(k)}\psi(x)=&\int_0^\infty e^{-\epsilon s_k}ds_k\int_0^\infty e^{-\epsilon s_{k-1}}ds_{k-1}\cdots\int_{0}^\infty e^{-\epsilon s_1}ds_1\int
d^3\xi_1\cdots d^3\xi_kd^3q \mathcal{V}(\xi,k)\\
&e^{i(x\cdot (\xi_k+q)+(s_k\xi_k^2+\cdots +s_1\xi_1^2)+2(s_k\xi_k+\cdots+s_1\xi_1)\cdot q )}\frac{\hat{\psi}(q)}{(2\pi)^{\frac{3}{2}}}.
\end{align*}
The $\s_x^p$ estimates of $I_\epsilon^{(k)}$ is based on bounding following operator
\begin{equation}
\mathcal{J}_\epsilon(\xi):=\int_0^\infty ds e^{-\epsilon s+i(s|\xi|^2+2s\xi\cdot P)}, \quad\mathcal{J}_\epsilon^{(k)}(\xi_1,\cdots,\xi_k):=\Pi_{j=1}^k \mathcal{J}_{\epsilon}(\xi_j),\text{ for }\xi_j\in \mathbb{R}^3.
\end{equation}
Now we have to recall the definition of the operator $T_\epsilon(\eta)$(see equation \eqref{Tetae})and then rewrite $I_\epsilon^{(k+1)}$ as
\begin{equation}
I_\epsilon^{(k+1)}=\int d^3\xi_1\cdots d^3\xi_{k}\mathcal{V}(\xi,k)T_{\epsilon}(\xi_{k})\mathcal{J}_\epsilon^{(k)}(\xi), \text{ for }\xi=(\xi_1,\cdots,\xi_k).
\end{equation}
We have the following representation and estimates  for $ \int d^3\xi f(\xi) \mathcal{J}_\epsilon(\xi)$.
\begin{lemma}\label{sJ1}Assume $f(\xi)\in C_\xi^2(\mathbb{R}^3)$. For $1\leq p< \infty$, $\int d^3\xi f(\xi)\mathcal{ J}_\epsilon(\xi): \beta(|P|>32M)\s_x^p\to \s_x^p$ and $\int d^3\xi f(\xi)\mathcal{ J}_\epsilon(\xi): \beta(|P|>32M)C_0\to C_0$; preserves the support of the frequency and for $\psi$ in the given space,
\begin{equation}
\int d^3\xi f(\xi)\mathcal{J}_\epsilon(\xi)\psi(x)=\int d^3\xi f(\xi)Q_0\psi(x)+\sum\limits_{j=1}^3\sum\limits_{l=0}^2 \int d^3\xi \partial_{\xi\cdot e_j}^l[f(\xi)]Q_{3(j-1)+l+1}\psi(x)
\end{equation}
for some operator $Q_j=Q_j(\xi,\epsilon): \s_x^q\to \s_x^q$, $1\leq q\leq \infty$, with $\|Q_j(\xi,\epsilon)\beta(|P|>32M)\|_{\s_x^q\to \s_x^q}\lesssim \frac{1}{M}.$ Moreover, for $\psi\in \s_x^p$,
\begin{equation}
\| Q_l^{(*)}(\xi)\psi(x)\|_{\s^p_x}:=\|\sup\limits_{\epsilon\geq 0}| Q_l(\xi,\epsilon)\psi(x)|\|_{\s^p_x}\lesssim \frac{1}{M}\|\psi(x)\|_{\s^p_x}.
\end{equation}
\end{lemma}
\begin{remark}
Here we regard $f(\xi)$ as a multiplier.
\end{remark}
\begin{proof}Since $\mathcal{J}_\epsilon$ is a composition of translation operators, $\int d^3\xi f(\xi)\mathcal{ J}_\epsilon(\xi)$ preserves the support of the frequency. Now we would like to get a  detailed formula. We choose $\psi\in \beta(|P|>32M)\mathcal{S}_x$. According to the similar transformation used  for $I_\epsilon^{(k)}\psi$, we have
\begin{equation}
\int d^3\xi f(\xi)\mathcal{ J}_\epsilon(\xi)\psi=\frac{1}{(2\pi)^{\frac{3}{2}}}\int d^3\xi d^3q\int_0^\infty ds f(\xi) e^{-\epsilon s+i(x\cdot q+s|\xi|^2+2s q\cdot \xi)}\hat{\psi}(q).
\end{equation}
Recall that $\{e_1,e_2,e_3\}$ is a basis in $\mathbb{R}^3$. Let $\xi_j=\xi\cdot e_j$. We claim that for all $\xi\neq 0,$
\begin{align*}
\beta(|P|>32M)=&\left[\sum\limits_{j=1}^3\beta(|\xi_j+P_j|>2M)\beta_j(\xi,P,2M)+\beta(||\xi|+2P\cdot\hat{\xi}|>2M)\gamma(\xi,P,2M)\right]\times\\
&\beta(|P|>32M)=:\beta_{1,1}+\beta_{1,2}+\beta_{1,3}+\beta_{1,4}
\end{align*}
where
\begin{equation}
\beta_j( \xi,P,2M):=\Pi_{l=1}^{j-1}\beta(|\xi_l+P_l|\leq 2M), \text{ for }j=1,2,3, \text{ with }\Pi_{l=1}^0=1
\end{equation}
\begin{equation}
\gamma(\xi,P,2M):=\Pi_{j=1}^3 \beta( |\xi_j+P_j|\leq2M).
\end{equation}
We prove the claim first.
\subsubsection{Proof of the claim}
\begin{proof}
In fact,
\begin{align*}
1= &\sum\limits_{j=1}^3\beta(|\xi_j+P_j|>2M)\beta_j(\xi,P,2M)\\
&+\beta(||\xi|+2P\cdot\hat{\xi}|>2M)\gamma(\xi,P,2M)+\beta(||\xi|+2P\cdot\hat{\xi}|\leq2M)\gamma(\xi,P,2M).
\end{align*}
Then in order to prove the claim, since for $q\in\mathbb{R}^3$, $\beta(|q|>32M)$ implies $|q|>16M$, it suffices to prove that
\begin{equation}
\{q : |q|>16M, |q_j+\xi_j|\leq 2M, j=1,\cdots,3 \}\bigcap \{ ||\xi|+2q\cdot \hat{\xi} |\leq 2M\}=\emptyset.\label{claimid}
\end{equation}
Assume that $|\xi_j+q_j|\leq 2M$, $|q|>16M$. Then
\begin{equation}
|\xi+q|\leq\sqrt{\sum\limits_{j=1}^3(\xi_j+q_j)^2} \leq 2\sqrt{3}M<4M<\frac{|q|}{4},\label{Ni.1}
\end{equation}
which implies
\begin{equation}
|\xi|\geq |q|-|\xi+q|>\frac{3|q|}{4}, \text{ and}\quad |\xi|\leq |q|+|\xi+q|<\frac{7|q|}{4}.\label{Ni.2}
\end{equation}
Then according to equation \eqref{Ni.1}, \eqref{Ni.2},
\begin{equation}
|\xi^2+2\xi\cdot q|=|(\xi+q)^2-q^2|\geq \frac{15q^2}{16}>\frac{15|\xi||q|}{28}>\frac{ 60|\xi|M}{7}>2|\xi|M.
\end{equation}
Hence,
\begin{equation}
| |\xi|+2q\cdot \hat{\xi}|>2M
\end{equation}
which yields equation \eqref{claimid}. Then when multiplied by $\beta(|q|>32M)$, $\beta(||\xi|+2q\cdot\hat{\xi}|\leq2M)\gamma(\xi,q,2M)$ drops and therefore  the claim follows.
\end{proof}
So $\psi$ can be written as a sum of $4$ parts:
\begin{equation}
\psi=\beta_{1,1}\psi+\beta_{1,2}\psi+\beta_{1,3}\psi+\beta_{1,4}\psi=:\psi_1+\psi_2+\psi_3+\psi_4.\label{psij1}
\end{equation}
For $\psi_j$, $j=1,2,3$,
\begin{equation}
\psi_j(x)=\beta(|\xi_j+P_j|>2M)\beta_j(\xi,P,2M)\psi(x)=:\beta(|\xi_j+P_j|>2M)\psi_{j,1}(x)\label{psij}.
\end{equation}
Recalling the definition of $E_{n,M}$ in equation \eqref{gg2},
\begin{equation}
\|\psi_{j,1}(x)\|_{\s^p_x}\leq E_{3,2M}^{j-1}\|\psi(x)\|_{\s^p_x}, \text{ and }\|\psi_j(x)\|_{\s^p_x}\leq E_{3,2M}^j\|\psi(x)\|_{\s^p_x}.
\end{equation}
Since $\beta(|\xi_j+P_j|>2M)$ implies $|\xi_j+q_j|>M$($q$ denotes the argument of $\hat{\psi}$), for $s\geq \frac{1}{M}$ we do integration by parts in $\xi_j$ twice, by setting
\begin{equation}
e^{is(\xi_j^2+2\xi_j q_j)}=\frac{1}{2is(\xi_j+q_j)}\partial_{\xi_j}[ e^{is(\xi_j^2+2\xi_jq_j)}]
\end{equation}
and have
\begin{align*}
&\int d^3\xi f(\xi)\beta_j(|\xi_j+q_j|>2M)e^{i s(\xi_j^2+2\xi_j q_j)}\\=&\frac{1}{(2is)^2}\int d^3\xi \partial_{\xi_j}[\frac{1}{(\xi_j+q_j)}\times \partial_{\xi_j}[\frac{f(\xi)\beta(|\xi_j+q_j|>2M)}{\xi_j+q_j}]]e^{i s(\xi_j^2+2\xi_j q_j)}
\end{align*}
with
\begin{align}
&\partial_{\xi_j}[\frac{1}{(\xi_j+q_j)}\times \partial_{\xi_j}[\frac{f(\xi)\beta(|\xi_j+q_j|>2M)}{\xi_j+q_j}]]\\
=&\partial_{\xi_j}^2[f(\xi)]\frac{\beta(|\xi_j+q_j|>2M)}{(\xi_j+q_j)^2}+f(\xi)\partial_{\xi_j}[ \frac{1}{(\xi_j+q_j)}\times \partial_{\xi_j}[\frac{\beta(|\xi_j+q_j|>2M)}{\xi_j+q_j}]]+\\
&\partial_{\xi_j}[f(\xi)]\left[\partial_{\xi_j}[ \frac{\beta(|\xi_j+q_j|>2M)}{(\xi_j+q_j)^2}]+\frac{1}{(\xi_j+q_j)}\partial_{\xi_j}[\frac{\beta(|\xi_j+q_j|>2M)}{\xi_j+q_j}]\right]\\
=:&\partial_{\xi_j}^2[f(\xi)]\mathscr{F}[J_{2}](\xi_j+q_j)+f(\xi)\mathscr{F}[J_{0}](\xi_j+q_j)+\partial_{\xi_j}[f(\xi_j)]\mathscr{F}[J_{1}](\xi_j+q_j).\label{Jk}
\end{align}
Then take the integral over $q$ and we have
\begin{align*}
\int d^3\xi &f(\xi)\mathcal{J}_\epsilon(\xi)\psi_j(x)=\sum\limits_{l=0}^2\int d^3\xi \int_{\frac{1}{\sqrt{M}}}^\infty \frac{ds}{(2is)^2}\partial_{\xi_j}^l[f(\xi)]e^{-\epsilon s+is\xi^2}\int dk J_{l}(k)e^{-i\xi_jk}\psi_{j,1}(x+2s\xi-ke_j)\\
&+\int d^3\xi\int_0^{\frac{1}{\sqrt{M}}}ds f(\xi)e^{-\epsilon s+is\xi^2}\psi_j(x+2s\xi)=:\sum\limits_{l=0}^2\int d^3\xi\partial_{\xi_j}^l[f(\xi)]Q_{3(j-1)+l+1}(\xi,\epsilon)\psi(x)
\end{align*}
where for $\psi\in L^q$, $1\leq q\leq \infty$, $j=1,2,3$, $l=1,2$,
\begin{align}
Q_{3(j-1)+0+1}&(\xi,\epsilon)\psi(x):=\int_0^{\frac{1}{M}}ds e^{-\epsilon s+is\xi^2}\psi_j(x+2s\xi)+\\
&\int_{\frac{1}{M}}^\infty \frac{ds}{(2is)^2}e^{-\epsilon s+is\xi^2}\int dk J_{0}(k)e^{-i\xi_jk}\psi_{j,1}(x+2s\xi-ke_j),
\end{align}
\begin{equation}
Q_{3(j-1)+l+1}(\xi,\epsilon)\psi(x):=\int_{\frac{1}{M}}^\infty  \frac{ds}{(2is)^2}e^{-\epsilon s+is\xi^2}\int dk J_{l}(k)e^{-i\xi_jk}\psi_{j,1}(x+2s\xi-ke_j)
\end{equation}
and for the definition of $\psi_j, \psi_{j,1}$, see equation \eqref{psij}. For $\psi_4$,
\begin{equation}
\psi_4=\beta(||\xi|+2q\cdot\hat{\xi}|>2M)\gamma(\xi,q,2M)\psi=:\beta(||\xi|+2\hat{\xi}\cdot P|>2M)\psi_{4,1},
\end{equation}
with $\|\psi_{4,1}(x)\|_{\s^p_x}\leq E_{3,2M}^3\|\psi(x)\|_{\s^p_x}$. For $\int d^3\xi f(\xi)\mathcal{J}_\epsilon(\xi)\psi_4$, we take the integral over $s$ directly.
\begin{equation}
\int_0^\infty ds e^{-\epsilon s+is(\xi^2+2\xi\cdot q)}\beta(||\xi|+2\hat{\xi}\cdot q|>2M)=\frac{-\beta(||\xi|+2\hat{\xi}\cdot q|>2M)}{|\xi|(-\frac{\epsilon}{|\xi|} +i(|\xi|+2\hat{\xi}\cdot q))}=:\frac{1}{|\xi|}\mathscr{F}[J_{4,\epsilon/|\xi|}](||\xi|+2q\cdot\hat{\xi}|).
\end{equation}
Let
\begin{equation}
J_{4,\epsilon}(\lambda):=\mathscr{F}_k^{-1}[\frac{-\beta(|k|>2M )}{-\epsilon +ik} ](\lambda).
\end{equation}
Then
\begin{equation}
\int_0^\infty ds e^{-\epsilon s+is(\xi^2+2\xi\cdot q)}\beta(||\xi|+2\hat{\xi}\cdot q|>2M)=\frac{-\beta(||\xi|+2\hat{\xi}\cdot q|>2M)}{|\xi|(-\frac{\epsilon}{|\xi|} +i(|\xi|+2\hat{\xi}\cdot q))}=\frac{1}{|\xi|}\mathscr{F}[J_{4,\epsilon/|\xi|}](||\xi|+2q\cdot\hat{\xi}|).
\end{equation}
In this case, since $|\xi+q|\leq 2\sqrt{3}M$, $|\xi|\geq |q|-2\sqrt{3}M>M>1$. Then
\begin{equation}
\int d^3\xi f(\xi)\mathcal{J}_\epsilon(\xi)\psi_4=\int d^3\xi f(\xi)Q_{0}(\xi,\epsilon)\psi(x),
\end{equation}
where
\begin{equation}
Q_{0}(\xi,\epsilon)\psi(x):=\frac{\beta(|\xi|>1)}{2|\xi|}\int dk J_{4,\epsilon/|\xi|}(k/2)e^{-i|\xi|k/2}\psi_{4,1}(x-k\hat{\xi}).
\end{equation}
Due to Lemma \ref{g}, we have
\begin{equation}
\|J_j(k)\|_{\s^1_k(\mathbb{R})}\lesssim \frac{1}{M^2},\text{  and }\|J_{4,\epsilon}(k)\|_{\s^1_k(\mathbb{R})}\lesssim \frac{1}{M}, j=0,1,2.\label{12.27.1}
\end{equation}
Hence, combining with \ref{gg2} and equation \eqref{12.27.1}, for $1\leq q\leq\infty$,
\begin{equation}
 \|Q_{l} (\xi,\epsilon)\|_{\s_x^q\to \s_x^q}\lesssim \frac{1}{M}, \text{ for } l=0,1,2, \cdots,9.
\end{equation}
According to the expression of $Q_l(\xi,\epsilon), l=1,\cdots,9$ and Lemma \ref{g},
\begin{equation}
\| Q_l^{(*)}(\xi)\psi(x)\|_{\s^p_x}:=\|\sup\limits_{\epsilon\geq 0}| Q_l(\xi,\epsilon)\psi(x)|\|_{\s^p_x}\lesssim \frac{1}{M}\|\psi(x)\|_{\s^p_x}
\end{equation}
and finish the proof.
\end{proof}
Now we will do the $\s_x^p$ estimates for $I_\epsilon^{(k+1)}\beta(|P|>M)$. We will show that for some sufficiently large $M>0$, $\|I_\epsilon^{(k+1)}\beta(|P|>M)\|_{\s_x^p\to \s_x^p}\leq \frac{C^k}{M^k}$ uniformly in $\epsilon\in [0,1]$. Then according to the same process, $L_x^p$ boundedness of $I^{(*,k+1)}\beta(|P|>M)$ follows as a corollary. We will consider $s_l,\xi_l, $ with $l=1,\cdots,k+1$. When $l=1,\cdots,k$ and when we look at $\xi_l,s_l$, we have to deal with
\begin{equation}
\int d^3\xi_l \hat{V}(\xi_{l+1}-\xi_l)\hat{V}(\xi_l-\xi_{l-1})\mathcal{J}_{\epsilon}(\xi_l)\psi(x).\label{inj}
\end{equation}
Applying Lemma \ref{sJ1} to \eqref{inj}, we obtain that \eqref{inj} is equal to
\begin{equation}
\sum\limits_{j_l=0}^9 \int d^3\xi_l Q_{j_l,1}(\xi_l)[\hat{V}(\xi_{l+1}-\xi_l)\hat{V}(\xi_l-\xi_{l-1})]Q_{j_l}(\xi_l,\epsilon)\psi(x)
\end{equation}
where $Q_{0,1}:=identity$ and for $j_l=1,\cdots,9$,
\begin{equation}
Q_{j_l,1}(\xi_l):=\partial_{\xi_l\cdot e_{r_l}}^{m_l}, \text{ with } m_l:=[j_l-1]_3, \quad r_l:=\frac{j_l-1-m_l}{3}+1.\label{Qj1}
\end{equation}
Now we need to introduce some notation. For $j=(j_1,\cdots,j_k)\in \{0,\cdots,9\}^k:=\alpha^k, \xi=(\xi_1,\cdots,\xi_k)\in \mathbb{R}^{3k}, \epsilon>0, k\in\mathbb{N}^+$, define
\begin{equation}
Q_{j}(\xi,\epsilon,k):=\Pi_{l=1}^{k}Q_{j_l}(\xi_l,\epsilon),~ Q_{j,1}(\xi,k):=\Pi_{l=1}^kQ_{j_l,1}(\xi_l).
\end{equation}
\begin{remark}
Here, since $Q_{j_l}(\xi_l,\epsilon)$ commutes with $Q_{j_{l^\prime}}(\xi_{l^\prime},\epsilon)$ and $Q_{j_l,1}(\xi_l)$ commutes with $Q_{j_{l^\prime},1}(\xi_{l^\prime})$ for $l\neq l^\prime$, there is no confusion about $\Pi_{l=1}^kQ_{j_l}(\xi_l,\epsilon)$ and $\Pi_{l=1}^kQ_{j_l,1}(\xi_l)$.
\end{remark}
Then for $\psi\in \beta(|P|>32M)\mathcal{S}_x$,
\begin{equation}
I_\epsilon^{(k+1)}\psi(x)=\sum\limits_{j\in \alpha^{k}}\int d^3\xi_1\cdots d^3\xi_{k}Q_{j,1}(\xi,k)[\mathcal{V}_{k}(\xi)T_{\epsilon}(\xi_{k})]Q_{j}(\xi,\epsilon,k)\psi(x)\label{Niink3}
\end{equation}
where recall that
\begin{equation}
\partial_{\xi_j}^l[T_\epsilon(\xi)]=\int_0^\infty dt e^{itH_0}V(x)\partial_{\xi_j}^l[e^{i\xi_j\cdot x}]e^{-itH_0}
\end{equation}
which is equivalent to the potential $(ix\cdot e_m)^l V(x)e^{i\xi_j\cdot x}$. Now let us look at the $L_x^p$ estimates of $I_\epsilon^{(k)}$ on $\beta(|P|>32)\mathcal{S}_x$.
\begin{lemma}\label{inIkh}  If $V(x)$ satisfies the assumptions of Theorem \ref{intime} and
\begin{equation}
|||V(x)|||_{in}:=\max(\|\hat{V}_a(\xi)\|_{\s^1_\xi}, K_m),
\end{equation}
then for $\psi(x)\in \beta(|P|> 32M)\s_x^p$, $k\geq 1$, $M>1$, there exists some constant $C>0$ such that
\begin{equation}
\|I_\epsilon^{(k+1)}\psi(x)\|_{\s^p_x}\lesssim \frac{C^k|||V(x)|||_{in}^{k+1}}{M^k}\|\psi(x)\|_{\s^p_x}
\end{equation}
for $1\leq p\leq \infty$, $\epsilon \in [0,1]$.
\end{lemma}
\begin{proof}For $p\neq \infty$, choose $\psi\in \beta(|P|>32M)\mathcal{S}_x$. For $l=0,1,2, j=1,2,3$, due to Corollary \ref{T1} and Lemma \ref{sJ1},
\begin{equation}
\|\partial_{\xi_{k}\cdot e_j}^l[T_{\epsilon}(\xi_{k})]Q_j(\xi,\epsilon,k)\psi(x)\|_{\s^p_x}\lesssim K_m\|Q_j(\xi,\epsilon,k)\psi(x)\|_{\s^p_x}\leq \frac{C^kK_m}{M^k}\|\psi(x)\|_{\s^p_x}.\label{Niink}
\end{equation}
The expression
\begin{equation}
\int d^3\xi_1\cdots d^3\xi_kQ_{j,1}(\xi)[\mathcal{V}(\xi,k)T_{\epsilon}(\xi_k)]
\end{equation}
is the sum of $L$ many terms ($L\leq 4^k$) with each term having the form:
\begin{equation}
\frac{1}{(2\pi)^{\frac{3k}{2}}}\int d^3\xi_1\cdots d^3\xi_k P_{\xi_1\cdot e_{j_1}}^{l_1}[\hat{V}(\xi_1)]\cdots P_{\xi_k\cdot e_{j_k}}^{l_k}[\hat{V}(\xi_k-\xi_{k-1})]\partial_{\xi_k\cdot e_{j_{k+1}}}^{l_{k+1}}[T_{\epsilon}(\xi_k)],
\end{equation}
for $j_m\in \{1,2,3\}$, $l_m\in \{0,1,2,3,4\}$, $m=1,\cdots,k$, $l_{k+1}\in \{0,1,2\}$, $j_{k+1}\in \{1,2,3\}$.
According to equation \eqref{Niink} and Lemma \ref{power},
\begin{equation}
\left\| \int d^3\xi_1\cdots d^3\xi_{k}Q_{j,1}(\xi,k)[\mathcal{V}_{k}(\xi)T_{\epsilon}(\xi_{k})]Q_{j}(\xi,\epsilon,k)\psi(x)\right\|_{\s^p_x}\lesssim \frac{C_2^{k}|||V(x)|||_{in}^kK_m}{M^{k}}\|\psi(x)\|_{\s^p_x}.\label{Niink2}
\end{equation}
Then according to equation \eqref{Niink2} and \eqref{Niink3},
\begin{equation}
\|I_\epsilon^{(k+1)}\psi(x)\|_{\s^p}\lesssim \sum\limits_{j\in \alpha^k}  \frac{C_2^{k}|||V(x)|||_{in}^kK_m}{M^{k}}\|\psi(x)\|_{\s^p_x}\lesssim \frac{C^{k}|||V(x)|||_{in}^{k+1}}{M^{k}}\|\psi(x)\|_{\s^p_x}
\end{equation}
for some constant $C>0$. Then by B.L.T. theorem, we get the conclusion for $1\leq p<\infty$. For $p=\infty$, we work on $\beta(|P|>32M) C_0$ first. Then by using duality twice, we get the estimates for $p=\infty$.
\end{proof}
\begin{corollary} If $V(x)$ satisfies the assumptions in Theorem \ref{intime} and
\begin{equation}
|||V(x)|||_{in}:=\max(\|\hat{V}_a(\xi)\|_{\s^1_\xi}, K_m),
\end{equation}
then for $\psi(x)\in \beta(|P|> 32M)L_x^p$, $k\geq 1$, $M>1$, there exists some constant $C>0$ such that
\begin{equation}
\|I^{(*,k+1)}\psi(x)\|_{\s^p_x}\lesssim \frac{C^k|||V(x)|||_{in}^{k+1}}{M^k}\|\psi(x)\|_{\s^p_x}
\end{equation}
for $1\leq p\leq \infty$.
\end{corollary}
\begin{proof}It follows from the same proof of Lemma \ref{inIkh} by replacing $I_\epsilon^{(k+1)}\psi(x)$ with $I^{(*,k+1)}\psi(x)$.
\end{proof}
Now we prove Theorem \ref{intime}.
 \begin{proof}According to Lemma \ref{inIkh}, we have that for $M>C|||V(x)|||_{in}$ and for $\psi\in \mathcal{S}$, $1\leq p\leq \infty$, any $\epsilon\in [0,1]$,
 \begin{equation}
 \| \Omega_\epsilon\beta(|P|>32M)\psi(x)\|_{\s^p_x}\lesssim \left(1+\frac{|||V(x)|||_{in}}{1-\frac{C|||V(x)|||_{in}}{\sqrt{M}}}\right)E_3\|\psi(x)\|_{\s^p_x}
 \end{equation}
 and with $\Omega^{(*)}\beta(|P|>32M)\psi:=\sup\limits_{\epsilon\in (0,1]} \Omega_\epsilon\beta(|P|>32M)\psi$,
 \begin{equation}
\| \Omega^{(*)}\beta(|P|>32M)\psi(x)\|_{\s^p_x}\lesssim \left(1+\frac{|||V(x)|||_{in}}{1-\frac{C|||V(x)|||_{in}}{\sqrt{M}}}\right)E_3\|\psi(x)\|_{\s^p_x}
 \end{equation}
which completes the proof of $\Omega_\epsilon\beta(|P|>32M) \to \Omega_0\beta(|P|>32M)=\Omega\beta(|P|>32M) $ in strong $\s^p$-sense,  provided that the almost everywhere convergence of $\Omega_\epsilon\beta(|H_0|>M)$ to $\Omega\beta(|H_0|>M)$ is a consequence of the $\s^p$ boundedness of $\Omega^{(*)}\beta(|H_0|>M)$ and of Theorem 2.1.14 in \cite{G2008}. By duality, we get the same result for $\beta(|P|>32M)\Omega^*$ and we finish the proof.
 \end{proof}
 \begin{remark}From the proof, based on such definition of $\Omega_\epsilon$, we can see that the result comes from that $\Omega_\epsilon: \s^p\to \s^p$, is bounded uniformly in $\epsilon\in [0,1]$.
\end{remark}
Step further, we get asymptotic completeness on high frequency subspace.
\begin{corollary}If $V(x)$ satisfies the condition in Theorem \ref{intime}, the Schr\"odinger equation has asymptotic completeness on high frequency subspace.
\end{corollary}
\begin{corollary}\label{Oct.25.1}If $V(x)$ satisfies the assumptions in Theorem \ref{intime} and
\begin{equation}
|||V(x)|||_{in}:=\max(\|\hat{V}_a(\xi)\|_{\s^1_\xi}, K_m),
\end{equation}
then for $\psi(x)\in \s_x^p$, $k\geq 1$, $M>1$, there exists some constant $C>0$ such that
\begin{equation}
\|\beta(|P|> 32M)\left(I^{(*,k+1)}\right)^*\psi(x)\|_{\s^p_x}\lesssim \frac{C^k|||V(x)|||_{in}^{k+1}}{M^k}\|\psi(x)\|_{\s^p_x}
\end{equation}
for $1\leq p\leq \infty$, where
\begin{equation}
\left(I^{(*,k+1)}\right)^*:=\max\limits_{\epsilon>0}\left( I_\epsilon^{(k+1)}\right)^*.
\end{equation}
\end{corollary}
\begin{proof} According to Lemma \ref{inIkh}, for $1\leq p<\infty$, by duality, we get the conclusion. When $p=\infty$, choose $\phi\in \s^1_x, \psi\in \s^\infty_x$
\begin{align}
&\left|\langle \phi(x),\beta(|P|> 32M)\left(I^{(*,k+1)}\right)^*\psi(x) \rangle_{L^2_x}\right|\\
=& \left|\langle I^{(*,k+1)}\beta(|P|> 32M)\phi(x),\psi(x) \rangle_{L^2_x}\right|\\
\leq & \|I^{(*,k+1)}\beta(|P|> 32M) \|_{\s^1_x\to \s^1_x}\|\phi(x)\|_{\s^1_x}\|\psi(x)\|_{\s^\infty_x}.
\end{align}
So we conclude that for $\psi\in \s^\infty_x$,
\begin{equation}
\|\beta(|P|> 32M)\left(I^{(*,k+1)}\right)^*\psi(x)\|_{\s^\infty_x}\lesssim \frac{C^k|||V(x)|||_{in}^{k+1}}{M^k}\|\psi(x)\|_{\s^\infty_x}.
\end{equation}

\end{proof}

 \begin{corollary}\label{norm1}If $V(x)$ satisfies the assumptions in Theorem \ref{intime}, and if $P_c$ is bounded on $\s^p$ for all $1\leq p\leq \infty$, there exists $M=M(V(x))>0$ such that
 \begin{equation}
 \sup\limits_{T\in \mathbb{R}}\| P_cU(0,T)e^{-iTH_0}\beta(|P|>M)\|_{\s_x^p\to \s_x^p}<C.
 \end{equation}
 \end{corollary}
 \begin{proof} Due to Theorem \ref{main3}, there exists $M>0$ such that
  \begin{equation}
 \Omega\beta(|P|>M)=s\text{-}\lim\limits_{\epsilon\downarrow 0}\Omega_\epsilon\beta(|P|>M).
 \end{equation}
 Then
 \begin{equation}
 \lim\limits_{\epsilon\downarrow0} (f, \int_0^\infty dt e^{-\epsilon t}\Omega^\prime(t)\beta(|P|>M)g)_{L^2_x}=(f, \lim\limits_{\epsilon\downarrow0} \int_0^\infty dt e^{-\epsilon t}\Omega^\prime(t)\beta(|P|>M)g)_{L^2_x}.
 \end{equation}
 Let
 \begin{equation}
 a(T,f,g):=(f,  P_cU(0,T)e^{-iTH_0}\beta(|P|>M)g)_{L^2_x}.
 \end{equation}
Then for each $f\in \s^p$, $g\in \s^q$, $a(T,f,g)$ is continuous in $T$ since for $t_1,t_2\in \mathbb{R}$,
\begin{equation}
\|\int_{t_1}^{t_2} dt \Omega^\prime(t)\beta(|P|>M)\|_{\s_x^p\to \s_x^p}<\infty
\end{equation}
and goes to 0 as $t_1\to t_2$. Due to Theorem \ref{main3}, we have $\lim\limits_{T\to \infty}a(T,f,g)$ exists for each pair $f,g$. Combining with the continuity, for each $g\in \s^q$,
\begin{equation}
\sup\limits_{T\in \mathbb{R}^+}|a(T,f,g)|<C(f,g).
\end{equation}
By Principle of uniform boundedness,
\begin{equation}
\sup\limits_{\|g \|_{\s^q}\leq 1}\sup\limits_{T\in \mathbb{R}^+}|a(T,f,g)|<C(f),
\end{equation}
that is,
\begin{equation}
\sup\limits_{T\in \mathbb{R}^+}\| P_cU(0,T)e^{-iTH_0}\beta(|P|>M)f\|_{\s^p_x}< C(f).
\end{equation}
Then by Principle of uniform boundedness again and duality,
\begin{equation}
\sup\limits_{T\in \mathbb{R}^+}\| P_cU(0,T)e^{-iTH_0}\beta(|P|>M)\|_{\s^p_x\to \s^p_x}<C.
\end{equation}
Similarly, we have
 \begin{equation}
\sup\limits_{T\in \mathbb{R}^-}\|P_cU(0,T)e^{-iTH_0}\beta(|P|>M)\|_{\s^p_x\to \s^p_x}<C.
\end{equation}
Thus,
\begin{equation}
\sup\limits_{T\in \mathbb{R}}\| P_cU(0,T)e^{-iTH_0}\beta(|P|>M)\|_{\s^p_x\to \s^p_x}<C.
\end{equation}

 \end{proof}
\begin{corollary}If $V(x)$ satisfies the assumptions in Theorem \ref{intime} and $V(x)$ is sufficiently small, then \begin{equation}
\Omega=s\text{-}\lim\limits_{\epsilon \downarrow0}\Omega_\epsilon, \text{ in }\s_x^p, 1\leq p\leq\infty
\end{equation}
and $\Omega^*,\Omega: \s_x^p\to \s_x^p$ are bounded.
\end{corollary}
\begin{proof}
In $I^{(k)}_\epsilon$, for $s_j,\xi_j$, we have to deal with
\begin{equation}
 \int d^3\xi_j\int_0^\infty ds_j\hat{V}(\xi_{j+1}-\xi_j)\hat{V}(\xi_j-\xi_{j-1})e^{-s_j\epsilon+is_j(\xi_j^2+2\xi_j\cdot P)}.
\end{equation}
We do change of variables $s_j\to u_j=s_j|\xi_j|$, $j=1,\cdots, k-1$. For $u_j\leq 1$, we leave as is. For $u_j>1$, we do integration by parts in $|\xi_j|$ twice by setting
\[
e^{iu_j|\xi_j|}=\frac{1}{iu_j}\partial_{|\xi_j|}[e^{iu_j|\xi_j|}].
\]
For $j=k$, we apply Corollary \ref{T1} and for $I_\epsilon^{(k)}$,
\begin{equation}
\| I_\epsilon^{(k)}\|_{\s^p_x\to \s^p_x}\leq C^k |||V(x) |||_{in}^k, \text{ for some }C , \text{ independent on }V(x).
\end{equation}
and
\begin{equation}
\| I^{(*,k)}\|_{\s^p_x\to \s^p_x}\leq C^k |||V(x) |||_{in}^k.
\end{equation}
Then if $ |||V(x) |||_{in}$ is sufficiently small,  the conclusion follows.
\end{proof}

\section{\texorpdfstring{$\s^p$}{\s^p} boundedness of wave operator for some time-dependent potentials}\label{section 3}
In this section, we begin the analysis of  time-dependent potentials. We will show the $\s^p$ boundedness of the wave operator on the high frequency subspace for a class of Mikhlin-type potentials $V(x,t)$ satisfying
\begin{equation}
\sup\limits_{t\in \mathbb{R}}\frac{(1+|t|)^a}{a!}\sum\limits_{l,j=0}^2\sum\limits_{m,r=1}^3|\frac{\partial^a}{\partial t^a}\left[\partial_{\xi\cdot e_r}^l\partial_{\xi\cdot e_m}^j\hat{V}(\xi,t)\right]|\leq c^a \hat{V}_0(\xi), \text{ for all }a\in \mathbb{N}, \text{ some }c\geq 1 \label{pp1}
\end{equation}
with $\hat{V}_0(\xi)\in \s^1_\xi (\mathbb{R}^3)\cap \s^\infty_\xi(\mathbb{R}^3)$.
\subsection{\texorpdfstring{$\s^p$}{Lp} boundedness for \texorpdfstring{$I\sK$}{I}}
We show $\s^p$ boundedness of $I\sK$ when $V(x,t)$ satisfies the condition
\begin{equation}
|||V(x,t)|||_{W1}:=\|\sup\limits_{t\in \mathbb{R}}4\pi\sum\limits_{l,j=0}^2\sum\limits_{m=1}^3(|t|+1)^l|\partial_{\xi\cdot e_m}^j\partial_t^l\hat{V}(\xi, t)|\|_{\s^1_\xi\cap \s^\infty_\xi}<\infty.\label{12.9.eq1}
\end{equation}
\begin{lemma}\label{reWV} If $  V(x,t)$ satisfies assumption \eqref{pp1}, then
\begin{equation}
||| (x\cdot e_j)^lV(x,t)e^{i\eta\cdot x}|||_{W1} <4\pi(1+c+2c^2)\|\hat{V}_0(\xi)\|_{\s^1_\xi\cap \s^\infty_\xi}
\end{equation}
for any $\eta\in \mathbb{R}^3$, $j=1,2,3$, $l=0,1,2$.
\end{lemma}
\begin{proof}Due to assumption \eqref{pp1} and the definition of $||| \cdot |||_{W1}$,
\begin{equation}
||| (x\cdot e_j)^lV(x,t)e^{i\eta\cdot x}|||_{W1} \leq  4\pi\|(0!c^0+1!c+2!c^2)\hat{V}_0(\xi-\eta)\|_{\s^1_\xi\cap \s^\infty_\xi}=4\pi(1+c+2c^2)\|\hat{V}_0(\xi)\|_{\s^1_\xi\cap \s^\infty_\xi}.
\end{equation}
\end{proof}
\begin{theorem}\label{3I1}If $V(x,t)$ satisfies the assumption \eqref{12.9.eq1}, then $I_\epsilon: \s_x^p\to \s_x^p$ is uniformly bounded in $\epsilon\in [0,1]$ for $1\leq p\leq \infty$.
\end{theorem}
\begin{proof} By the same transformation in equation \eqref{Ig1}, we get
\begin{equation}
I_\epsilon\psi(x)=\frac{1}{(2\pi)^{\frac{3}{2}}}\int_{S^2} d\sigma(\xi)\int_0^\infty du\int_0^\infty d|\xi||\xi|\hat{V}(\xi,\frac{u}{|\xi|})e^{-\epsilon \frac{u}{|\xi|}+i(x\cdot \xi +u|\xi|)}\psi(x+2u\hat{\xi}).
\end{equation}
Rewrite $I_\epsilon\psi(x)$ as
\begin{align*}
I_\epsilon\psi(x)=&\frac{1}{(2\pi)^{\frac{3}{2}}}\int d^3\xi\int_0^\infty\chi(|x\cdot \hat{\xi}+u|\leq 1) du \frac{\hat{V}(\xi,\frac{u}{|\xi|})}{|\xi|}e^{-\epsilon \frac{u}{|\xi|}+i(x\cdot \xi +t\xi^2)}\psi(x+2u\hat{\xi})+\\
&\frac{1}{(2\pi)^{\frac{3}{2}}}\int_{S^2}d\sigma(\xi)\int_0^\infty du\int_0^\infty d|\xi|  \chi(|x\cdot \hat{\xi}+u|>1) |\xi|\hat{V}(\xi,\frac{u}{|\xi|})e^{-\epsilon \frac{u}{|\xi|}+i(x\cdot \xi +u|\xi|)}\psi(x+2u\hat{\xi})\\
:=&I_{1\epsilon}\psi(x)+I_{2\epsilon}\psi(x).
\end{align*}
For $I_{1\epsilon}\psi(x)$, due to Lemma \ref{c1} for any $\hat{\xi}$ direction($\chi(|x\cdot\hat{\xi}+u|\leq 1)f(\frac{u}{|\xi|})\in \s^1_u$),
\begin{equation}
\| I_{1\epsilon}\psi(x)\|_{\s^p_x}\lesssim \|\sup\limits_{u\in\mathbb{R}}\frac{|\hat{V}(\xi,\frac{u}{|\xi|})|}{|\xi|}\|_{\s^1_\xi}\|\psi(x)\|_{\s^p_x}\lesssim  |||V(x,t)|||_{W1}\|\psi(x)\|_{\s^p_x}\label{timeIeq1}
\end{equation}
where we use the inequality
\begin{align*}
&\|\frac{\sup\limits_{u\in\mathbb{R}}|\hat{V}(\xi,\frac{u}{|\xi|})|}{|\xi|}\|_{\s^1_\xi}=\| \frac{\chi(|\xi|\geq 1)\sup\limits_{u\in\mathbb{R}}|\hat{V}(\xi,\frac{u}{|\xi|})|}{|\xi|}\|_{\s^1_\xi}+\| \frac{\chi(|\xi|< 1)\sup\limits_{u\in\mathbb{R}}|\hat{V}(\xi,\frac{u}{|\xi|})|}{|\xi|}\|_{\s^1_\xi} \\
\leq &\|\sup\limits_{u\in\mathbb{R}}|\hat{V}(\xi,\frac{u}{|\xi|})|\|_{\s^1_\xi}+\int_{S^2}d\sigma(\xi)\int_0^1 (d|\xi|)|\xi|\|\sup\limits_{u\in\mathbb{R}}| \hat{V}(\xi, \frac{u}{|\xi|})|\|_{\s^\infty_\xi}\leq ||| V(x,t)|||_{W1} .
\end{align*}
For $I_{2\epsilon}\psi(x)$, we do integration by parts in $|\xi|$ in the same way as time-independent case and we have
\begin{align*}
&I_{2\epsilon}\psi(x)=\frac{-1}{(2\pi)^{\frac{3}{2}}}\int_0^\infty du\int_{S^2}d\sigma(\xi)\frac{\chi(|x\cdot\hat{\xi}+u|>1)}{i(x\cdot\hat{\xi}+u)}\psi(x+2u\hat{\xi})\\
&\left[\int_0^{\frac{1}{\sqrt{|x\cdot\hat{\xi}+u|}}} d|\xi|\partial_{|\xi|}[|\xi|\hat{V}(\xi,\frac{u}{|\xi|})e^{-\epsilon \frac{u}{|\xi|}}]e^{i(x\cdot\xi+u|\xi|)}+\int_{\frac{1}{\sqrt{|x\cdot\hat{\xi}+u|}}}^\infty d|\xi|\partial_{|\xi|}[|\xi|\hat{V}(\xi,\frac{u}{|\xi|})e^{-\epsilon \frac{u}{|\xi|}}]e^{i(x\cdot\xi+u|\xi|)}\right]\\
&:=I_{21\epsilon}\psi(x)+I_{22\epsilon}\psi(x)
\end{align*}
where we throw away the boundary terms both near infinity and near $0$ due to our assumptions:
\begin{equation}
\|\sup\limits_{t\in \mathbb{R}}|\hat{V}(\xi,t)|\|_{\s^\infty_\xi} \implies  |\xi|\hat{V}(\xi,\frac{u}{|\xi|})\vert_{|\xi|=0}=0
\end{equation}
and due to the definition of $|||\cdot |||_{W1}$,
\begin{align*}
\hat{V}_0(\xi):=\|\sup\limits_{t\in \mathbb{R}}|\hat{V}(\xi,t)|\|_{\s^1_\xi} \implies &\exists r_n(r_n\to \infty \text{ as }n\to \infty) s.t. \\
&r_n|\hat{V}(r_n\hat{\xi},\frac{u}{r_n})|\leq r_n\hat{V}_0(r_n\hat{\xi})\to 0 ,\text{ as }n\to \infty.
\end{align*}
For $I_{21\epsilon}\psi(x)$, we have
\begin{align*}
\| I_{21\epsilon}\psi(x)\|_{\s^p_x}\leq &\left\|\int_0^\infty \frac{du}{(2\pi)^{3/2}}\int_{S^2}d\sigma(\xi)\frac{\chi(|x\cdot\hat{\xi}+u|>1)}{|x\cdot\hat{\xi}+u|^{\frac{3}{2}}}\|\partial_{|\xi|}[|\xi|\hat{V}(\xi, \frac{u}{|\xi|})e^{-\epsilon \frac{u}{|\xi|}}] \|_{\s^\infty_{|\xi|}[0,1]}|\psi(x+2u\hat{\xi})| \right\|_{\s^p_x}\\
(\text{Lemma \ref{c1} })\lesssim &\int_{S^2}d\sigma(\xi)\|\sup\limits_{u\in \mathbb{R}^+}|\partial_{|\xi|}[|\xi|\hat{V}(\xi,\frac{u}{|\xi|})e^{-\epsilon \frac{u}{|\xi|}}]| \|_{\s^\infty_{|\xi|}[0,1]}\|\psi(x)\|_{\s^p_x}\\
\lesssim &|||V(x,t)|||_{W1}\|\psi(x)\|_{\s^p_x}
\end{align*}
where from the second line to the third line, we use ($\partial_1[\hat{V}(\xi,t)]:=\partial_{|\xi|}[\hat{V}(\xi,t)]$, $\partial_2[\hat{V}(\xi,t)]:=\partial_{t}[\hat{V}(\xi,t)]$.)
\begin{align*}
&\| \sup\limits_{u\in \mathbb{R}^+}|\partial_{|\xi|}[|\xi|\hat{V}(\xi,\frac{u}{|\xi|})e^{-\epsilon \frac{u}{|\xi|}}]|\|_{\s^\infty_{|\xi|}[0,1]}\\
\leq&\|\sup\limits_{u\in \mathbb{R}^+} \left[(1+\frac{\epsilon u}{|\xi|})|\hat{V}(\xi,\frac{u}{|\xi|})|e^{-\epsilon \frac{u}{|\xi|}}+|\xi||\partial_1[\hat{V}(\xi, \frac{u}{|\xi|})]|+\frac{|u|}{|\xi|}|\partial_2[\hat{V}(\xi, \frac{u}{|\xi|})]|\right]\|_{\s^\infty_{|\xi|}[0,1]}\\
\lesssim& |||V(x,t)|||_{W1}.
\end{align*}
For $I_{22\epsilon}\psi(x)$, we do integration by parts in $|\xi|$ in the same way again, and have
\begin{align*}
&I_{22\epsilon}\psi(x)=\frac{-1}{(2\pi)^{\frac{3}{2}}}\int_0^\infty du\int_{S^2}d\sigma(\xi)\frac{\chi(|x\cdot\hat{\xi}+u|>1)}{(i(x\cdot\hat{\xi}+u))^2}\psi(x+2u\hat{\xi})\\
&\left[ \partial_{|\xi|}[|\xi|\hat{V}(\xi,\frac{u}{|\xi|})e^{-\epsilon \frac{u}{|\xi|}}]e^{i(x\cdot\xi+u|\xi|)}\vert_{|\xi|=\frac{1}{\sqrt{|x\cdot\hat{\xi}+u|}}}^{|\xi|=\infty}-\int_{\frac{1}{\sqrt{|x\cdot\hat{\xi}+u|}}}^\infty d|\xi| \partial^2_{|\xi|}[|\xi|\hat{V}(\xi,\frac{u}{|\xi|})e^{-\epsilon \frac{u}{|\xi|}}]e^{i(x\cdot\xi+u|\xi|)}\right].
\end{align*}
Then similarly, take absolute value in the integral, use Lemma \ref{c1} and compute the $L^p_x$ norm of $I_{22\epsilon}\psi(x)$:
\begin{equation}
\| I_{22\epsilon}\psi(x)\|_{\s^p_x}\lesssim |||V(x,t)|||_{W1}\|\psi(x)\|_{\s^p_x} .
\end{equation}
Then we have
\begin{equation}
\|I_{2\epsilon}\psi(x)\|_{\s^p_x}\lesssim |||V(x,t)|||_{W1}\|\psi(x)\|_{\s^p_x}.\label{timeIeq2}
\end{equation}
Hence, according to equation \eqref{timeIeq1} and equation \eqref{timeIeq2}, we obtain
\begin{equation}
\|I_\epsilon\psi(x)\|_{\s^p_x}\lesssim |||V(x,t)|||_{W1}\|\psi(x)\|_{\s^p_x}.
\end{equation}
\end{proof}
\begin{corollary}\label{3T1}Let
\begin{equation}
T_\epsilon^{[k]}(\eta)\psi(x)=\int_0^\infty dt e^{iH_0t}f^{[k]}(t)e^{-\epsilon t}(x\cdot e_m)^lV(x,t)e^{ix\cdot \eta}e^{-iH_0t}\psi(x),
\end{equation}
for $\psi\in L^p, a_j\geq 0, \eta\in \mathbb{R}^3,$ $k\in \mathbb{N}^+, l=0,1,2, e_m\in S^2,$ with
\begin{equation}
f^{[k]}(t)=\Pi_{j=1}^k f_{j}(a_j+t),~ a_j\geq 0, \quad  \sup\limits_{t\in \mathbb{R}^+}|t|^a|f^{(a)}_{j}(t)|\leq C_j, \text{ for }a=0,1,2, \text{ and for some }C_j>1.
\end{equation}
If $V(x,t)$ satisfies the condition \eqref{pp1}, then $T_\epsilon^{[k]}: \s_x^p\to \s_x^p$ is uniformly bounded in $\epsilon \in [0,1]$ for all $1\leq p\leq\infty$ and
\begin{equation}
\|T^{[k]}_\epsilon(\eta)\|_{\s_x^p\to \s_x^p}\lesssim c^2k^2 (\Pi_{j=1}^kC_j )\|\hat{V}_0(\xi)\|_{\s^1_\xi\cap \s^\infty_\xi}.
\end{equation}
\end{corollary}
\begin{proof} Replace $ V(x,t)$ with $ V(x,t)e^{i\eta\cdot x}f^{[k]}(t) $ in the proof of Theorem \ref{3I1}. Since for $t\geq 0$,
\begin{equation}
\left|t^j\frac{d^j[f_l(t+a_j)]}{dt}\right|\leq |(t+a_j)^j \frac{d^j[f_l(t+a_j)]}{dt^j}|\leq C_l, \text{ for }j=0,1,2, \quad l=1,\cdots ,k,
\end{equation}
based on Leibniz formula,
\begin{equation}
\left|t^j\frac{d^j[f^{[k]}(t)]}{dt}\right|\leq  k^j\Pi_{l=1}^kC_l,\text{ for }j=0,1,2, a\geq 0.\label{kTmain}
\end{equation}
Then for $l=0,1,2$,
\begin{align}
&4\pi\sum\limits_{u=0}^2\sum\limits_{r=1}^3(|t|+1)^u |\partial_{\xi\cdot e_r}^j\partial_t^u [f^{[k]}(t)\partial_{\xi\cdot e_m}^l[\hat{V}(\xi-\eta,t)]]|\\\leq& \sum\limits_{u=0}^2\sum\limits_{r=1}^3\sum\limits_{l_1=0}^u \begin{pmatrix}
u\\
l_1
\end{pmatrix}4\pi (|t|+1)^{l_1}|\partial_{\xi\cdot e_r}^j\partial_t^{l_1}\partial_{\xi\cdot e_m}^l[\hat{V}(\xi-\eta, t)]|\times k^2\Pi_{l=1}^kC_l.\label{kTmain2}
\end{align}
Hence, due to Lemma \ref{reWV} and equation \eqref{kTmain2},
\begin{equation}
||| f^{[k]}(t)(x\cdot e_m)^lV(x,t)e^{ix\cdot \eta} |||_{W1}\lesssim c^2k^2(\Pi_{l=1}^kC_l)\|\hat{V}_0(\xi)\|_{\s^1_\xi\cap \s^\infty_\xi}<\infty.
\end{equation}
Apply Theorem \ref{3I1} and we finish the proof.
\end{proof}

\subsection{\texorpdfstring{$\s_x^p$}{Lxp} boundedness for \texorpdfstring{$I_\epsilon^{(k)}$}{Ie(k)} on high frequency space}
In this section, we use the following notation. For $\alpha\in \{0,1\}^k$, let
\begin{equation}
\mathcal{V}(\xi,s,k):=\frac{1}{(2\pi)^{\frac{3k}{2}}}\Pi_{j=1}^k\hat{V}(\xi_j-\xi_{j-1},\sum\limits_{l=j}^ks_l)
\end{equation}
for $\xi=(\xi_1,\cdots,\xi_k)\in \mathbb{R}^{3k}, \xi_0=0, s=(s_1,\cdots,s_k)\in \mathbb{R}^k$. For $\psi\in \s_x^q$, $1\leq q\leq \infty$, $j=1,2,3$, $l=1,2$, let
\begin{equation}
Q_{3(j-1)+l+1}^1(\xi,\epsilon,s)\psi(x):=\frac{\chi(s>\frac{1}{M})}{(2is)^2}e^{-\epsilon s+is\xi^2}\int dk J_l(k)e^{-i\xi_jk}\psi_{j,1}(x+2s\xi-ke_j),
\end{equation}
\begin{align}
Q_{3(j-1)+0+1}^{1}(\xi,\epsilon,s)\psi(x):=&\chi(s\leq \frac{1}{M}) e^{-\epsilon s+is\xi^2}\psi_j(x+2s\xi)+\\
&\frac{\chi(s>\frac{1}{M})}{(2is)^2}e^{-\epsilon s+is\xi^2}\int dk J_0(k)e^{-i\xi_jk}\psi_{j,1}(x+2s\xi-ke_j).
\end{align}
Here we recall the definition of $J_l, \psi_j,\psi_{j,1}$, see \eqref{Jk}, \eqref{psij}, \eqref{psij1}. Then
\begin{equation}
\int_0^\infty ds Q_{3(j-1)+0+1}^{1}(\xi,\epsilon,s)\psi(x)=Q_{3(j-1)+0+1}(\xi,\epsilon)\psi(x),
\end{equation}
\begin{equation}
\int_0^\infty  ds Q_{3(j-1)+l+1}^1(\xi,\epsilon,s)\psi(x)=Q_{3(j-1)+l+1}(\xi,\epsilon)\psi(x).
\end{equation}
We immediately have the following lemma:
\begin{lemma}\label{sbd}For $j=,1,2,3,$ $l=1,2$, $1\leq p\leq\infty$,
\begin{equation}
\int_0^\infty ds \left\| Q_{3(j-1)+0+1}^{1}(\xi,\epsilon,s)\right\|_{\s^p_x\to \s^p_x}\lesssim  \frac{1}{M},
\end{equation}
\begin{equation}
\int_0^\infty ds  \left\| Q_{3(j-1)+l+1}^{1}(\xi,\epsilon,s)\right\|_{\s^p_x\to \s^p_x}\lesssim \frac{1}{M}.
\end{equation}
\end{lemma}
\begin{proof} This follows directly from the proof of Lemma \ref{sJ1}.
\end{proof}
Now we can get the $\s_x^p$ estimates for $I^{(k)}_\epsilon$:
\begin{lemma}\label{exlnk}If $V(x,t)$ satisfies condition \eqref{pp1}, then for $M\geq 1$, when $\psi\in \beta(|P|>32M)\mathcal{S}_x$, $\epsilon\geq 0$,
\begin{equation}
\| I_\epsilon^{(k)}\psi(x)\|_{\s^p_x}\lesssim \frac{ C^kc^{3k+2}k^3\|\hat{V}_0(\xi)\|_{\s^1_\xi\cap \s^\infty_\xi}^{k}}{M^{k-1}}\|\psi(x)\|_{\s^p_x},
\end{equation}
and
\begin{equation}
\| \beta(|P|>32M)\left(I_\epsilon^{(k)}\right)^*\|_{\s^p_x\to \s^p_x}\lesssim \frac{ C^kc^{3k+2}k^3\|\hat{V}_0(\xi)\|_{\s^1_\xi\cap \s^\infty_\xi}^{k}}{M^{k-1}},
\end{equation}
for $1\leq p\leq\infty, k\geq 2.$
\end{lemma}
\begin{proof}According to the same transformation in $t_j$ in section $2$, we can rewrite $I_\epsilon^{(k)}\psi(x)$ as
\begin{align*}
I_\epsilon^{(k)}\psi(x)=&\sum\limits_{\gamma\in \{0,1\}^{k-1}}\int_0^\infty ds_k\cdots\int_{0}^\infty ds_1\int
d^3\xi_1 \cdots d^3\xi_k d^3q\beta^{\gamma}(\xi,q,k)\mathcal{V}(\xi,s,k)\\
&e^{-s_k\epsilon-\cdots-s_1\epsilon+i(x\cdot (\xi_k+q)+s_k(\xi_k^2+2q\cdot\xi_k)+\cdots+s_1(\xi_1^2+2\xi_1\cdot q))}\frac{\hat{\psi}(q)}{(2\pi)^{\frac{3}{2}}}=: \sum\limits_{\gamma\in \{0,1\}^{k-1}}I_{\gamma\epsilon}^{(k)}\psi(x),
\end{align*}
where
\begin{equation}
\beta^{\gamma_j}(|\xi_j+q|>2M)=\begin{cases} \beta(|\xi_j+q|>2M)&\text{ if }\gamma_j=0\\  \beta(|\xi_j+q|\leq 2M)&\text{ if }\gamma_j=1\end{cases},\quad\beta^\gamma(\xi,q,k)=\Pi_{j=1}^{k-1}\beta^{\gamma_j}(|\xi_j+q|>2M).
\end{equation}
For $I_{\gamma\epsilon}^{(k)}\psi(x)$, if $\gamma_j=0$ for all $j=1,\cdots,k-1$, the transformation we will take is the same as that in time-independent case. After such a transformation, we use Corollary \ref{3T1} instead of Corollary \ref{T1} and get that in this case,
\begin{equation}
\| I_{\gamma\epsilon}^{(k)}\psi(x)\|_{\s^p_x}\leq \frac{c^2k^2C^k\| \hat{V}_0(\xi)\|_{\s^1_\xi\cap \s^\infty_\xi}^k}{M^{k-1}}\|\psi(x)\|_{\s^p_x}
\end{equation}
for some constant $C>0$. The rest of the task is to deal with  $I_{\gamma\epsilon}^{(k)}\psi(x)$ when there exists some $j$ such that $\gamma_j=1$. In this case, let
\begin{equation}
\{j_1,\cdots, j_r\}:=\{j: |\xi_j+q|\leq 2M \text{ and }j\in\{1,\cdots,k-1\}\}, \text{ with }j_1<\cdots<j_r,
\end{equation}
where $r$ denotes the number of such $j$ with $|\xi_j+s_j|\leq 2M$, $j\leq k-1$.

 In the following, we will use some transformation to get a desired upper bound for such $I_{\gamma\epsilon}^{(k)}\psi(x)$. This transformation  is slightly different from that in time-independent case. \\
\textbf{Transformation :}\par
We do the transformation for $\xi_{l}, s_l$, with $l\in \{j_1,\cdots,j_r\}$ first. Recall that when $|\xi_l+q|\leq 2M$, $||\xi_l|+2q\cdot\hat{\xi}_l|>2M$. We begin with $j_1$. Look at the integral over $s_{j_1}$
\begin{equation}
\int_0^\infty ds_{j_1}e^{-\epsilon s_{j_1}+is_{j_1}(\xi_{j_1}^2+2\xi_{j_1}\cdot q)}\mathcal{V}(\xi,s,k).
\end{equation}
We do integration by parts in $s_{j_1}$ variable by setting
\begin{equation}
e^{-\epsilon s_{j_1}+is_{j_1}(\xi_{j_1}^2+2\xi_{j_1}\cdot q)}=\frac{1}{-\epsilon +i(\xi_{j_1}^2+2\xi_{j_1}\cdot q )}\partial_{s_{j_1}}[ e^{-\epsilon s_{j_1}+is_{j_1}(\xi_{j_1}^2+2\xi_{j_1}\cdot q)}]
\end{equation}
and get two terms: boundary term
\begin{equation}
\frac{-1}{-\epsilon +i(\xi_{j_1}^2+2\xi_{j_1}\cdot q )}=\frac{-1}{-\epsilon +i(\xi_{j_1}^2+2\xi_{j_1}\cdot q )}\int_0^\infty ds_{j_1}\delta(s_{j_1})e^{-\epsilon s_{j_1}+is_{j_1}(\xi_{j_1}^2+2\xi_{j_1}\cdot q)}\mathcal{V}(\xi,s,k)\label{bdterm}
\end{equation}
and integral term
\begin{equation}
\frac{-1}{-\epsilon +i(\xi_{j_1}^2+2\xi_{j_1}\cdot q )}\int_0^\infty ds_{j_1} e^{-\epsilon s_{j_1}+is_{j_1}(\xi_{j_1}^2+2\xi_{j_1}\cdot q)}\partial_{s_{j_1}}[\mathcal{V}(\xi,s,k) ].\label{interm}
\end{equation}
For the boundary term, if $r=1$, we stop. Otherwise, we move to $j_2$ and do the same transformation in $s_{j_2}$. For the integral term, we keep taking integration by parts in $s_{j_1}$ in the same way. We keep doing such transformation for the boundary terms and integration terms for $r+2$  times, and the terms with $\delta(s_{j_1})\cdots \delta(s_{j_r})$ are left out. For the rest $j\in \{1,\cdots,k-1\}$, the transformation is the same as that in time-independent case. To be precise, here are the full set of steps:
\begin{enumerate}
\item \textbf{Transformation for }$\{j_1,\cdots, j_r\}$:\\
\textbf{Step one:} set $l=1,$ $m=0$ and
\begin{equation}
F=\beta^{\gamma}(\xi,q,k)\mathcal{V}(\xi,s,k)e^{-s_k\epsilon-\cdots-s_1\epsilon+i(x\cdot (\xi_k+q)+s_k(\xi_k^2+2q\cdot\xi_k)+\cdots+s_1(\xi_1^2+2\xi_1\cdot q))}.
\end{equation}
\textbf{Step two:} set $m=m+1$ and in $\int_0^\infty ds_{j_l}F$, take integration by parts in $s_{j_l}$ variable by setting
\begin{equation}
e^{-\epsilon s_{j_l}+is_{j_l}(\xi_{j_l}^2+2\xi_{j_l}\cdot q)}=\frac{1}{-\epsilon +i(\xi_{j_l}^2+2\xi_{j_l}\cdot q )}\partial_{s_{j_l}}[ e^{-\epsilon s_{j_l}+is_{j_l}(\xi_{j_l}^2+2\xi_{j_l}\cdot q)}]
\end{equation}
and get two terms: boundary term $-\int_0^\infty ds_{j_l}\delta(s_{j_l})F_1$ and integral term $-\int_0^\infty ds_{j_l}F_2$. For example, when $l=1$, see \eqref{bdterm} and \eqref{interm}. For boundary term, we go to \textbf{Step three} and go to \textbf{Step four} for integral term. \\
\textbf{Step three:} for boundary term $-\int_0^\infty ds_{j_l}\delta(s_{j_l})F_1$, if $l<r$ and $m<r+2$, set $F=F_1$, $l=l+1$ and move back to \textbf{Step two}. Otherwise, (($l<r$ and $m=r+2$) or ($l=r$)) we stop taking transformation on the boundary term.\\
\textbf{Step four:} for integral term, if $m<r+2$, set $F=F_2 $ and move back to \textbf{Step two}. Otherwise, $m=r+2$ and we stop taking transformation on the integral term. \par

\textbf{After these transformation, we get no more than $2^{r+2}$ many sub-terms.} Each term has the form of (we call the case when $m=r+2$, type $1$)
\begin{align*}
(-1)^{r+2}&\int_0^\infty ds_1\cdots \int_0^\infty ds_k \int d^3qd^3\xi_1\cdots d^3\xi_k \delta(s_{j_1})\cdots \delta(s_{j_{m-1}})\partial_{s_{j_1}}^{l_1}\cdots  \partial_{s_{j_m}}^{l_m}[\mathcal{V}](\xi,s,k)\times\\
& 1/\left[ (i(\xi_{j_1}^2+2\xi_{j_1}\cdot q))^{l_1+1}\times \cdots \times (i(\xi_{j_{m-1}}^2+2\xi_{j_{m-1}}\cdot q))^{l_{m-1}+1}\times (i(\xi_{j_m}^2+2\xi_{j_m}\cdot q))^{l_m}\right]\times\\
&\beta^\gamma(\xi,q,k)e^{-s_k\epsilon-\cdots-s_1\epsilon+i(x\cdot (\xi_k+q)+s_k(\xi_k^2+2q\cdot\xi_k)+\cdots+s_1(\xi_1^2+2\xi_1\cdot q))}\frac{\hat{\psi}(q)}{(2\pi)^{\frac{3}{2}}}
\end{align*}
with $m-1+\sum\limits_{u=1}^ml_u=r+2,1\leq m\leq k-1, l_u\geq 0$, or of(we call the case when $m=r+1$, type $2$)
\begin{align*}
(-1)^{r+1}&\int_0^\infty ds_1\cdots \int_0^\infty ds_k \int d^3qd^3\xi_1\cdots d^3\xi_k \delta(s_{j_1})\cdots \delta(s_{j_{r}})\partial_{s_{j_1}}^{l_1}\cdots  \partial_{s_{j_r}}^{l_{r}}[\mathcal{V}](\xi,s,k)\times\\
& 1/\left[ (i(\xi_{j_1}^2+2\xi_{j_1}\cdot q))^{l_1+1}\times \cdots \times (i(\xi_{j_{r}}^2+2\xi_{j_{r}}\cdot q))^{l_{r}+1}\right]\times\\
&\beta^\gamma(\xi,q,k)e^{-s_k\epsilon-\cdots-s_1\epsilon+i(x\cdot (\xi_k+q)+s_k(\xi_k^2+2q\cdot\xi_k)+\cdots+s_1(\xi_1^2+2\xi_1\cdot q))}\frac{\hat{\psi}(q)}{(2\pi)^{\frac{3}{2}}}
\end{align*}
with $\sum\limits_{u=1}^{r}l_u=1, l_u\geq 0$, or of(we call the case when $m=r$, type $3$)
\begin{align*}
(-1)^{r}&\int_0^\infty ds_1\cdots \int_0^\infty ds_k \int d^3qd^3\xi_1\cdots d^3\xi_k \delta(s_{j_1})\cdots \delta(s_{j_{r}})\mathcal{V}(\xi,s,k)\times\\
& 1/\left[ (i(\xi_{j_1}^2+2\xi_{j_1}\cdot q))\times \cdots \times (i(\xi_{j_{r}}^2+2\xi_{j_{r}}\cdot q))\right]\times\\
&\beta^\gamma(\xi,q,k)e^{-s_k\epsilon-\cdots-s_1\epsilon+i(x\cdot (\xi_k+q)+s_k(\xi_k^2+2q\cdot\xi_k)+\cdots+s_1(\xi_1^2+2\xi_1\cdot q))}\frac{\hat{\psi}(q)}{(2\pi)^{\frac{3}{2}}}.
\end{align*}
Here each $ 1/(\xi_{j_u}^2+2\xi_{j_u}\cdot q)$ will give us a factor $C_1/M$ for some fixed constant $C_1>0$.
\item \textbf{Transformation for the rest $j\in \{1,\cdots,k-1\}-\{j_1,\cdots,j_r\}$ }:\\
When it comes to these $j$, for each term, we do the same transformation as before and will gain at least $\frac{C_2}{M}$($C_2$ is some fixed constant) for each $j$ with this property. And according to the definition of $r,$ we have $k-1-r$  such $j$ and will gain $\frac{C_2^{k-1-r}}{M^{k-1-r}}$ from the transformation here.

\end{enumerate}
\textbf{Estimates for all three types:} the estimates are based on how we deal with $j=k$. For type $1$, we do nothing for $\xi_k,s_k$ and defer its $L_x^p$ estimates to the end. \\
\textbf{Estimates for type $2$}: for type $2$, after the transformation to case when $|\xi_j+q|>2M$, it becomes the sum of no more than $81^{k}$ many terms since for
\begin{equation}
\partial_{\xi_j\cdot e_m}^l[\hat{V}(\xi_j-\xi_{j-1}, \sum\limits_{a=j}^ks_a)]Q_r, m\in\{1,2,3\}, j\in\{1,\cdots,k\}, l\in \{0,1,2\}, r\in\{1,\cdots,9\},
\end{equation}
there are $81^k$ many cases. Here for $Q_r$, see Lemma \ref{sJ1}. For each term, when it comes to $\xi_k,s_k$, we have to face
\begin{equation}
\int_0^\infty ds_k \int d^3\xi_k\partial_{s_{j_1}}^{l_1}\cdots  \partial_{s_{j_r}}^{l_{r}}[ f^{[k-1]}(\xi,s)\partial_{\xi_{k-1}\cdot e_v}^{w}[\hat{V}(\xi_k-\xi_{k-1},s_k)]] e^{iH_0s_k}e^{i\xi_k\cdot Q}e^{-iH_0s_k}
\end{equation}
for some direction $e_v$, some $w\in \{0,1,2\}$, with
\begin{equation}
f^{[k-1]}(\xi,s)=\partial_{\xi_{k-1}\cdot e_{m,k-1}}^{w_{k-1}}[\hat{V}(\xi_{k-1}-\xi_{k-2},\sum\limits_{a=k-1}^ks_a)]\times \cdots\times\partial_{\xi_{1}\cdot e_{m,1}}^{w_{1}}[\hat{V}(\xi_{1}-\xi_{0},\sum\limits_{a=1}^ks_a)]\label{fk}
\end{equation}
for some $w_j\in\{0,1,2\}, e_{m,j}\in\{1,2,3\}$.
Since for type $2$, $\sum\limits_{u=1}^r l_u=1$, we have
\begin{align*}
\partial_{s_{j_1}}^{l_1}\cdots  \partial_{s_{j_r}}^{l_{r}}&[f^{[k-1]}(\xi,s)\partial_{\xi_{k-1}\cdot e_v}^w[\hat{V}(\xi_k-\xi_{k-1},s_k)]]=\partial_{s_{j_u}}[f^{[k-1]}(\xi,s)\partial_{\xi_{k-1}\cdot e_v}^w[\hat{V}(\xi_k-\xi_{k-1},s_k)]]\\
=&\sum\limits_{a=1}^{j_u}f^{[k-1]}_a (\xi,s)\partial_{\xi_{k-1}\cdot e_v}^w[\hat{V}(\xi_k-\xi_{k-1},s_k)], \text{ for some }u\in\{1,\cdots,r\},
\end{align*}
where the difference between $f^{[k-1]}_a$ and $f^{[k-1]}$ is that they have a different $a$th factor, that is, in $f^{[k-1]}_a$, for the $a$th factor, it has
\begin{equation}
\partial_{s_{j_u}} \partial_{\xi_{a}\cdot e_{m,a}}^{w_{a}}[\hat{V}(\xi_{a}-\xi_{a-1},\sum\limits_{b=a}^ks_b)]
\end{equation}
instead of
\begin{equation}
 \partial_{\xi_{a}\cdot e_{m,a}}^{w_{a}}[\hat{V}(\xi_{a}-\xi_{a-1},\sum\limits_{b=a}^ks_b)].
\end{equation}
Since for $b=0,1$, $j=0,1,2$, $a=1,\cdots,k$,
\begin{equation}
\sup\limits_{s_k\in \mathbb{R}^+}|s_k|^j|\partial_{s_k}^{j}\partial_{s_{j_u}}^{b} \partial_{\xi_{j}\cdot e_{m,j}}^{w_{j}}[\hat{V}(\xi_{j}-\xi_{j-1},\sum\limits_{b=j}^ks_b)]
 |\leq c^3\hat{V}_0(\xi_j-\xi_{j-1}), \label{12.21.1}
\end{equation}
we can apply Corollary \ref{3T1}, Lemma \ref{power} and have
\begin{align*}
\| \text{type }2\|_{\s^p_x}\lesssim &\frac{C_3^kc^2k\times k^281^k (c^3\|\hat{V}_0(\xi)\|_{\s^1_\xi\cap \s^\infty_\xi})^k}{M^{r+1+(k-r-1)}} \|\psi(x)\|_{\s^p_x}
\end{align*}
where we have another $k$ since $j_u\leq k-1<k$. Therefore
\begin{equation}
\| \text{type }2\|_{\s^p_x}\lesssim \frac{ C_4^kc^{3k+2}k^3\|\hat{V}_0(\xi)\|_{\s^1_\xi\cap \s^\infty_\xi}^{k}}{M^{k}}\|\psi(x)\|_{\s^p_x}.
\end{equation}
\textbf{Estimates for type $3$}: for type 3, similarly, after the transformation to case when $|\xi_j+q|>2M$, it becomes the sum of no more than $9^{k}$ many terms. For each term, when it comes to $\xi_k,s_k$, we have to face the operator
\begin{equation}
\int_0^\infty ds_k \int d^3\xi_kf^{[k-1]}(\xi,s)\partial_{\xi_{k-1}\cdot e_v}^w[\hat{V}_{\alpha_k}(\xi_k-\xi_{k-1},s_k) ]e^{iH_0s_k}e^{i\xi_k\cdot Q}e^{-iH_0s_k}
\end{equation}
with $f^{[k-1]}$ satisfying equation \eqref{fk}. Due to inequality \eqref{12.21.1}, Lemma \ref{power} again, we have
\begin{equation}
\| \text{type }3\|_{\s^p_x}\lesssim \frac{ C_5^kc^281^k k^2(c^{3}\|\hat{V}_0(\xi)\|_{\s^1_\xi\cap \s^\infty_\xi})^k}{M^rM^{k-1-r}}\|\psi(x)\|_{\s^p_x}
\end{equation}
and therefore
\begin{equation}
\| \text{type }3\|_{\s^p_x}\lesssim \frac{ C_6^kc^{3k+2}k^2\|\hat{V}_0(\xi)\|_{\s^1_\xi\cap \s^\infty_\xi}^{k}}{M^{k-1}}\|\psi(x)\|_{\s^p_x}.
\end{equation}
\textbf{Estimates for type $1$}: it requires the following lemma:
 \begin{lemma}\label{key12.21.1}For $1\leq j_1<\cdots<j_m<k$, $\mathbb{N}=\{0,1,\cdots\}$, let
 \begin{equation}
 \mathcal{L}_{m}:=\Pi_{l=1}^m f(s_{j_l}+s_{j_l+1}+\cdots+s_k)
 \end{equation}
 and for $\gamma\in \mathbb{N}^m$,
 \begin{equation}
  \mathcal{L}_m^{\gamma}:=\Pi_{l=1}^m \frac{1}{\gamma_l! }f^{(\gamma_l)}(s_{j_l}+s_{j_l+1}+\cdots+s_k).
 \end{equation}
 If $l_1+\cdots+l_m\leq k+1$, then
\begin{equation}
\partial_{s_{j_1}}^{l_1}\cdots\partial_{s_{j_m}}^{l_m}[\mathcal{L}_m]=\sum\limits_{\gamma\in \mathbb{N}^m, |\gamma|=l_1+\cdots+l_m} c_\gamma \mathcal{L}_m^\gamma\label{12.21.n1}
\end{equation}
with
\begin{equation}
\sum\limits_{\gamma\in \mathbb{N}^m, |\gamma|=l_1+\cdots+l_m}|c_\gamma|\leq (2k)^{l_1+\cdots+l_m}.
\end{equation}
\end{lemma}
 \begin{proof}Let
 \begin{equation}
 \mathcal{M}:=\Pi_{l=1}^m f(s+a_l), \text{ for }a_l>0
 \end{equation}
 and for $\gamma\in \mathbb{N}^m$,
  \begin{equation}
 \mathcal{M}^\gamma:=\Pi_{l=1}^m \frac{1}{\gamma_l!}f^{(\gamma_l)}(s+a_l).
 \end{equation}
 Since
 \begin{equation}
 \partial_{s}[\mathcal{M}^{\gamma}]=\sum\limits_{l=1}^m (\gamma_l+1) \mathcal{M}^{\eta(l)}
 \end{equation}
 for $\eta(l)\in \mathbb{N}^m$, with
 \begin{equation}
 \gamma_j=\eta(l)_j, \quad j\in \{1,\cdots,l-1,l+1,\cdots,m\},\quad \gamma_l+1=\eta(l)_l,
 \end{equation}
 then $  \partial_{s}[\mathcal{M}^{\gamma}]$ can be regarded as the sum of
 \begin{equation}
 \sum\limits_{l=1}^{m}(\gamma_l+1)=m+ \sum\limits_{l=1}^{m}\gamma_l
 \end{equation}
 many terms with each term having the form of $\mathcal{M}^\eta$ with
\begin{equation}
\gamma_{j_0}+1=\eta_{j_0}, \quad \gamma_j=\eta_j, \quad j\in \{1,\cdots, m\} -\{j_0\}, \text{ for some }j_0\in \{1,\cdots,m\}.
\end{equation}
Then $\partial_{s_{j_1}}^{l_1}\cdots\partial_{s_{j_m}}^{l_m}[\mathcal{L}_m] $ can be regarded as the sum of no more than
\begin{equation}
\Pi_{j=0}^{l_1+\cdots+l_m-1}(m+j)
\end{equation}
many terms, with each term having the form of $\mathcal{M}^\eta$ with $|\eta|=l_1+\cdots+l_m.$  Since $m\leq k-1,$  therefore
 \begin{equation}
 \Pi_{j=0}^{l_1+\cdots+l_m-1}(m+j)\leq (2k)^{l_1+\cdots+l_m},
 \end{equation}
 we have
 \begin{equation}
 \sum\limits_{\gamma\in \mathbb{N}^m, |\gamma|=l_1+\cdots+l_m}|c_\gamma|\leq (2k)^{l_1+\cdots+l_m}
 \end{equation}
 and finish the proof.
 \end{proof}
 Then for type $1$, we do transformation in the following order: take the integral over $s_{j_l}$ for $l\leq m-1$, use Lemma \ref{key12.21.1} and condition \eqref{pp1}, use
 \begin{equation}
\sup\limits_{t\in \mathbb{R}}\frac{1}{a!}\sum\limits_{l,j=0}^2\sum\limits_{m,r=1}^3|\frac{\partial^a}{\partial t^a}\left[\partial_{\xi\cdot e_r}^l\partial_{\xi\cdot e_m}^j\hat{V}(\xi,t)\right]|\leq \frac{c^a \hat{V}_0(\xi)}{(1+|t|)^a} \text{ and } \frac{1}{(1+s+a)^{j}}\leq \frac{1}{(1+s)^j}, \text{ for }s,a,j>0,
 \end{equation}
take the integral over $\xi_1,\cdots,\xi_k, s_{j}$(such $s_j$ with $|\xi_j+q|>2M$) and we have
 \begin{align*}
 \| \text{type }1 \|_{\s^p_x}\leq &\int_0^\infty ds_{j_m}\cdots \int_0^\infty ds_{j_r}\int_0^\infty ds_k  (2k)^{r+3-m}\times\\
 &\frac{c^{r+3-m}\|\hat{V}_0(\xi)\|^k_{\s^1_\xi }}{(1+s_{j_m}+\cdots+s_{j_r}+s_k)^{r+3-m}}\frac{C_1^{r+2}}{(2\pi)^{3k/2}M^{r+2}}\times \frac{81^kC_2^{k-1-r}}{M^{k-1-r}}\|\psi(x)\|_{\s^p_x}.
 \end{align*}
 Since
\begin{align*}
&\int_0^\infty ds_{j_m}\cdots \int_0^\infty ds_{j_r}\int_0^\infty ds_k  (2k)^{r+3-m}\frac{1}{(1+s_{j_m}+\cdots+s_{j_r}+s_k)^{r+3-m}}\\
=&\frac{(2k)^{r+3-m}}{(r+2-m)!}\leq 2k e^{2k},
\end{align*}
we have
\begin{equation}
\| \text{type }1\|_{\s^p_x}\leq \frac{2c^{k+1}kC_3^{k+1}\|\hat{V}_0(\xi)\|_{\s^1_\xi\cap \s^\infty_\xi}^k}{M^{k+1}}\|\psi(x)\|_{\s^p_x}.
\end{equation}
\textbf{Estimates for $I_{\epsilon}^{(k)}\psi(x)$:} combining the estimates for type $1$, type $2$ and type $3$, we have
\begin{equation}
\| I_{\gamma\epsilon}^{(k)}\psi(x)\|_{\s^p_x}\lesssim \frac{c^{3k+2}k^3C_4^k\|\hat{V}_0(\xi)\|_{\s^1_\xi\cap \s^\infty_\xi}^{k}}{M^{k-1}}\|\psi(x)\|_{\s^p_x}.
\end{equation}
Hence,
\begin{equation}
\| I_\epsilon^{(k)}\psi(x)\|_{\s^p_x}\lesssim \frac{c^{3k+2}k^3C^k\|\hat{V}_0(\xi)\|_{\s^1_\xi\cap \s^\infty_\xi}^{k}}{M^{k-1}}\|\psi(x)\|_{\s^p_x}.
\end{equation}
Similarly,
\begin{equation}
\| \beta(|P|>32M)\left(I_\epsilon^{(k)}\right)^*\|_{\s^p_x\to \s^p_x}\lesssim \frac{c^{3k+2}k^3C^k\|\hat{V}_0(\xi)\|_{\s^1_\xi\cap \s^\infty_\xi}^{k}}{M^{k-1}}.
\end{equation}
\end{proof}
Now we can go to prove Theorem \ref{main3}.
\begin{proof}The proof is the same as Theorem \ref{intime} by applying Lemma \ref{exlnk}, Theorem \ref{3I1} instead.
\end{proof}
Similarly, we get asymptotic completeness on high frequency subspace.
\begin{corollary}If $V(x,t)$ satisfies the condition in Theorem \ref{main3}, the Schr\"odinger equation has asymptotic completeness on high frequency subspace.
\end{corollary}
   
Now let us think about
\begin{equation}
\Omega_T:= s\text{-}\lim\limits_{t\to \infty }U(T+t,T)e^{-itH_0}, \text{ on }\s^2, \text{ for }T\geq 0.\label{OmegaT}
\end{equation}
Assume
\begin{equation}
\Omega_T(t)=U(T+t,T)e^{-itH_0}.
\end{equation}
\begin{equation}
\Omega_{T,\epsilon}=I+(-i)\int_0^\infty dt e^{-\epsilon t}\Omega_T(t) e^{itH_0}V(x,t+T)e^{-itH_0}.
\end{equation}
By the same argument, we also have its $\s^p$ boundedness on high-frequency subspace:
\begin{corollary}\label{cor3}If $V(x,t)$ satisfies condition \eqref{pp1}, there exists $M=M(V(x,t))>0$ such that for all $1\leq p\leq\infty$, 
\begin{equation}
\Omega_T\beta(|H_0|>M^2)=s\text{-}\lim\limits_{\epsilon \downarrow0} \Omega_{T,\epsilon}\beta(|H_0|>M^2), \text{ on }\s^p,
\end{equation}
and $ \beta(|H_0|>M^2)\Omega_T^*, \Omega_T\beta(|H_0|>M^2)$ are bounded on $\s^p$.
\end{corollary}
\begin{proof}Since $\Omega_T$ is obtained by replacing $V(x,t)$ with $V(x,T+t)$ in $\Omega$ and since
\begin{equation}
\frac{(1+t)^a}{a!}\leq \frac{(1+t+T)^a}{a!}, \text{ for }t,T\geq 0,
\end{equation}
then following the same argument in Theorem \ref{main3}, the conclusion follows.
\end{proof}
Similarly, we have the following corollary:
\begin{corollary}If $V(x,t)$ satisfies the assumptions in Theorem \ref{main3}, and if $P_c$ is bounded on $\s^p$ for all $1\leq p\leq \infty$, there exists $M=M(V(x,t))>0$, such that
 \begin{equation}
 \sup\limits_{T\in \mathbb{R}}\| P_cU(T,0)e^{-iTH_0}\beta(|P|>M)\|_{\s_x^p\to \s_x^p}<C.
 \end{equation}
 \end{corollary}
This can be extended to the case when
\begin{equation}
V(x,t)=\chi(|t|<T_0)B(x,t)+\chi(|t|\geq T_0)V_1(x,t),
\end{equation}
with  $V_1(x,t) $ satisfying the assumption in Theorem \ref{3I1}, $\hat{B}(\xi,t)\in \s^\infty_t\s^1_\xi $. This application is based on the following operators
\begin{equation}
I_{\epsilon}^{(k)}(T_0):=\int_{T_0}^\infty dt_k \int_{t_k}^\infty dt_{k-1}\cdots \int_{t_2}^\infty e^{-\epsilon t_1}dt_1\sK_{t_k}(V(x,t_k))\cdots \sK_{t_1}(V(x,t_1))
\end{equation}
and
\begin{equation}
J_{\epsilon}^{(k)}(T_0):=\int_0^{T_0}dt_k \int_{t_k}^{T_0}dt_{k-1}\cdots \int_{t_2}^{T_0} e^{-\epsilon t_1}dt_1 \sK_{t_k}(V(x,t_k))\cdots \sK_{t_1}(V(x,t_1))
.
\end{equation}
Then
\begin{equation}
I_{\epsilon}^{(k)}=\sum\limits_{j=0}^k J_{\epsilon}^{(j)}(T_0)I_{\epsilon}^{(k-j)}(T_0).\label{relationshipeq}
\end{equation}

\begin{corollary}\label{8.24.1}If $V_1(x,t)$ satisfies the assumptions in Theorem \ref{main3}, $\hat{B}(\xi,t)\in \s^\infty_t \s^1_\xi $, then there exists some large $M$ such that for all $1\leq p\leq \infty$, $\Omega \beta(|P|>32M):\s^p_x\to \s^p_x$ is bounded.
\end{corollary}
\begin{proof} Similarly, we have that for $\psi\in \beta(|P|>32M)\mathcal{S}_x$,
\begin{align*}
I_\epsilon^{(k)}(T_0)\psi(x)=&\int_0^\infty e^{-\epsilon s_k}ds_k\cdots\int_{0}^\infty e^{-\epsilon s_1}ds_1\int
d^3\xi_1\cdots d^3\xi_kd^3q e^{i(x\cdot (\xi_k+q)+2(s_k\xi_k+\cdots+s_1\xi_1)\cdot q )}\\
\mathcal{V}(\xi,k)&e^{i(s_k\xi_k^2+\cdots +s_1\xi_1^2)}\frac{\hat{\psi}(q)}{(2\pi)^{\frac{3}{2}}}-\int_0^{T_0}e^{-\epsilon s_k}ds_k\cdots\int_{0}^\infty e^{-\epsilon s_1}ds_1\int
d^3\xi_1\cdots d^3\xi_kd^3q \mathcal{V}(\xi,k)\frac{\hat{\psi}(q)}{(2\pi)^{\frac{3}{2}}}\\
&e^{i(x\cdot (\xi_k+q)+(s_k\xi_k^2+\cdots +s_1\xi_1^2)+2(s_k\xi_k+\cdots+s_1\xi_1)\cdot q )}:=L_1^{(k)}\psi(x)+L_2^{(k)}\psi(x).
\end{align*}
We apply Lemma \ref{exlnk} to $L_1^{(k)}\psi(x)$ and have
\begin{equation}
\| L_1^{(k)}\psi(x)\|_{\s^p_x}\leq \frac{(2k)^3C_{V_1}^k}{\sqrt{M}^{k-1}}\|\psi(x)\|_{\s_x^p}, \text{ for some }C_{V_1}>0.
\end{equation}
For $L_2^{(k)}\psi(x)$, according to the proof of Lemma \ref{exlnk}, we do the same transformation for $\xi_j,s_j$, $j=1,\cdots,k-1$ while we do nothing for $s_k,\xi_k$. Similarly, in the end, we will get
\begin{equation}
\| L_2^{(k)}\psi(x)\|_{\s^p_x}\leq \frac{T_0(2k)^3D_{V_1}^k}{\sqrt{M}^{k-1}}\|\psi(x)\|_{\s_x^p},\text{ for some }D_{V_1}>0.
\end{equation}
Hence,
\begin{equation}
\| I_\epsilon^{(k)}(T_0)\psi(x)\|_{\s^p_x}\leq \frac{(2k)^3(1+T_0)(D_{V_1}+C_{V_1})^k}{\sqrt{M}^{k-1}}\|\psi(x)\|_{\s_x^p}.
\end{equation}
According to the same proof of Corollary \ref{exp1}, we have that for $\psi\in \s^q$,
\begin{equation}
\|J_{\epsilon}^{(k)}(T_0)\psi(x) \|_{\s^q_x}\leq \frac{T_0^k \|\hat{V}(\xi, t)\|_{\s^\infty_t\s^1_\xi }^k}{k!}\leq \frac{ T_0^k\|\hat{B}(\xi, t)\|_{\s^\infty_t\s^1_\xi }^k}{k!}\|\psi(x)\|_{\s^p_x}.
\end{equation}
Then for $\psi\in \beta(|P|>32M)\mathcal{S}_x$,
\begin{align*}
\|I_{\epsilon}^{(k)}\psi(x) \|_{\s^p_x}\leq &\sum\limits_{j=0}^k \frac{\mathcal{M}^j}{j!} \frac{(1+T_0)(2k-2j)^3\mathcal{M}^{k-j}}{\sqrt{M}^{k-j-1}}\|\psi(x)\|_{\s^p_x}\leq \frac{(1+T_0)(2k)^3\mathcal{M}^k}{\sqrt{M}^{k-1}}(\sum\limits_{j=0}^\infty \frac{\sqrt{M}^j}{j!})\|\psi(x)\|_{\s^p_x}\\\leq& \frac{(1+T_0)(2k)^3\mathcal{M}^k}{\sqrt{M}^{k-1}}\times\exp(\sqrt{M})\|\psi(x)\|_{\s^p_x},
\end{align*}
where
\begin{equation}
\mathcal{M}:=\max\left( T_0\|\hat{B}(\xi,t)\|_{\s^\infty_t\s^1_\xi }, D_{V_1}+C_{V_1}\right).
\end{equation}
Then choose $M$ large enough to make
\begin{equation}
\sum\limits_{k=1}^\infty \frac{k^3\mathcal{M}^k}{\sqrt{M}^{k-1}}<\infty
\end{equation}
and then we get the conclusion.
\end{proof}
\begin{corollary}\label{3Tu}If $V(x,t)$ satisfies the assumption in Theorem \ref{main3}, and if $P_c$ is bounded on $\s^p$ for all $1\leq p\leq \infty$, then when $M>0$ is sufficiently large,
\begin{equation}
\sup\limits_{T\in \mathbb{R}}\| P_cU(T,0)e^{-iTH_0}\beta(|P|>M)\|_{\s_x^p\to \s_x^p}<\infty, \text{ for }1\leq p\leq\infty.
\end{equation}
Therefore,
\begin{equation}
\sup\limits_{T\in \mathbb{R}}|T|^{3/2}\| P_cU(T,0)\beta(|P|>M)\|_{\s_x^p\to \s_x^{p^\prime}}<\infty, \text{ for }1\leq p\leq 2.
\end{equation}
\end{corollary}
\begin{proof}The proof is the same as that of Corollary \ref{norm1}.
\end{proof}

\subsection{Examples}
In this subsection, we consider the potential $V(x,t)$ satisfying
\begin{equation}
V(x,t)=\sum\limits_{j=0}^\infty V_j(x)\frac{1}{(1+t)^j}, \text{ for }t>\frac{T_0}{2}, \text{ for some }T_0>0.\label{8.24.c1}
\end{equation}
If
\begin{equation}
 \sum\limits_{b=0}^\infty \frac{2^b}{(1+T_0)^b}\sum\limits_{l,j=0}^2\sum\limits_{m,r=1}^3|\partial_{\xi\cdot e_r}^l\partial_{\xi\cdot e_m}^j\hat{V}_a(\xi)|\in \s^1_\xi \cap \s^\infty_\xi, \label{8.24.c2}
\end{equation}
and $\hat{V}(\xi, t)\in \s^\infty_t(0,T_0)\s^1_\xi $, then we choose $B(x,t)=\chi(t< T_0)V(x,t)$ and $V_1(x,t)=\chi(t\geq T_0)V(x,t)$ with
\begin{align*}
&\frac{(1+t)^a}{a!}\left|\sum\limits_{l,j=0}^2\sum\limits_{m,r=1}^3|\frac{\partial^a}{\partial t^a}\left[\partial_{\xi\cdot e_r}^l\partial_{\xi\cdot e_m}^j\hat{V}(\xi,t)\right]|\right|\leq  \sum\limits_{b=0}^\infty \frac{\begin{pmatrix}b+a-1\\a\end{pmatrix}}{(1+t)^{b}}\sum\limits_{l,j=0}^2\sum\limits_{m,r=1}^3|\partial_{\xi\cdot e_r}^l\partial_{\xi\cdot e_m}^j\hat{V}_a(\xi)|\\
\leq &2^a  \sum\limits_{b=0}^\infty \frac{2^b}{(1+T_0)^b}\sum\limits_{l,j=0}^2\sum\limits_{m,r=1}^3|\partial_{\xi\cdot e_r}^l\partial_{\xi\cdot e_m}^j\hat{V}_a(\xi)| .
\end{align*}
Then we can choose $c=2$ and
\begin{equation}
\hat{V}_0(\xi)= \sum\limits_{b=0}^\infty \frac{2^b}{(1+T_0)^b}\sum\limits_{l,j=0}^2\sum\limits_{m,r=1}^3|\partial_{\xi\cdot e_r}^l\partial_{\xi\cdot e_m}^j\hat{V}_a(\xi)| .
\end{equation}
Apply Corollary \ref{8.24.1} and we have the following corollary:
\begin{corollary}\label{ext1}Assume $V(x,t)$ has the form of \eqref{8.24.c1} and satisfies condition \eqref{8.24.c2},  then $\Omega \beta(|P|> M): \s_x^p\to \s_x^p$ is bounded for some sufficiently large $M$.
\end{corollary}
Now we are considering the potential $V(x,t)$ satisfying
\begin{equation}
V(x,t)=\sum\limits_{j=0}^\infty V_j(x)f_j(t),\label{8.24.c3}
\end{equation}
when $t>\frac{T_0}{2}$ for some $T_0>0$. If $\hat{V}(\xi,t)\in \s^\infty_t(0,T_0/2)\s^1_\xi $ and if for $b\in \mathbb{N}$,
\begin{equation}
\sup\limits_{t\in [T_0/2,\infty)}\frac{(t+1)^b}{b!}|f^{(b)}_j(t)|\leq c_j^b, \quad \sum\limits_{a=0}^\infty c_a^b \sum\limits_{l,j=0}^2\sum\limits_{m,r=1}^3|\partial_{\xi\cdot e_r}^l\partial_{\xi\cdot e_m}^j\hat{V}_a(\xi)| <\infty,\label{12.28.1}
\end{equation}
we will get a similar result:
\begin{corollary}\label{timesum1}Assume $V(x,t)$ has the form of \eqref{8.24.c3} and satisfies condition \eqref{12.28.1}, then $\Omega \beta(|P|> M): \s_x^p\to \s_x^p$ is bounded for some sufficiently large $M$.
\end{corollary}

Here are some other examples.
\begin{example}[quench potentials] A quench potential has the form of $V(x,t)=\chi(t\geq d)V_1(x)$ or $V(x,t)=\beta(t> 2d)V_1(x)$ for some $d>0$. If $\sum\limits_{l,j=0}^2\sum\limits_{m,r=1}^3|\partial_{\xi\cdot e_r}^l\partial_{\xi\cdot e_m}^j\hat{V}_1(\xi)|\in \s^1_\xi\cap \s^\infty_\xi$, then $\Omega\beta(|P|>M)$ is bounded on $\s_x^p$ for some sufficiently large $M$.
\end{example}
\begin{proof} Choose $B(x,t)=V(x,t), $ $T_0=d, c=1$, $V_1(x,t)=V_1(x)$. When we take the derivative with respect to $t$, it is $0$ and of course satisfies the condition \eqref{pp1}.
\end{proof}
\begin{example}[Hyperbolic potentials] A hyperbolic potential has the form of $V(x,t)=\tanh(t)V_1(x)+V_0(x)$. If $\sum\limits_{l,j=0}^2\sum\limits_{m,r=1}^3|\partial_{\xi\cdot e_r}^l\partial_{\xi\cdot e_m}^j\hat{V}_a(\xi)|\in \s^1_\xi\cap \s^\infty_\xi$, $a=0,1$, then $\Omega\beta(|P|>M)$ is bounded on $\s_x^p$ for some sufficiently large $M$.
\end{example}
\begin{proof} Since for $a\in \mathbb{N}^+$, $t\geq 1$,
\begin{equation}
\frac{(1+t)^j}{j!}\frac{d^j}{dt^j}[\tanh t]=\frac{(1+t)^j}{j!}\frac{d^j}{dt^j}[1-2e^{-2t}\sum\limits_{l=0}^\infty (-1)^le^{-2lt}]=-\sum\limits_{l=0}^\infty(-1)^l\frac{[-2(l+1)(1+t)]^j}{j!}e^{-2(l+1)t},
\end{equation}
we can choose $c=4$ and
\begin{equation}
\sup\limits_{t\in [1,\infty)}\frac{(1+t)^j}{j!}|\frac{d^j}{dt^j}[\tanh t]|\leq 4^j\sum\limits_{l=0}^\infty e^{-(l+1)t}=\frac{4^je^{-t}}{1-e^{-t}}<4^j.
\end{equation}
For $t\in [0,1)$, it satisfies the condition for some time. By Corollary \ref{timesum1}, we get the result.
\end{proof}

\section{Moving and self-similar potentials}\label{section 4}
A fundamental class of time dependent potentials is moving potentials, of the form $\sum_i  V_i(x-c_i(t))$. They appear naturally in charge transfer models, soliton dynamics, models of Atom$+$Radiation and more. The mathematical analysis of such potentials has been carried out for certain classes, mostly when
\begin{equation}
c_i(t)=ct+f(t)
\end{equation}
with $f(t)$ decaying fast, see \cite{RSS2005} and \cite{P2004}. More general movement was considered in \cite{BS2011},\cite{BS2012} and \cite{BS2019}, but it was limited to ONE potential term. Moreover it was assumed that the velocity goes to zero, or random in other cases. The more difficult cases when the movement is not linear is treated in this section. But the case $c(t)=t$ does not satisfy our condition, if there is another potential added. For more information about charge transfer models, see \cite{Chen2016}, \cite{Cai2003} and \cite{CL1999}.\par
After the posting of this work, further works appeared  in this direction. In particular, the Local Decay and Strichartz Estimates for Charge transfer Hamiltonians were established in the work \cite{Chen-J}. In \cite{S-Wu} a proof of the existence and completeness of the wave operators for Charge Transfer Hamiltonians with general time dependent potentials (and moving with constant velocities)
was proved by a new approach to the scattering problem. But this approach did not use or yield Local Decay estimates.

We prove Theorem \ref{dilationg}(the self-similar example) first.
\begin{theorem}If $V(x,t)$ is defined in equation \eqref{type21} and satisfies condition \eqref{710.3}, then
\begin{equation}
\lim\limits_{T\to \pm \infty}\| U(0,T)e^{-iTH_0}-\Omega \|_{\s^p\to \s^p}= 0, \quad \|\Omega\|_{\s^p\to \s^p}\leq \exp\left(\frac{\| h(t)\|_{\s^1_t(0,\infty)}}{(2\pi)^{\frac{3}{2}}}\right).
\end{equation}
\end{theorem}
\begin{proof}In this case, since
\begin{equation}
\sK_t(V(x,t))=\frac{1}{(2\pi)^{\frac{3}{2}}}\int d^3\xi \hat{V}_1(\xi,t)e^{iH_0t}e^{i\xi \cdot g(t)x}e^{-iH_0t}+\sum\limits_{j=1}^\infty f_j(t) e^{iH_0t}e^{ia_j \cdot g_j(t)x}e^{-iH_0t}.
\end{equation}
According to the same computation in section 1 and the proof of Corollary \ref{exp1}, we have that for $T_0\in [0,\infty]$,
\begin{equation}
\|\sum\limits_{k=0}^\infty i^kI(T_0)^{(k)}\|_{\s_x^p\to \s_x^p}\leq \exp\left(\frac{\| h(t)\|_{\s^1_t[0,\infty)}}{ (2\pi)^{\frac{3}{2}}}\right)
\end{equation}
where
\begin{equation}
I(T_0)^{(j)}:=\int_0^{T_0} dt_1\int_{0}^{t_{1}} dt_{2}\cdots\int_0^{t_{j-1}}dt_j\sK_{t_j}(V(x,t_j)) \cdots \sK_{t_1}(V(x,t_1)) ,
\end{equation}
and as $T_0\to \infty$,
\begin{equation}
\|\sum\limits_{k=0}^\infty I(\infty)^{(k)}-\sum\limits_{k=0}^\infty I(T_0)^{(k)}\|_{\s_x^p\to \s_x^p}\leq \frac{\| h(t)\|_{\s^1_t[T_0,\infty)}}{(2\pi)^{\frac{3}{2}}}\times \exp\left(\frac{\| h(t)\|_{\s^1_t[0,\infty)}}{ (2\pi)^{\frac{3}{2}}}\right)\to 0.
\end{equation}
Then $\sum\limits_{k=0}^\infty I(T_0)^{(k)}\to \sum\limits_{k=0}^\infty I(\infty)^{(k)}$ in norm. Then
\begin{equation}
\Omega=\sum\limits_{k=0}^\infty I(\infty)^{(k)},\quad \|\Omega\|_{\s_x^p\to \s_x^p}\leq \exp\left(\frac{\| h(t)\|_{\s^1[0,\infty)}}{(2\pi)^{\frac{n}{2}}}\right).
\end{equation}
\end{proof}
\begin{corollary}\label{4Tu}If $V(x,t)$ satisfies the assumption in Theorem \ref{dilationg}, then
\begin{equation}
\sup\limits_{T\in \mathbb{R}}\| U(0,T)e^{-iTH_0}\|_{\s_x^p\to \s_x^p}<\infty, \text{ for }1\leq p\leq\infty.
\end{equation}
Therefore,
\begin{equation}
\sup\limits_{T\in \mathbb{R}}|T|^{3/2}\| U(T,0)\|_{\s_x^p\to \s_x^{p^\prime}}<\infty, \text{ for }1\leq p\leq 2.
\end{equation}
\end{corollary}
\begin{proof}The proof is the same as that of Theorem \ref{dilationg}.
\end{proof}
Here is an example where $f(t)$ does not even have a limit in $\mathbb{R}^3$ as $t\to \pm\infty$ and it is not just limited to one potential:
\begin{example}Assume a potential has the form of $V(x,t)=V_1(x-\sin(\ln(1+|t|))v)+V_0(x)$ for some $v\in \mathbb{R}^3$. Then if $\sum\limits_{l,j=0}^2\sum\limits_{m,r=1}^3|\partial_{\xi\cdot e_r}^l\partial_{\xi\cdot e_m}^j\hat{V}_a(\xi)|\in \s^1_\xi\cap \s^\infty_\xi$, $a=0,1$, and the support of $\hat{V}_1$ is contained in a ball $B_R$ centered at the origin with a radius $R$,  then $\Omega\beta(|P|>M)$ is bounded on $\s_x^p$ for some sufficiently large $M$.
\end{example}
\begin{proof} In this case,
\begin{equation}
\hat{V}(\xi,t)=\hat{V}_0(\xi)+\hat{V}_1(\xi)e^{-\sin(\ln(1+|t|))i\xi\cdot v}.
\end{equation}
For $t\geq 0$, $a\in \mathbb{N}^+$,
\begin{equation}
\left|\partial_t^a  \partial_{\xi\cdot e_r}^l\partial_{\xi\cdot e_m}^j[\hat{V}(\xi,t)]\right|\leq \sum\limits_{b=0}^4(R|v|)^b\left|\partial_t^a[\sin(\ln(1+t))^{b}e^{-\sin(\ln(1+t))i\xi\cdot v}]\right| \sum\limits_{l,j=0}^2\sum\limits_{m,r=1}^3|\partial_{\xi\cdot e_r}^l\partial_{\xi\cdot e_m}^j\hat{V}_1(\xi)|.
\end{equation}
Since for $a_1,a_2,a_3\in \mathbb{R}$,
\begin{align*}
\frac{d}{dt}&[e^{(a_1i-a_2)\ln(1+t)-i\sin(\ln(1+t))a_3}]=(a_1i-a_2)e^{(a_1i-a_2-1)\ln(1+t)-i\sin(\ln(1+t))a_3}\\
&-\frac{i}{2}e^{((a_1+1)i-a_2-1)\ln(1+t)-i\sin(\ln(1+t))a_3}-\frac{i}{2}e^{((a_1-1)i-a_2-1)\ln(1+t)-i\sin(\ln(1+t))a_3},
\end{align*}
we can regard it as the sum of $|a_1|+|a_2|+1$ many terms with each term having the form of
\[
\pm e^{(b_1i-a_2-1)\ln(1+t)-i\sin(\ln(1+t))a_3}, \pm ie^{(b_1i-a_2-1)\ln(1+t)-i\sin(\ln(1+t))a_3}
\]
with $|b_1-a_1|=0$ or $|b_1-a_1|=1$. Hence, for $b\in \{-4,-3,\cdots,3,4\}$,
\begin{equation}
\frac{(1+t)^a}{a!}\left|\frac{d^a}{dt^a} [e^{bi\ln(1+t)-\sin(\ln(1+t))iv\cdot \xi}]\right|\leq \frac{1}{a!}\Pi_{j=0}^{a-1}(|b|+1+2j)\leq 2^{2a+3}.
\end{equation}
Then there exists a constant $C$ independent on $a$ such that
\begin{equation}
\left|\partial_t^a  \partial_{\xi\cdot e_r}^l\partial_{\xi\cdot e_m}^j[\hat{V}(\xi,t)]\right|\leq \sum\limits_{b=0}^4(R|v|)^b\left|\partial_t^a[\sin(\ln(1+t))^{b}e^{-\sin(\ln(1+t))i\xi\cdot v}]\right| \leq C \frac{(4|v|R)^a}{(1+t)^a}
\end{equation}
which implies $V(x,t)$ satisfies condition \eqref{pp1} and finish the proof.

\end{proof}
In the following, we apply the same argument as in previous sections, to prove decay estimates for potentials $V(x-\sqrt{1+|t|}v)$ on high frequency subspace for $v\in \mathbb{R}^3,$ which satisfies assumption \ref{sqrtc}.
\begin{remark}Here $\sqrt{1+|t|}$ is crucial since $\sqrt{1+|t|}$ is not Mikhlin-type anymore, and the derivative of $\sqrt{1+t}(t>0)$ is not in $\s^2_t(0,\infty)$.
\end{remark}
We stick to $t>0$. Let
\begin{equation}
\mathcal{G}_{\leq 2M}(\eta, t):=\beta(|P|\leq 2M)e^{itH_0}e^{i\eta\cdot x}e^{-itH_0}, \text{ for }\eta\in \mathbb{R}^3,
\end{equation}
\begin{equation}
\mathcal{G}_{> 2M}(\eta, t):=\beta(|P|> 2M)e^{itH_0}e^{i\eta\cdot x}e^{-itH_0},
\end{equation}
\begin{align}
\mathscr{G}_M(\xi^j,&t_{k+j+1}, s^j,k)= \mathcal{G}_{\leq 2M}(\xi_{k+j}-\xi_{k+j-1}, t_{k+j+1}+s_{k+j})\times\\
\mathcal{G}_{> 2M}(\xi_{k+j-1}-\xi_{k+j-2},& t_{k+j+1}+\sum\limits_{l=k-1}^{k}s_{l+j})\cdots\mathcal{G}_{> 2M}(\xi_{1+j}-\xi_{j}, t_{k+j+1}+\sum\limits_{l=1}^{k}s_{l+j}),
\end{align}
for $\xi\in \mathbb{R}^{3(k+j)}, s\in \mathbb{R}^{k+j}, t_{k+j+1}\in \mathbb{R}$, $j\in\mathbb{N}$, with $\xi_0=0$,
\begin{equation}
 \mathcal{V}(\xi,t_{k+1},s,k):=\Pi_{j=1}^k \mathscr{F}[V(x-\left(\sqrt{1+t_{k+1}+\sum\limits_{l=j}^{k}s_l}\right)v)](\xi_j-\xi_{j-1})
\end{equation}
and let
\begin{align}
&\mathscr{J}^{(k+1)}_{M,\epsilon}:=\frac{1}{(2\pi)^{3k/2}}\int d^3\xi_1\cdots d^3\xi_k\int_0^{\infty}dt_{k+1}e^{-\epsilon t_{k+1}}U(0,t_{k+1})e^{i\xi_k\cdot x}V(x-\sqrt{1+t_{k+1}}v)\times\\
&e^{-it_{k+1}H_0}\int_0^\infty e^{-\epsilon s_k}ds_k\int_{0}^\infty e^{-\epsilon s_{k-1}}ds_{k-1}\cdots \int_{0}^\infty e^{-\epsilon s_1}ds_1 \mathcal{V}(\xi,t_{k+1},s,k)\mathscr{G}_M(\xi^0,t_{k+1},s^0,k),
\end{align}
\begin{align}
\mathscr{K}^{(k)}(T):=&\int_0^Tdt_1\int_0^{t_1}dt_2\cdots \int_0^{t_{k-1}}dt_k e^{it_kH_0}V(x-\sqrt{1+|t_k|}v)e^{-it_kH_0}\beta(|P|>2M)\cdots \\
&e^{it_1H_0}V(x-\sqrt{1+|t_1|}v)e^{-it_1H_0}\beta(|P|>2M).
\end{align}
Its proof is based on following lemma and the estimates for $\mathscr{J}_{M,\epsilon}^{(k+1)}, \mathscr{K}^{(k)}(T)$:
\begin{lemma}[Representation formula 2]\label{refor2}For $\xi_i\in \mathbb{R}^n$, $i=1,\cdots, k$ $(k \in\mathbb{N}^+)$, $\psi(x)\in \s_x^p(\mathbb{R}^n)$, we have
\begin{align*}
&\mathcal{G}_{\leq 2M}(\xi_{k}-\xi_{k-1},t_{k})\mathcal{G}_{> 2M}(\xi_{k-1}-\xi_{k-2},t_{k-1})\cdots\mathcal{G}_{> 2M}(\xi_{k-1}-\xi_{k-2},t_{k-1}) \psi(x)\\
=&\frac{1}{(2\pi)^{\frac{n}{2}}}\int d^nq e^{i(x\cdot (\xi_k+q)+t_k(\xi_k^2+2q\cdot \xi_k)+(t_{k-1}-t_k)(\xi_{k-1}^2+2q\cdot \xi_{k-1})+\cdots +(t_1-t_2)(\xi_1^2+2\xi_1\cdot q))}\times\\
&\beta(|\xi_{k}+q|\leq 2M)\Pi_{j=1}^{k-1}\beta(|\xi_{j}+q|> 2M)\hat{\psi}(q).
\end{align*}
\end{lemma}
\begin{proof}It follows directly from
\begin{equation}
f(|P|)e^{ix\cdot\xi}=e^{ix\cdot \xi}f(|P+\xi|)
\end{equation}
and Lemma \ref{refor}.
\end{proof}
\begin{lemma}\label{flowJ}If $V(x-\sqrt{1+t}v)$ satisfies assumption \eqref{sqrtc}, then when $M$ is large enough,
\begin{equation}
\sup\limits_{T\in \mathbb{R}}|T|^{3/2} \| \mathscr{J}_{M,\epsilon}^{(k+1)}e^{iTH_0}\|_{\s^{1}_x\to \s^\infty_x}\leq \frac{k^5(C|||V(x)|||_p)^k}{\sqrt{M}^{k}}
\end{equation}
for some constant $C$.
\end{lemma}
\begin{proof}Due to Lemma \ref{refor2}, for $s_k,\xi_k$, we have a factor $\beta(|\xi_k+P|\leq 2M)$. We deal with them first. \textbf{Step one:} in this case, we have to face
\begin{equation}
\int_0^\infty ds_k e^{is_k(\xi_k^2+2\xi_k\cdot q)-\epsilon s_k}\left(\Pi_{l=1}^k e^{-i\sqrt{1+\sum\limits_{l=j}^{k+1}s_l}v\cdot (\xi_j-\xi_{j-1})}\right).
\end{equation}
We do the same transformation as before, that is,
\[
e^{is_k(\xi_k^2+2\xi_k\cdot q)-\epsilon s_k}=\frac{1}{i(\xi_k^2+2\xi_k\cdot q)-\epsilon}\partial_{s_k}[e^{is_k(\xi_k^2+2\xi_k\cdot q)-\epsilon s_k} ].
\]
Then we will get two terms: boundary term
\[
\frac{1}{i(\xi_k^2+2\xi_k\cdot q)-\epsilon}e^{is_k(\xi_k^2+2\xi_k\cdot q)-\epsilon s_k}\Pi_{l=1}^k e^{-i\sqrt{1+\sum\limits_{l=j}^{k+1}s_l}v\cdot (\xi_j-\xi_{j-1})}
 \vert_{s=0}
\]
and the integral term
\[
\frac{1}{i(\xi_k^2+2\xi_k\cdot q)}\int_0^\infty ds_ke^{is_k(\xi_k^2+2\xi_k\cdot q)-\epsilon s_k}\partial_{s_k}\left(\Pi_{l=1}^k e^{-i\sqrt{1+\sum\limits_{l=j}^{k+1}s_l}v\cdot (\xi_j-\xi_{j-1})}\right).
\]
For the integral term, we keep doing this transformation until we reach $\partial_{s_k}^5$($\partial_{s_k}^{5}$ will bring no more than $(2k)^5$ many terms with each term controlled by $1/(1+s_k+t_{k+1})^{5}$). \textbf{Step two: }we keep doing transformation for the boundary terms. For each boundary term, we break it into two terms ($\mathcal{G}_{\leq 2M}(\xi_{k+1}-\xi_{k}, t_{k+1})$ and $\mathcal{G}_{> 2M}(\xi_{k+1}-\xi_{k}, t_{k+1})$). \textbf{Step three: }for the term with $\mathcal{G}_{> 2M}(\xi_{k+1}-\xi_{k}, t_{k+1})$, we keep using Duhamel's formula
\begin{equation}
\mathbbm{1}\cdots+i\int_0^\infty dt_{k+2}U(t_{k+2},0)\cdots.
\end{equation}
For the $\mathbbm{1}$ term, it has the same form as $I_\epsilon^{(k+1)}e^{iTH_0}$. For the integral term, we break it into two terms ($\mathcal{G}_{\leq 2M}(\xi_{k+2}-\xi_{k+1}, t_{k+2})$ and $\mathcal{G}_{> 2M}(\xi_{k+2}-\xi_{k+1}, t_{k+2})$). We keep doing this until we gain $\mathcal{G}_{\leq 2M}(\xi_{k+j}-\xi_{k+j-1},t_{k+j})$ for some $j\in \mathbb{N}^+$(\textbf{type one}) or there is no $U(0,t_{k+j})$(\textbf{type two}) in it. \textbf{Step four:} for the term with $\mathcal{G}_{\leq 2M}(\xi_{k+j}-\xi_{k+j-1}, t_{k+j})$, we use Duhamel's formula one more time. Then for the integral term, after changes of variables $t_{k+l}=t_{k+j+1}+\sum\limits_{m=k+l}^{k+j}s_m$, $l=1,\cdots,j$, we get
\begin{align*}
&\int d^3\xi_1\cdots d^3\xi_{k+j}\int_0^{\infty}dt_{k+j+1}e^{-\epsilon t_{k+j+1}}U(t_{k+j+1},0)e^{i\xi_{k+j+1}\cdot x}V(x-\sqrt{1+t_{k+j+1}}v)e^{-it_{k+j+1}H_0}\\
&\int_0^\infty e^{-\epsilon s_{k+j}}ds_{k+j}\cdots \int_{0}^\infty e^{-\epsilon s_1}ds_1\delta(s_k)\partial_{s_k}^{b_k}[\mathcal{V}(\xi,t_{k+j+1},s,k+j)]\mathscr{G}_M(\xi,t_{k+j+1},s,k,j)\frac{1}{(2\pi)^{3(k+j)/2}}\\
&\times (-1)^{b_k+1}/(i(\xi^2_k+2\xi_k\cdot P)^{b_k+1})\beta(|P|>32M),
\end{align*}
for some $b_k\in\{0,1,2,3,4\}$, where
\begin{equation}
\mathscr{G}_M(\xi,t_{k+j+1},s,k,j):=\mathscr{G}_{M}(\xi^k,t_{k+j+1},s^k,j)\mathscr{G}_M(\xi^0,t_{k+j+1}+\sum\limits_{l=k+1}^{k+j}s_{l},s^0,k).
\end{equation}
Then for $\xi_{k+j},s_{k+j}$, we do the same transformation as $\xi_k,s_k$ except that for $\xi_{k+j},s_{k+j}$, we stop integration by parts until we gain $\partial_{s_{k+j}}^{b_{k+j}}$ with $b_{k+j}=5-b_k$. For the boundary terms, we do the same transformation as step two to step four except that we stop until we gain $\partial_{s_{k+j_1}}^{b_{k+j_1}}\cdots \partial_{s_{k+j_l}}^{b_{k+j_l}}$ with $b_{k+j_1}+\cdots+b_{k+j_l}=5$. After these transformations, we will get many terms having  the following form: \\
\textbf{case one:}
\begin{align*}
&\int d^3\xi_1\cdots d^3\xi_{k+|j|}\int_0^{\infty}dt_{k+|j|+1}e^{-\epsilon t_{k+|j|+1}}U(t_{k+|j|+1},0)e^{i\xi_{k+|j|+1}\cdot x}V(x-\sqrt{1+t_{k+|j|+1}}v)e^{-it_{k+|j|+1}H_0}\\
&\int_0^\infty e^{-\epsilon s_{k+|j|}}ds_{k+|j|}\cdots \int_{0}^\infty e^{-\epsilon s_1}ds_1\delta(s_k)\delta(s_k+j_1)\cdots \delta(s_k+j_1+\cdots+j_{l-1}) \\
&\partial_{s_k}^{b_k}\cdots\partial_{s_k+j_1+\cdots+j_{l}}^{b_{k+l}}[\mathcal{V}(\xi,t_{k+|j|+1},s,k+|j|)]\mathscr{G}_M(\xi,t_{k+|j|+1},s,k,j,l)\frac{1}{(2\pi)^{3(k+|j|)/2}}\beta(|P|>32M)\\
&\times (-1)^{l+b_k+\cdots+b_l}\times1/(i(\xi_{k+|j|}^2+2\xi_{k+|j|}\cdot P))^{b_{k+l}}\times\Pi_{m=0}^{l-1}1/(i(\xi_{k+j_1+\cdots j_{m}}^2+2\xi_{k+j_1+\cdots j_{m}}\cdot P))^{b_{k+m}+1}
\end{align*}
for $b_k+\cdots+b_{k+l}=5$, $b_{k+m}\in\mathbb{N}, m=0,\cdots,m$, where $j=(j_1,\cdots,j_l)\in \mathbb{N}^l$,
\begin{align}
&\mathscr{G}_M(\xi,t_{k+|j|+1},s,k,j,l):=\mathscr{G}_{M}(\xi^{k+j_1+\cdots+j_{l-1}},t_{k+|j|+1},s^{k+j_1+\cdots+j_{l-1}},j_l)\times\cdots\\
&\mathscr{G}_{M}(\xi^{k},t_{k+|j|+1}+\sum\limits_{m=k+j_1+1}^{k+|j|}s_m,s^{k},j_1)\mathscr{G}_M(\xi^0,t_{k+|j|+1}+\sum\limits_{l=k+1}^{k+|j|}s_{l},s^0,k);
\end{align}
\textbf{case two:}
\begin{align*}
&\int d^3\xi_1\cdots d^3\xi_{k+|j|}\int_0^{\infty}dt_{k+|j|+1}e^{-\epsilon t_{k+|j|+1}}e^{i(k+|j|+1)H_0}e^{i\xi_{k+|j|+1}\cdot x}V(x-\sqrt{1+t_{k+|j|+1}}v)e^{-it_{k+|j|+1}H_0}\\
&\int_0^\infty e^{-\epsilon s_{k+|j|}}ds_{k+|j|}\cdots \int_{0}^\infty e^{-\epsilon s_1}ds_1\delta(s_k)\delta(s_k+j_1)\cdots \delta(s_k+j_1+\cdots+j_{l}) \\
&\partial_{s_k}^{b_k}\cdots\partial_{s_k+j_1+\cdots+j_{l}}^{b_{k+l}}[\mathcal{V}(\xi,t_{k+|j|+1},s,k+|j|)]\mathscr{G}_M(\xi,t_{k+|j|+1},s,k,j,l)\frac{1}{(2\pi)^{3(k+|j|)/2}}\\
&\times\Pi_{m=0}^{l}1/(i(\xi_{k+j_1+\cdots j_{m}}^2+2\xi_{k+j_1+\cdots j_{m}}\cdot P))^{b_{k+m}+1}\beta(|P|>32M)
\end{align*}
for $b_k+\cdots+b_{k+l}\leq 4$.\par
Now we deal with $\xi_j,s_j$ with $\beta(|\xi_j+q|> 2M)$. In this case, we do the same transformation as before except that for $s_j\geq 1/\sqrt{M}$, after taking integration by parts in $\xi_{j,l}:=\xi_j\cdot e_l$ for some direction $e_l$, we may gain
\begin{equation}
\frac{iv_l\xi_{j,l}(\sqrt{1+t_{k+1}+s_k+\cdots+s_{j+1}}-\sqrt{1+t_{k+1}+s_k+\cdots+s_{j}}) }{s_j}e^{-i\sqrt{1+\sum\limits_{l=j}^{k+1}s_l}v\cdot (\xi_j-\xi_{j-1}) }
\end{equation}
which means for some terms, we can only gain
\[
\frac{1}{\sqrt{1+t_{k+1}+s_k+\cdots+s_{j+1}}+\sqrt{1+t_{k+1}+s_k+\cdots+s_j}}
\]
since we have
\begin{equation}
\mathscr{F}[V(x-\sqrt{1+s}v)](\xi)=\hat{V}(\xi)e^{-i\sqrt{1+s}v\cdot \xi}.
\end{equation}
For these terms, we keep doing the same transformation until we gain
\[
\frac{1}{(\sqrt{1+t_{k+1}+s_k+\cdots+s_{j+1}}+\sqrt{1+t_{k+1}+s_k+\cdots+s_j})^{a}}\times \frac{1}{s_j^b} \text{ for }a/2+b>1,\text{ for some } a, b\in\mathbb{N}
\]
which means we do this transformation for no more than $3$ times. In the end, we deal with $t_{k+|j|+1}$. For case two, we have to face
\begin{equation}
D(T):=\int_0^\infty dt_{k+|j|+1} e^{it_{k+|j|+1}H_0}e^{ix\cdot \xi_{k+|j|}}V(x-\sqrt{1+t_{k+|j|+1}}v)e^{i(T-t_{k+|j|+1})H_0}
\end{equation}
since due to Lemma \ref{refor2}, other parts are reduced to be translation. We need following lemma:
\begin{lemma}\label{D(T)}If $\hat{V}(\xi)\in \s^1_\xi$ and $V(x)\in \s^1_x$, then
\begin{equation}
\sup\limits_{T\in \mathbb{R}} |T|^{3/2}\| D(T)\|_{\s_x^1\to \s_x^\infty} \leq C(\|\hat{V}(\xi)\|_{\s^1_\xi}+\|V(x) \|_{\s^1_x}), \text{ for some }C>0.
\end{equation}
\end{lemma}
\begin{proof}For $t_{k+|j|+1}\in (1, T-1)\cup (T+1,\infty)$, we use
\begin{equation}
\|e^{it_{k+|j|+1}H_0}e^{ix\cdot \xi_{k+|j|}}V(x-\sqrt{1+t_{k+|j|+1}}v)e^{i(T-t_{k+|j|+1})H_0}\|_{\s^1_x\to \s^\infty_x}\leq \frac{\|V(x)\|_{\s^1_x} }{|t_{k+|j|+1}|^{3/2}|T-t_{k+|j|+1}|^{3/2}}
\end{equation}
while for $t_{k+|j|+1}\in (0,1]\cup [T-1,T+1]$, we use cancellation lemma \ref{cancell}. Then the result follows.
\end{proof}
After all these transformations, based on Lemma \ref{D(T)}, we will gain no more than $C_1^{k+|j|}$ many terms for some $C_1$. Then for each term, we will gain at least $C_1^{k+|j|}||| V(x)|||_p^{k+|j|+1}/\sqrt{M}^{k+|j|} $. Hence, we have
\begin{equation}
\sup\limits_{T\in \mathbb{R}}|T|^{3/2}\|\text{case two}^{(k+|j|+1)}(T)\|_{\s_x^1\to \s_x^\infty}\leq \frac{(k+|j|)^4C^{k+|j|+1}|||V(x)|||_p^{k+|j|+1}}{\sqrt{M}^{k+|j|}}
\end{equation}
where $(k+|j|)^4$ comes from that for $a:=b_k+\cdots+b_{k+l}\leq 4$,
\begin{align*}
&|  \partial_{s_k}^{b_k}\cdots\partial_{s_k+j_1+\cdots+j_{l}}^{b_{k+l}}[\Pi_{m=1}^{k+|j|} e^{-i(t_{k+|j|+1}+s_{k+|j|}+\cdots+s_m)(\xi_m-\xi_{m-1})\cdot v}]|\\
\leq &\frac{C(k+|j|)^4 \max\limits_{j=1\cdots,k+|j|}(|\xi_j-\xi_{j-1}|+1)^4}{(1+t_{k+|j|+1})^a}, \text{ for some }C>0.
\end{align*}
For case one, we need following lemma:
\begin{lemma}\label{case1.1}If $V\in \s^\infty_t\s^1_x\cap \s^2_x $ and $\hat{V}(\xi,t)\in \s^\infty_t\s^1_\xi $,  then
\begin{equation}
B:=\sup\limits_{|s-t|\geq 1}\|U(s,t)\|_{\s_x^1\to \s_x^\infty}<\infty.
\end{equation}
\end{lemma}
\begin{proof}By using Duhamel's formula twice,
\begin{align*}
&U(s,t)=e^{-i(t-s)H_0}+(-i)\int_0^{t-s} du e^{-i[(t-s)-u]H_0}V(x,s+u)e^{-iuH_0}-\int_0^{t-s}du\int_0^{u}dwe^{-i[(t-s)-u]H_0}\times\\
&V(x,s+u)U(s+w,s+u)V(x,s+w)e^{-iwH_0}=: A_1+A_2+\int_0^{t-s}du\int_0^udw A_3(u,w,s,t).
\end{align*}
For the first two terms, it is clear when $V(x,t)\in \s^\infty_t\s^1_x$ and $\hat{V}(\xi,t)\in \s^\infty_t\s^1_\xi $. For the last term, when $u\leq 1$, we use
\begin{equation}
\sup\limits_{|a|\leq 1}\| U(s+w, a+s+w)e^{iaH_0}\|_{\s_x^\infty\to \s_x^\infty}<C, \text{ for some constant }C.
\end{equation}
So in the following, we stick to $u\geq 1$. When there is no singularity, since $U(s,t)$ is unitary on $\s_x^2$, we have
\begin{equation}
\| A_3(u,w,s,t)\|_{\s^1_x\to \s^\infty_x}\leq \frac{\|V(x,t)\|_{\s^2_x \s^\infty_t}}{|w|^{3/2}|t-s-u|^{3/2}}
\end{equation}
and then it is integrable over $\int_0^{t-s}du \int_{0}^u dw$ when there is no singularity. When there is a singularity for $1/w$, we use
\begin{equation}
U(s+w, s+u)V(x,s+w)e^{-iwH_0}=U(s+w, s+u+w) [U(s+u+w,s+u)e^{-iwH_0}] [e^{iwH_0}V(x,s+w)e^{-iwH_0}].
\end{equation}
Since Corollary \ref{exp1} tells us $ U(s+u+w,s+u)e^{-iwH_0}: \s_x^p\to \s_x^p$, is bounded by $e$ if $w$ is small enough, we have
\begin{equation}
\| A_3(u,w,s,t)\|_{\s^1_x\to \s^\infty_x}\leq  \frac{C_1(B+1)\|\hat{V}(\xi,t)\|_{\s^\infty_t\s^1_\xi }\|V(x,t)\|_{\s^\infty_t\s^1_x}}{|t-s-u|^{3/2}}
\end{equation}
for some constant $C$, where we use
\begin{equation}
\|U(s+w, s+u+w)  \|_{\s^1_x\to \s^\infty_x}\leq B+ \frac{1}{u^{3/2}}\sup\limits_{|a|\leq 1}\| U(s+w, a+s+w)e^{iaH_0}\|_{\s_x^\infty\to \s_x^\infty}
\end{equation}

Then this part can be controlled by $\int_0^{c_1}dw C_2B$.  We choose $c_1$ small enough such that $C_2c_1<1/4$. Similarly, when there is a singularity for $ 1/(t-s-u)$, we use
\begin{align*}
&e^{-i[(t-s)-u]H_0}V(x,s+u)U(s+w,s+u)=[e^{-i[(t-s)-u]H_0}V(x,s+u)e^{i[(t-s)-u]H_0}] \times\\
&[e^{-i[(t-s)-u]H_0}U(s+w, s+w+u-(t-s))]U(s+w+u-(t-s),s+u)
\end{align*}
and then
\begin{align*}
\int_{t-s-c_2}^{t-s} du\int_0^{u}dw&\| A_3(u,w,s,t)\|_{\s^1_x\to \s^\infty_x}\leq  \int_{t-s-c_2}^{t-s} du\int_{c_1}^{u}dw \frac{C_3B}{|w|^{3/2}}+\\
&\int_{t-s-c_2}^{t-s} du\int_{u-1)}^{u}dw \frac{C_3}{(|t-s-w|^{3/2})|w|^{3/2}}\leq C_4(B+1)(c_2+c_2^{1/2}).
\end{align*}
Then we can choose $c_2$ small enough such that $C_4(c_2+c_2^{1/2})<1/4 $. If we have a singularity both for $1/w$ and $1/(t-s-u)$, then we use
\begin{align*}
&e^{-i[(t-s)-u]H_0}V(x,s+u)U(s+w,s+u)V(x,s+w)e^{-iwH_0}=[e^{-i[(t-s)-u]H_0}V(x,s+u)e^{i[(t-s)-u]H_0}] \times\\
&[e^{-i[(t-s)-u]H_0}U(s+w, s+w+u-(t-s))]U(s+w+u-(t-s),s+u+w)\times\\
&[U(s+u+w,s+u)e^{-iwH_0}][e^{iwH_0}V(x,s+w)e^{-iwH_0}].
\end{align*}
Then we get
\begin{equation}
\| A_3(u,w,s,t)\|_{\s^1_x\to \s^\infty_x}\leq \frac{C_5B}{|t-s|^{3/2}}\leq C_5B.
\end{equation}
Then we choose $c_3$ small enough in $\int_{t-s-c_3}^{t-s}du\int_0^{c_1}dw$ such that $ c_3c_1 C_5<1/4$.
So we have that for each pair $s,t$ with $|s-t|\geq 1$,
\begin{equation}
\| U(s,t)\|_{\s_x^1\to \s_x^\infty}\leq 3/4B+C.\label{sup}
\end{equation}
Take the supremum over $\{(s,t): |s-t|\geq 1\}$ on the left in equation \eqref{sup} and we have
\begin{equation}
B\leq 4C.
\end{equation}
Then the conclusion follows.
\end{proof}
Due to Lemma \ref{case1.1}, we have
\begin{align*}
&\sup\limits_{T\in \mathbb{R}}|T|^{3/2}\int_0^\infty \frac{dt_{k+|j|+1}}{(1+t_{k+|j|+1})^{3/2}}\times \\
&\| U(t_{k+|j|+1},0) e^{i\xi_{k+|j|}\cdot x}V(x-\sqrt{1+t_{k+|j|+1}}v)e^{i(T-t_{k+|j|+1})H_0}\|_{\s^1\to \s^\infty}<\infty
\end{align*}
where we have $1/(1+t_{k+|j|+1})^{3/2}$ since from $b_k+\cdots+b_l=5$, we gain $1/(1+t_{k+|j|+1}+s_{k+|j|})^{5/2}$. After taking the integral over $s_{t+|k|}$, we have $1/(1+t_{k+|j|+1})^{3/2}$. Hence,
\begin{equation}
\sup\limits_{T\in \mathbb{R}}|T|^{3/2}\|\text{case one}^{(k+|j|+1)}(T)\|_{\s_x^1\to \s_x^\infty}\leq \frac{(k+|j|)^5C^{k+|j|+1}|||V(x)|||_p^{k+|j|+1}}{\sqrt{M}^{k+|j|}}.
\end{equation}
Fix $|j|$. For case one, $l\in \{0,1,\cdots,|j|\}$ and for each $l$ and there are $\begin{pmatrix} 5+l\\ l\end{pmatrix}\leq 2^{5+|j|}$ many solutions of $(b_k,b_{k+1},\cdots, b_{k+l})\in \mathbb{N}^{l+1}$ satisfying
$$
b_k+b_{k+1}+\cdots+b_{k+1}=5.
$$
So for $k+|j|$, there are no more than $ j\times 2^{5+|j|}$ many case one terms. For case two, $l\in \{0,1,\cdots,|j|\}$ and for each $l$ and there are $\begin{pmatrix} b+l\\ l\end{pmatrix}\leq 2^{4+|j|}$ many solutions of $(b_k,b_{k+1},\cdots, b_{k+l})\in \mathbb{N}^{l+1}$ satisfying
$$
b_k+b_{k+1}+\cdots+b_{k+1}=b, \text{ for }b=0,1,2,3,4.
$$
So there are no more than $ 5j\times 2^{4+|j|}$ many case one terms. Thus,
\begin{align*}
\sup\limits_{T\in \mathbb{R}}|T|^{3/2}&\|\mathscr{J}_{M,T}^{(k+1)}\beta(|P|>32M) \|_{\s^1_x\to \s^\infty_x}\leq \sum\limits_{|j|=1}^\infty j\times 2^{5+|j|}\times \frac{(k+|j|)^5C^{k+|j|+1}|||V(x)|||_p^{k+|j|+1}}{\sqrt{M}^{k+|j|}} \\
&+5j\times 2^{4+|j|}\times \frac{(k+|j|)^4C^{k+|j|+1}|||V(x)|||_p^{k+|j|+1}}{\sqrt{M}^{k+|j|}}\leq \frac{k^5(C|||V(x)|||_p)^k}{\sqrt{M}^{k}}
\end{align*}
if $M$ is large enough.
\end{proof}
\begin{lemma}\label{KT}If $V(x-\sqrt{1+|t|}v)$ satisfies assumption \ref{sqrtc}, then
\begin{equation}
\sup\limits_{T\in \mathbb{R}}|T|^{3/2}\| \mathscr{K}^{(k)}(T)e^{iTH_0}\|_{\s^1_x\to \s^\infty_x}\leq \frac{(C|||V(x)|||_p)^{k}}{\sqrt{M}^{k-1}}, \text{ for }k\in\mathbb{N}^+. \label{KTk}
\end{equation}
\end{lemma}
\begin{proof}Apply Lemma \ref{refor2} and change of variables from $t_j \to t_{j}=s_j+\cdots+s_k$. For $\xi_j, s_j$, $j=1,\cdots,k-1$, it is the case when $\beta(|\xi_j+P|>2M)$. We do the same transformation as what we do in the proof of Lemma \ref{flowJ}. Then for each $j$, we will gain $C|||V(x)|||_p/\sqrt{M}$. For $s_k$, we apply Lemma \ref{D(T)} and then get the estimate \eqref{KTk}.
\end{proof}
Now we can prove its decay estimate. According to the definition of $\mathscr{J}_{M,\epsilon}^{(k+1)}, \mathscr{K}^{(k)}(T)$, we have
\begin{equation}
s\text{-}\lim\limits_{T\to \infty}\mathscr{D}(T):=s\text{-}\lim\limits_{T\to \infty} U(T,0)e^{-iTH_0}- \sum\limits_{k=1}^\infty i^{k+1} \mathscr{J}_{M,\epsilon}^{(k+1)}-\sum\limits_{k=1}^\infty i^{k}\mathscr{K}^{(k)}(T)=\mathbbm{1}.
\end{equation}
Then we have the following result.
\begin{lemma}If $V(x-\sqrt{1+|t|}v)$ satisfies assumption \ref{sqrtc}, we have
\begin{equation}
\sup\limits_{T\in \mathbb{R}^+} \|\mathscr{D}(T)\|_{\s_x^p\to \s_x^p}<\infty, \text{ for }1\leq p\leq\infty.
\end{equation}
\end{lemma}
\begin{proof}
The proof is the same as that of Corollary \ref{norm1}.
\end{proof}
Then the decay estimate follows.
\begin{proof}
For $T\geq 0$, it follows from
\begin{equation}
U(0,T)=\mathscr{D}(T)+\sum\limits_{k=1}^\infty i^{k+1} \mathscr{J}_{M,\epsilon}^{(k+1)}+\sum\limits_{k=1}^\infty i^{k}\mathscr{K}^{(k)}(T)
\end{equation}
and Lemma \ref{KT}, Lemma \ref{flowJ}. For $T<0$, it follows in the same way.
\end{proof}

\section{Application to NLS equations}\label{section 7}
\subsection{$\s^\infty$ boundedness for Hartree-type NLS}

We prove Theorem \ref{HartreeNLS} by proving an example. 
\subsubsection{$\s^\infty$ boundedness for some specific Hartree NLSs and the proof for Theorem \ref{HartreeNLS}}\label{GWS3}
In this section, we start with an example. Consider Hartree NLS equations
\begin{equation}
i\partial_t\psi(t)=H_0\psi(t)\pm\lambda [f*|\psi(t)|^2](x)\psi(t), \quad \psi(0)=\psi_0\text{ for }f(x,t)\in C_t\s^2_x\label{NLS1}.
\end{equation}
We prove Theorem \ref{2020LinftyNLS}. In other word, we show that $\psi(t)$ is bounded in $\s^\infty_x$ uniformly in $t\in (-\infty,-c]\cup [c,\infty)$ for any $c>0$ if $\psi_0\in \s_x^1\cap \s_x^2$. We reach this result by establishing its advanced CL: 
\begin{lemma}[Advanced CL]\label{2021can1}If $\psi(t)\in C_t([-T,T]) \s^2_x\cap \s^{8/3}_t([-T,T])\s^4_x $, then 
\eq
 \int_{-T}^T dt\|\sK_t(f*|\psi(t)|^2)\|_{\s^p_x\to \s^p_x}\lesssim T^{1/4} \|f(x)\|_{\s^2_x}\| \psi(t)\|_{\s^{8/3}_t([-T,T])\s^4_x}^2.\label{ineqcan1}
\eeq
In addition, 
\eq
 \int_{-T}^T dt\|\sK_t(f*|\psi(t)|^2)\|_{\s^4_t\s^p_x\to \s^p_x}\lesssim 1. \label{ineqcan}
\eeq
\end{lemma}
We defer the proof of Lemma \ref{2021can1} to the end of the section. We also have to show that the solution $\psi(t)$ to \eqref{NLS1} satisfies the assumption of Lemma \ref{2021can1}:
\begin{lemma}\label{2020NLSsol}If $\psi_0\in \s^2_x$, then for any $T>0, a\in \mathbb{R}$, 
\eq
\| \psi(t)\|_{ \s^{8/3}_t([-T+a,T+a])\s^4_x}\lesssim_{T, \|\psi_0\|_{\s^2_x}} 1. 
\eeq

\end{lemma}
The proof of Lemma \ref{2020NLSsol} is based on the construction of solution to \eqref{NLS1} by using CL and iteration scheme and we defer the proof to the end of this section. 

In the end, all result can be extended to the perturbed NLS. 

We are back to prove Theorem \ref{2020LinftyNLS}. We stick to $t\geq 0,$ $f(x,t)=f(x)$ and for $t<0$, the results follow from time reversal symmetry. The case for time-dependent $f$ will follow in the same way.
\begin{proof}[Proof of Theorem \ref{2020LinftyNLS}]We stick to $t\geq 1$ and the case for $t\geq c>0$ will follow in the same argument. By using Duhamel's formula, rewrite $\psi(t)$ as
\begin{multline}
\psi(t)=e^{-itH_0}\psi_0(x)+(-i)\int_0^{t-1/10}ds_1 e^{-i(t-s_1)H_0} [f*|\psi(s_1)|^2](x)\psi(s_1)+\\
(-i)\int_{t-1/10}^tds_1 e^{-i(t-s_1)H_0} [f*|\psi(s_1)|^2](x)\psi(s_1)\\
=:\psi_1(t)+\psi_2(t)+\psi_3(t).
\end{multline}
For $\psi_1(t),$ its $\s^\infty_x$ boundedness follows from the decay estimates of $e^{-itH_0}$. For $\psi_2(t)$, we have
\begin{multline}
\| \psi_2(t)\|_{\s^\infty_x}\lesssim \int_0^{t-1/10}ds_1 \frac{1}{|t-s_1|^{3/2}}\|  [f*|\psi(s_1)|^2](x) \psi(s_1)\|_{\s^1_x}\\
(\text{H\"older's inequality})\lesssim   \int_0^{t-1/10}ds_1 \frac{1}{|t-s_1|^{3/2}} \|   [f*|\psi(s_1)|^2](x)\|_{\s^2_x}\|\psi(s_1)\|_{\s^2_x}\\
\lesssim \int_0^{t-1/10}ds_1 \frac{1}{|t-s_1|^{3/2}}  \|f(x)\|_{\s^2_x}\| \psi(s_1)^2\|_{\s^1_x}\| \psi(s_1)\|_{\s^2_x}\\
(\text{H\"older's inequality})\lesssim \int_0^{t-1/10}ds_1 \frac{1}{|t-s_1|^{3/2}}  \|f(x)\|_{\s^2_x} \|\psi(s_1)\|_{\s^2_x}^3\\
\lesssim \int_0^{t-1/10}ds_1  \frac{1}{|t-s_1|^{3/2}} \|\psi_0\|_{\s^2_x}^3\\
\lesssim \|\psi_0\|_{\s^2_x}^3.
\end{multline}
For $\psi_3(t) $, we use Duhamel's formula again
\begin{multline}
\psi_3(t)=(-i)\int_{t-1/10}^t ds_1 e^{-i(t-s_1)H_0} [f*|\psi(s_1)|^2](x)e^{-is_1H_0}\psi_0(x)+\\
(-i)^2\int_{t-1/10}^t ds_1\int_0^{s_1-1/10}ds_2 e^{-i(t-s_1)H_0} [f*|\psi(s_1)|^2](x)e^{-i(s_1-s_2)H_0} [f*|\psi(s_2)|^2](x)  \psi(s_2)+\\
(-i)^2\int_{t-1/10}^t ds_1\int_{s_1-1/10}^{s_1}ds_2 e^{-i(t-s_1)H_0} [f*|\psi(s_1)|^2](x)e^{-i(s_1-s_2)H_0} [f*|\psi(s_2)|^2](x)  \psi(s_2)\\
=:\psi_{31}(t)+\psi_{32}(t)+\psi_{33}(t).
\end{multline}
For $\psi_{31}(t)$, using Lemma \ref{2021can1}, Lemma \ref{2020NLSsol} and the fact that $e^{-itH_0}\psi_0(x)\in \s^\infty_x$ for $t\geq a\geq\frac{1}{2}$, we have
\eq
\|\psi_{31}(t)\|_{\s^\infty_x}\lesssim_{\|\psi_0\|_{\s^2_x}} \|\psi_0\|_{\s^1_x}. 
\eeq
For $\psi_{32}(t)$, using Lemma \ref{2021can1}(regard $t-s_1$ variable as the time variable), Lemma \ref{2020NLSsol} and applying the same estimate for $\psi_2(t)$ to
\eq
\int_0^{s_1-1/10}ds_2 e^{-i(t-s_2)H_0} [f*|\psi(s_2)|^2](x)  \psi(s_2),
\eeq
we have
\eq
\|\psi_{32}(t)\|_{\s^\infty_x}\lesssim_{\|\psi_0\|_{\s^2_x}}1.
\eeq
For $\psi_{33}(t)$, we keep using Duhamel's formula in the same way twice. In the end, it is sufficient to deal with
\begin{multline}
\psi_4(t):=\int_{t-1/10}^t ds_1\int_{s_1-1/10}^{s_1}ds_2\int_{s_2-1/10}^{s_2}ds_3\int_{s_3-1/10}^{s_3}ds_4e^{-i(t-s_1)H_0} [f*|\psi(s_1)|^2](x)e^{-i(s_1-s_2)H_0} \cdots\\
e^{-i(s_3-s_4)H_0} [f*|\psi(s_4)|^2](x)\psi(s_4).
\end{multline}
We mainly use \eqref{ineqcan} in Lemma \ref{2021can1} 
\begin{multline}
\| \psi_4(t)\|_{\s^\infty_x}\lesssim \int_{t-1/10}^t ds_1\int_{s_1-1/10}^{s_1}ds_2\int_{s_2-1/10}^{s_2}ds_3\int_{s_3-1/10}^{s_3}ds_4\| \sK_{s_1-t}( [f*|\psi(s_1)|^2](x))\|_{\s^\infty_x\to\s^\infty_x }\times\\
\| \sK_{s_2-t}( [f*|\psi(s_2)|^2](x))\|_{\s^\infty_x\to\s^\infty_x }\| \sK_{s_3-t}( [f*|\psi(s_3)|^2](x))\|_{\s^\infty_x\to\s^\infty_x }\times\\
\|  e^{i(s_4-t)H_0} [f*|\psi(s_4)|^2](x)\psi(s_4)\|_{\s^\infty_x}\\
\lesssim_{\|\psi_0\|_{\s^2_x}} \int_{t-1/10}^t ds_1\int_{s_1-1/10}^{s_1}ds_2\int_{s_2-1/10}^{s_2}ds_3\int_{s_3-1/10}^{s_3}ds_4\| \sK_{s_1-t}( [f*|\psi(s_1)|^2](x))\|_{\s^\infty_x\to\s^\infty_x }\times\\
\| \sK_{s_2-t}( [f*|\psi(s_2)|^2](x))\|_{\s^\infty_x\to\s^\infty_x }\| \sK_{s_3-t}( [f*|\psi(s_3)|^2](x))\|_{\s^\infty_x\to\s^\infty_x }\times \frac{1}{|t-s_4|^{3/2}}\\
\lesssim_{\|\psi_0\|_{\s^2_x}}\int_{t-1/10}^t ds_1\int_{s_1-1/10}^{s_1}ds_2\int_{s_2-1/10}^{s_2}ds_3\| \sK_{s_1-t}( [f*|\psi(s_1)|^2](x))\|_{\s^\infty_x\to\s^\infty_x }\times\\
\| \sK_{s_2-t}( [f*|\psi(s_2)|^2](x))\|_{\s^\infty_x\to\s^\infty_x }\| \sK_{s_3-t}( [f*|\psi(s_3)|^2](x))\|_{\s^\infty_x\to\s^\infty_x }\frac{1}{|t-s_3|^{1/2}}\\
\lesssim_{\|\psi_0\|_{\s^2_x}}\int_{t-1/10}^t ds_1\int_{s_1-1/10}^{s_1}ds_2\| \sK_{s_1-t}( [f*|\psi(s_1)|^2](x))\|_{\s^\infty_x\to\s^\infty_x }\times\\
\| \sK_{s_2-t}( [f*|\psi(s_2)|^2](x))\|_{\s^\infty_x\to\s^\infty_x }\| \frac{\chi(s_3\in [s_2-1/10,s_2])}{|t-s_3|^{1/2}} \|_{\s^4_{s_3}}
\end{multline}
that is,
\begin{multline}
\| \psi_4(t)\|_{\s^\infty_x}\lesssim_{\|\psi_0\|_{\s^2_x}}\int_{t-1/10}^t ds_1\int_{s_1-1/10}^{s_1}ds_2\| \sK_{s_1-t}( [f*|\psi(s_1)|^2](x))\|_{\s^\infty_x\to\s^\infty_x }\times\\
\| \sK_{s_2-t}( [f*|\psi(s_2)|^2](x))\|_{\s^\infty_x\to\s^\infty_x }\frac{1}{|t-s_2|^{1/4}}\\
\lesssim _{\|\psi_0\|_{\s^2_x}}\int_{t-1/10}^t ds_1 \| \sK_{s_1-t}( [f*|\psi(s_1)|^2](x))\|_{\s^\infty_x\to\s^\infty_x }\| \frac{\chi(s_2\in [s_1-1/10,s_1])}{|t-s_2|^{1/4}} \|_{\s^4_{s_2}}\\
\lesssim_{ \|\psi_0\|_{\s^2_x}} 1.
\end{multline}
We finish the proof.

\end{proof}
Based on the proof of Theorem \ref{2020LinftyNLS}, we find that the proof only need the potential to be in $\s^2_x$ and it satisfies advanced CL. Thus, following a similar argument, we can extend the same result to a perturbed one:
\begin{proof}[Proof of Theorem \ref{HartreeNLS} part 1]\label{GWS0}If $\psi(t)$ exists in $\s^2_x$ and satisfies local Strichartz estimate, according to \ref{Vcondition1},\ref{Vcondition2},\eqref{Ncondition1}-\eqref{Ncondition3}, we  follow a similar argument of Theorem \ref{2020LinftyNLS} except that we may have to use Duhamel's formula for $N=[\frac{k_0'}{2}+1]+1$ times, in order to get the $\s^\infty_x$ boundedness result in Theorem \ref{HartreeNLS}, since when $N=[\frac{k_0'}{2}+1]+1$
\begin{multline}
\int_{t-1}^t ds_1\int_{t-1}^{s_1}ds_2\cdots \int_{t-1}^{s_{N-2}}ds_{N-1}|\int_{t-1}^{s_{N-1}}ds_N \frac{1}{|t-s_N|^{3/2}}|^{k_0'}\\
\lesssim \int_{t-1}^t ds_1\int_{t-1}^{s_1}ds_2\cdots \int_{t-1}^{s_{N-2}}ds_{N-1} \frac{1}{|t-s_{N-1}|^{\frac{k_0'}{2}}}\\
\lesssim_{k_0} \frac{1}{|t-s_1|^{\frac{k_0'}{2}-(N-1)}}\vert_{s_1=t-1}^{s_1=t}\lesssim_{k_0}1
\end{multline}
where
\eq
k_0=\min(k_1,k_2)\quad\text{ and }\quad\frac{1}{k_0'}+\frac{1}{k_0}=1.
\eeq
For $k_1,k_2$, see \ref{Vcondition1}, \eqref{Ncondition1}. So we have to show \eqref{generalNLS} has global well-posedness in $\s^2_x$ and local Strichartz estimate. We will show their proof in the following context, see \ref{GWS}.
\end{proof}

\begin{proof}[Proof of Lemma \ref{2021can1}] For \eqref{ineqcan1}, we only have to check that the Fourier transform of the potential is absolutely integrable,
\begin{multline}
\|\mathscr{F}[f*|\psi(t)|^2 ](\xi)\|_{\s^1_\xi}\sim\| \hat{f}(\xi) \mathscr{F}[|\psi(t)|^2](\xi) \|_{\s^1_\xi}\\
(\text{H\"older's inequality})\lesssim \|\hat{f}(\xi)\|_{\s^2_\xi}\| \mathscr{F}[|\psi(t)|^2](\xi)\|_{\s^2_\xi}\\
(\text{Plancherel theorem})\lesssim \|f(x)\|_{\s^2_x}\| \psi(t)\|_{\s^4_x}^2.
\end{multline}
Thus, 
\begin{multline}
\int_{-T}^T dt \| \sK_t( f*|\psi(t)|^2)\|_{\s^p_x\to \s^p_x}\lesssim \int_{-T}^T dt \|f(x)\|_{\s^2_x}\| \psi(t)\|_{\s^4_x}^2\\
(\text{H\"older's inequality})\lesssim  T^{1/4} \|f(x)\|_{\s^2_x}\| \psi(t)\|_{\s^{8/3}_t([-T,T])\s^4_x}^2.
\end{multline}
For \eqref{ineqcan1}, similarly, with $g(x,t)\in \s^4_t\s^p_x$, 
\begin{multline}
\int_{-T}^T dt \| \sK_t( f*|\psi(t)|^2)g(x,t)\|_{ \s^p_x}\lesssim \int_{-T}^T dt \|f(x)\|_{\s^2_x}\| \psi(t)\|_{\s^4_x}^2\|g(x,t)\|_{\s^p_x}\\
(\text{H\"older's inequality})\lesssim   \|f(x)\|_{\s^2_x}\| \psi(t)\|_{\s^{8/3}_t([-T,T])\s^4_x} \| g(x,t)\|_{\s^4_t\s^p_x}.
\end{multline}
We finish the proof.
\end{proof}
\begin{proof}[Proof of Lemma \ref{2020NLSsol}]It is sufficient to check the case when $a=0$ and $T>0$ sufficiently small. If we can get a boundedness only dependent on $\|\psi_0\|_{\s^2_x}$. Then we can extend the result to any other $a$ with the same $T$. For general finite $T>0$, we just have to use
\eq
\| \psi(t)\|_{\s^{8/3}_t([0,T])\s^4_x}\leq \sum\limits_{j=0}^N \| \psi(t)\|_{\s^{8/3}_t([T_j,T_{j+1}])\s^4_x}
\eeq
with $T_0=0, T_{N+1}=T$, where $N$ is sufficiently large number.

Now we go back to prove the case when $a=0$ and $T>0$ sufficiently small. It follows from an iteration scheme: set $\psi_1(t)=e^{-itH_0}\psi_0(x)$ and $\psi_{n+1}(t)$ satisfies 
\eq
\begin{cases}
i\partial_t\psi_{n+1}(t)=(-\Delta_x+f*|\psi_n(t)|^2)\psi_{n+1}(t)\\
\psi_{n+1}(0)=\psi_0
\end{cases},\quad t\in [0,T].
\eeq 
According to Lemma \ref{2021can1} and Strichartz estimates for $e^{itH_0}$, we have
\eq
\| \psi_{n+1}(t)\|_{\s^{8/3}_t \s^4_x([0,T]\times \mathbb{R}^3)}\leq \sum\limits_{j=0}^\infty \left(CT^{1/4} \|f(x)\|_{\s^2_x}\| \psi_n(t)\|_{\s^{8/3}_t([-T,T])\s^4_x}^2\right)^j\label{2021n1}
\eeq
and
\eq
\|\psi_{n+1}(t) \|_{\s^2_x}\leq \sum\limits_{j=0}^\infty \left(CT^{1/4} \|f(x)\|_{\s^2_x}\| \psi_n(t)\|_{\s^{8/3}_t([-T,T])\s^4_x}^2\right)^j\label{2021n2}
\eeq
for some constant $C>0$. From \eqref{2021n1}, we see if 
\eq
\| \psi_n(t)\|_{\s^{8/3}_t([-T,T])\s^4_x} \leq 2\|e^{-itH_0}\psi_0 \|_{\s^{8/3}_t\s^4_x}\leq 2C_{str}\|\psi_0\|_{\s^2_x}
\eeq
($C_{str}:=\|e^{itH_0}\|_{\s^2_x\to \s^{8/3}_t\s^4_x}$) and if we take $T>0$ small enough such that
\eq
4CT^{1/4} \|f(x)\|_{\s^2_x}C_{str}^2\|\psi_0\|_{\s^2_x}\leq \frac{1}{2},\label{2021T}
\eeq
then that
\eq
\| \psi_{n}(t)\|_{\s^{8/3}_t \s^4_x([0,T]\times \mathbb{R}^3)}\leq 2C_{str}\|\psi_0\|_{\s^2_x}
\eeq
implies 
\eq
\| \psi_{n+1}(t)\|_{\s^{8/3}_t \s^4_x([0,T]\times \mathbb{R}^3)}\leq 2C_{str}\|\psi_0\|_{\s^2_x}.
\eeq
Since
\eq
\| \psi_{1}(t)\|_{\s^{8/3}_t \s^4_x([0,T]\times \mathbb{R}^3)}\leq C_{str}\|\psi_0\|_{\s^2_x}\leq 2C_{str}\|\psi_0\|_{\s^2_x},
\eeq
we have for all $n=1,\cdots$, 
\eq
\| \psi_{n}(t)\|_{\s^{8/3}_t \s^4_x([0,T]\times \mathbb{R}^3)}\leq 2C_{str}\|\psi_0\|_{\s^2_x}\label{2021stri}
\eeq
if \eqref{2021T} is satisfied. Now we use standard contraction mapping argument to show $\psi_n$ converges both in $\s^2_x$ and $\s^{8/3}_t([0,T])\s^4_x$:
\begin{multline}
\| \psi_{n+1}(t)-\psi_n(t)\|_{\s^2_x}\leq \int_0^tds \|\sK_s( [f*|\psi_{n}(s)|^2](x))\|_{\s^2_x\to\s^2_x}\| \psi_{n+1}(s)-\psi_{n}(s)\|_{\s^2_x}+\\
\int_0^tds \|\sK_s( [f*(|\psi_{n}(s)|^2-|\psi_{n-1}(s)|^2)](x))\|_{\s^2_x\to\s^2_x}\| \psi_{n}(s)\|_{\s^2_x}\\
\leq CT^{1/4}\|f(x)\|_{\s^2_x}(2C_{str}\|\psi_0\|_{\s^2_x})^2 \sup\limits_{t\in [0,T]}\| \psi_{n+1}(t)-\psi_{n}(t) \|_{\s^2_x}+\\
CT^{1/4}\|f(x)\|_{\s^2_x}\times 4C_{str}\|\psi_0\|_{\s^2_x} \times 2\| \psi_0\|_{\s^2_x}\| \psi_{n}(t)-\psi_{n-1}(t) \|_{\s^{8/3}_t([0,T])\s^4_x}
\end{multline}
where we use
\eq
\| |\psi_{n}(t)|-|\psi_{n-1}(t)| \|_{\s^{8/3}_t([0,T])\s^4_x}\leq \| \psi_{n}(t)-\psi_{n-1}(t) \|_{\s^{8/3}_t([0,T])\s^4_x}.
\eeq
Then we have
\begin{multline}
\sup\limits_{t\in [0,T]}\| \psi_{n+1}(t)-\psi_n(t)\|_{\s^2_x}\leq CT^{1/4}\|f(x)\|_{\s^2_x}(2C_{str}\|\psi_0\|_{\s^2_x})^2 \sup\limits_{t\in [0,T]}\| \psi_{n+1}(t)-\psi_{n}(t) \|_{\s^2_x}+\\
CT^{1/4}\|f(x)\|_{\s^2_x}\times 4C_{str}\|\psi_0\|_{\s^2_x} \times 2\| \psi_0\|_{\s^2_x}\| \psi_{n}(t)-\psi_{n-1}(t) \|_{\s^{8/3}_t([0,T])\s^4_x}.\label{it1}
\end{multline}
Similarly, we have
\begin{multline}
\| \psi_{n+1}(t)-\psi_n(t)\|_{\s^{8/3}_t([0,T])\s^4_x}\leq C_{str}\int_0^tds \|\sK_s( [f*|\psi_{n}(s)|^2](x))\|_{\s^2_x\to\s^2_x}\| \psi_{n+1}(s)-\psi_{n}(s)\|_{\s^2_x}+\\
C_{str}\int_0^tds \|\sK_s( [f*(|\psi_{n}(s)|^2-|\psi_{n-1}(s)|^2)](x))\|_{\s^2_x\to\s^2_x}\| \psi_{n}(s)\|_{\s^2_x}\\
\leq C_{str}CT^{1/4}\|f(x)\|_{\s^2_x}(2C_{str}\|\psi_0\|_{\s^2_x})^2 \sup\limits_{t\in [0,T]}\| \psi_{n+1}(t)-\psi_{n}(t) \|_{\s^2_x}+\\
C_{str}CT^{1/4}\|f(x)\|_{\s^2_x}\times 4C_{str}\|\psi_0\|_{\s^2_x} \times 2\| \psi_0\|_{\s^2_x}\| \psi_{n}(t)-\psi_{n-1}(t) \|_{\s^{8/3}_t([0,T])\s^4_x}.\label{it2}
\end{multline}
Thus, by taking $T$ small enough such that we get
\eq
\sup\limits_{t\in [0,T]}\| \psi_{n+1}(t)-\psi_n(t)\|_{\s^2_x}\leq \frac{1}{3}\sup\limits_{t\in [0,T]}\| \psi_{n+1}(t)-\psi_n(t)\|_{\s^2_x}+ \frac{1}{3}\| \psi_{n}(t)-\psi_{n-1}(t)\|_{\s^{8/3}_t([0,T])\s^4_x}\label{it11}
\eeq
from \eqref{it1}, and 
\eq
\| \psi_{n+1}(t)-\psi_n(t)\|_{\s^{8/3}_t([0,T])\s^4_x}\leq \frac{1}{3}\sup\limits_{t\in [0,T]}\| \psi_{n+1}(t)-\psi_n(t)\|_{\s^2_x}+ \frac{1}{3}\| \psi_{n}(t)-\psi_{n-1}(t)\|_{\s^{8/3}_t([0,T])\s^4_x}\label{it12}
\eeq
from \eqref{it2}. Hence, according to \eqref{it11}, \eqref{it12}, we get
\eq
\| \psi_{n+1}(t)-\psi_n(t)\|_{\s^{8/3}_t([0,T])\s^4_x}\leq \frac{5}{6}\| \psi_{n}(t)-\psi_{n-1}(t)\|_{\s^{8/3}_t([0,T])\s^4_x}
\eeq
and
\eq
\sup\limits_{t\in [0,T]}\| \psi_{n+1}(t)-\psi_n(t)\|_{\s^2_x}\leq \frac{1}{2}\| \psi_{n}(t)-\psi_{n-1}(t)\|_{\s^{8/3}_t([0,T])\s^4_x}.
\eeq
Thus, by contraction mapping argument, we get that $\psi_n(t)$ converges to $\psi(t)$ in $\s^{8/3}_t([0,T])\s^4_x$ and therefore converges to $\psi(t)$ in $C_t([0,T])\s^2_x$. Thus,
\eq
\| \psi(t)\|_{\s^{8/3}_t([0,T])\s^4_x}\leq 2C_{str}\|\psi_0\|_{\s^2_x}
\eeq
due to \eqref{2021stri}. We finish the proof.

\end{proof}

\begin{proof}[Proof of Theorem \ref{HartreeNLS} part 2]\label{GWS} Based on the proof of Lemma \ref{2021can1} and Lemma \ref{2020NLSsol}, we can get the global wellposedness of \eqref{generalNLS} in $\s^2_x$(For $\s^2_x$, local wellposedness is equivalent to global wellposedness) and its local Strichartz estimates by using \ref{Vcondition1}, \eqref{Ncondition1} and \eqref{Ncondition2}. Here \ref{Vcondition1} is used to establish the local Strichartz estimates for $U_V(t,0)$ with $U_V(t,0)$, the semigroup generated by $H_0+V(x,t)$. We finish the proof of Theorem \ref{HartreeNLS}.
\end{proof}
\subsubsection{Typical examples}
Here are some typical examples: 
\begin{example}[Global wellposedness]\label{EX1}
When
\eq
\mathcal{N}(|\psi(t)|)=\pm \lambda[\frac{1}{|x|^{3/2-\delta}}*|\psi(t)|^2](x), \text{ for }\delta\in (0,\frac{3}{2}), \lambda>0,
\eeq
\eq
i\partial_t\psi(t)=(H_0+V(x,t))\psi(t)+\mathcal{N}(|\psi(t)|)\psi(t), \quad \psi(0)=\psi_0\label{GWPe1},
\eeq
with $V(x,t)$, satisfying \ref{Vcondition1}, \ref{Vcondition2},has global wellposedness in $\s^2_x$.
\end{example}
\begin{proof}
Compute its $\mathcal{F}\s^1_x$
\begin{multline}
\| \mathcal{N}(|\psi(t)|)\|_{\mathcal{F}\s^1_x}=\| \frac{1}{|\xi|^{3/2+\delta}}\mathscr{F}[|\psi(t)|^2](\xi) \|_{\s^1_\xi}\\
\leq \| \frac{\chi(|\xi|\leq 1)}{|\xi|^{3/2+\delta}}\mathscr{F}[|\psi(t)|^2](\xi) \|_{\s^1_\xi}+ \| \frac{\chi(|\xi|> 1)}{|\xi|^{3/2+\delta}}\mathscr{F}[|\psi(t)|^2](\xi) \|_{\s^1_\xi}\\
(\text{H\"older's inequality})\lesssim_\delta \|\psi(t) \|_{\s^2_x}^2+\|\frac{\chi(|\xi|> 1)}{|\xi|^{3/2+\delta}} \|_{\s^2_\xi}\|\mathscr{F}[|\psi(t)|^2](\xi) \|_{\s^2_\xi}\\
\lesssim_\delta \| \psi(t)\|_{\s^2_x}^2+\|\psi(t)\|_{\s^4_x}^2.
\end{multline}
Take $k_1=\frac{4}{3}$ and we have
\eq
\|\mathcal{N}(|\psi(t)|) \|_{\s^{4/3}_t([-T,T])\mathcal{F}\s^1_x}\lesssim_\delta \| \psi(t)\|_{C_t([-T,T])\s^2_x}^2+\|\psi(t)\|_{\s^{8/3}_t([-T,T])\s^4_x}^2.
\eeq
So \eqref{N1} is satisfied. Similarly,
\begin{multline}
\|\mathcal{N}(|\psi(t)|)-\mathcal{N}(|\phi(t)|)\|_{\mathcal{F}\s^1_x}= \| [\frac{1}{|x|^{3/2-\delta}}*(|\psi(t)|-|\phi(t)|)(|\psi(t)|+|\phi(t)|) ] \|_{\mathcal{F}\s^1_x}\\
\lesssim \| (|\psi(t)|-|\phi(t)|)(|\psi(t)|+|\phi(t)|)\|_{\s^1_x}+ \| (|\psi(t)|-|\phi(t)|)(|\psi(t)|+|\phi(t)|)\|_{\s^2_x}\\
\lesssim \| \psi(t)-\phi(t)\|_{\s^2_x}(\|\psi(t)\|_{\s^2_x}+\|\phi(t)\|_{\s^2_x})+\| \psi(t)-\phi(t)\|_{\s^4_x}(\|\psi(t)\|_{\s^4_x}+\|\phi(t)\|_{\s^4_x}).
\end{multline}
Then 
\begin{multline}
\int_{-T}^T dt \|\mathcal{N}(|\psi(t)|)-\mathcal{N}(|\phi(t)|)\|_{\mathcal{F}\s^1_x}\lesssim  T\| \psi(t)-\phi(t)\|_{C_t([-T,T])\s^2_x}(\|\psi(t)\|_{C_t(-T,T)\s^2_x}+\|\phi(t)\|_{C_t(-T,T)\s^2_x})\\
+T^{1/4}\| \psi(t)-\phi(t)\|_{\s^{8/3}_t([-T,T])\s^4_x}(\|\psi(t)\|_{\s^{8/3}_t([-T,T])\s^4_x}+\|\phi(t)\|_{\s^{8/3}_t([-T,T])\s^4_x}).
\end{multline}
So \eqref{N2} is satisfied. Thus, we have global wellposedness for \eqref{GWPe1}.
\end{proof}
\begin{example}[Global wellposedness and $\s^\infty$ boundedness]When
\eq
\mathcal{N}(|\psi(t)|)=\pm \lambda[\frac{e^{-c|x|}}{|x|^{3/2-\delta}}*|\psi(t)|^2](x), \text{ for }\delta\in (0,\frac{3}{2}), \lambda>0, c>0,
\eeq
\eq
i\partial_t\psi(t)=(H_0+V(x,t))\psi(t)+\mathcal{N}(|\psi(t)|)\psi(t), \quad \psi(0)=\psi_0\label{GWPe2},
\eeq
with $V(x,t)$, satisfying \ref{Vcondition1}, \ref{Vcondition2}, has global wellposedness in $\s^2_x$ and for any $c_0>0$, 
\eq
\sup\limits_{|t|\geq c_0}\|\psi(t)\|_{\s^\infty_x}\lesssim_{c_0, \|\psi_0\|_{\s^1_x\cap \s^2_x}}1.
\eeq
\end{example}
\begin{proof}Since
\eq
\mathcal{F}[\frac{e^{-c|x|}}{|x|^{\frac{3}{2}-\delta}}](\xi)\sim \frac{1}{\langle \xi\rangle^{\frac{3}{2}+\delta}},
\eeq
similarly, following the same estimate for Example \ref{EX1}, \eqref{N1}, \eqref{N2} are satisfied and we get global wellposedness in $\s^2_x$. In this case, according to H\"older's inequality, we have
\begin{multline}
\| \mathcal{N}(|\psi(t)|)\psi(t)\|_{\s^1_x}\lesssim \|  [\frac{e^{-c|x|}}{|x|^{3/2-\delta}}*|\psi(t)|^2](x)\|_{\s^2_x} \|\psi(t)\|_{\s^2_x}\\
\lesssim \|\frac{e^{-c|x|}}{|x|^{3/2-\delta}} \|_{\s^2_x} \|\psi(t)\|_{C([-T,T])\s^2_x}^3.
\end{multline}
So \ref{Ncondition3} is satisfied and we conclude \eqref{GWPe2} has global wellposedness in $\s^2_x$ and
\eq
\sup\limits_{|t|\geq c_0}\|\psi(t)\|_{\s^\infty_x}\lesssim_{c_0, \|\psi_0\|_{\s^1_x\cap \s^2_x}}1.
\eeq

\end{proof}
\subsection{Uniform $\s^p$ boundedness of wave operators for NLS equations for $2\leq p\leq\infty$}
In this section, we prove Theorem \ref{LinftyNLS} and Theorem \ref{LinftyNLSp}.

\subsubsection{$\s^\infty$ boundedness of $e^{itH_0}U(t,0)-1$}
We show $\s_x^\infty$ boundedness of $e^{itH_0}U(t,0)-1(\text{uniformly in }t\in [-\infty,\infty])$ on $\s^p_x\cap \mathcal{H}^1_x$ for $6<p\leq\infty$ by using the method of $ItT$ potential(ACL). If we only assume $\psi_0\in \mathcal{H}^1_x$ instead of $\psi_0\in \mathcal{H}^1_x\cap \s^p_x$, then $(e^{itH_0}U(t,0)-1)\psi_0$ is in $\s^\infty_x+\mathcal{F}\s^{1+\epsilon}_x$ for any $\epsilon\in (0,1)$, see Lemma \ref{Lem1}. As an application of Lemma \ref{Lem1}, we get a similar result for $U(t,0)-e^{-itH_0}$, see Corollary \ref{CLem1}. As an application of Theorem \ref{LinftyNLS}, we are able to get similar result for $U(t,0)$, see Lemma \ref{NLStinfty}. 
\begin{proof}[\textbf{Proof of Theorem \ref{LinftyNLS}}]Consider the $\s^\infty$ boundedness and begin with the case when $t=\infty$. Choose $\psi_0(x)\in \mathcal{H}^1_x$. Then due to \eqref{condition}, we have $\psi(t)\in \mathcal{H}^1_x$ uniformly in $t$. In the following context of the proof, $\psi(t)\in \mathcal{H}^1_x$ uniformly in $t\in \mathbb{R}$. We will give a proof for $\Omega_\pm^*-1$ and by replacing $\infty$ with $t$, we will get the same result for $e^{itH_0}U(t,0)-1$. According to Duhamel's formula, we have
\begin{multline}
i( \Omega_\pm^*-1)\psi_0(x)=\int_1^\infty ds e^{isH_0}\mathcal{N}({|\psi(s)|})\psi(s)+\int_0^1 ds e^{isH_0} [\beta(|P|> \frac{1}{s^{\frac{1}{2}+\frac{\epsilon}{2}}})\mathcal{N}(|\psi(s)|)]\psi(s)+\\
\int_0^1 ds e^{isH_0} [\beta(|P|\leq \frac{1}{s^{\frac{1}{2}+\frac{\epsilon}{2}}})\mathcal{N}({|\psi(s)|})]e^{-isH_0}\psi_0(x)+\int_0^1 ds e^{isH_0} [\beta(|P|\leq \frac{1}{s^{\frac{1}{2}+\frac{\epsilon}{2}}})\mathcal{N}({|\psi(s)|})]e^{-isH_0}\psi_1(s)\\
=:i\left[\J_{1}(\psi_0)+\J_{2}(\psi_0)+\J_{3}(\psi_0)+\J_{4}(\psi_0)\right],
\end{multline}
where 
\eq
\psi_1(s):=(-i)\int_0^sdu e^{iuH_0}\mathcal{N}({|\psi(s)|})\psi(u).
\eeq
For $\J_1(\psi_0)$, we have
\eq
\| \J_1(\psi_0) \|_{\s^\infty_x}\lesssim \int_1^\infty ds \frac{1}{s^{3/2}} \| \mathcal{N}({|\psi(s)|})\|_{\s^2_x}\| \psi(s)\|_{\mathcal{H}^1_x}\lesssim C(\|\psi_0(x)\|_{\mathcal{H}^1_x}).\label{nJ1}
\eeq
In order to estimate  
\eq
\int_0^1 ds e^{isH_0}\mathcal{N}({|\psi(s)|})\psi(s),
\eeq
we break it into $3$ pieces($\J_2(\psi_0), \J_3(\psi_0), \J_4(\psi_0)$) and estimate them separately. 

For $\J_2(\psi_0)$, we have
\begin{multline}
\|\J_{2}(\psi_0)\|_{\s^\infty_x}\lesssim \sum\limits_{l=1}^3\int_0^1ds\| e^{isH_0}\left[\frac{1}{P_l}\beta_l(|P|> \frac{1}{s^{\frac{1}{2}+\frac{\epsilon}{2}}})P_l[\mathcal{N}({|\psi(s)|})]\right]\psi(s)\|_{\s^\infty_x}\\
(\text{H\"older's inequality})\lesssim \sum\limits_{l=1}^3\int_0^1ds s^{\frac{1}{2}+\frac{\epsilon}{2}} \frac{1}{s^{3/2}}\| P_l[\mathcal{N}({|\psi(s)|})]\|_{\s^{6/5}_x}\|\psi(s)\|_{\s^{6}_x}\\
(\text{Since }\epsilon>0)\lesssim_\epsilon C(\|\psi_0(x)\|_{\mathcal{H}^1_x}) \label{nJ2}
\end{multline}
where $\epsilon>0$ will be chosen later(see \eqref{epsilon}, \eqref{epsilon2}), $\beta_l(P> \frac{1}{s^{\frac{1}{2}+\frac{\epsilon}{2}}})(l=1,2,3)$ is defined by 
\eq
\begin{cases}
\beta_1(P> \frac{1}{s^{\frac{1}{2}+\frac{\epsilon}{2}}}):=\beta(P_1 > \frac{1}{100s^{\frac{1}{2}+\frac{\epsilon}{2}}})\beta(P> \frac{1}{s^{\frac{1}{2}+\frac{\epsilon}{2}}})\\
\beta_2(P> \frac{1}{s^{\frac{1}{2}+\frac{\epsilon}{2}}}):= \beta(P_2 > \frac{1}{100s^{\frac{1}{2}+\frac{\epsilon}{2}}})\bar{\beta}(P_1 > \frac{1}{100s^{\frac{1}{2}+\frac{\epsilon}{2}}})\beta(P> \frac{1}{s^{\frac{1}{2}+\frac{\epsilon}{2}}})\\
\beta_3(P> \frac{1}{s^{\frac{1}{2}+\frac{\epsilon}{2}}}):=\beta(P_3 > \frac{1}{100s^{\frac{1}{2}+\frac{\epsilon}{2}}})\bar{\beta}(P_2 > \frac{1}{100s^{\frac{1}{2}+\frac{\epsilon}{2}}})\bar{\beta}(P_1 > \frac{1}{100s^{\frac{1}{2}+\frac{\epsilon}{2}}})\beta(P> \frac{1}{s^{\frac{1}{2}+\frac{\epsilon}{2}}})
\end{cases}.\label{2021bl}
\eeq
Here we also use
\eq
\|  \frac{1}{P_l}\beta(P_l> \frac{1}{100s^{\frac{1}{2}+\frac{\epsilon}{2}}})\|_{\s^{6/5}_x\to \s^{6/5}_x}\lesssim s^{\frac{1}{2}+\frac{\epsilon}{2}},
\eeq
see Lemma \ref{g}, and according to \eqref{condition},
\begin{multline}
\| P_l[\mathcal{N}({|\psi(s)|})]\|_{\s^{6/5}_x}\lesssim \| \mathcal{N}'({|\psi(s)|}) \times |P_l [\psi(s)]| \|_{\s^{6/5}_x}+\| \mathcal{N}'({|\psi(s)|}) \times |P_l [\psi^*(s)]| \|_{\s^{6/5}_x}\\
\lesssim C(\|\psi(s)\|_{\mathcal{H}^1_x})\lesssim  C(\|\psi_0(x)\|_{\mathcal{H}^1_x}).\label{H10}
\end{multline}
For $\J_{3}(\psi_0)$, we need the method of $ItT$.
\begin{lemma}[$ItT$ for NLS-1]\label{J3ItT}If $\psi_0\in \mathcal{H}^1_x\cap \s^{p}_x$ for some $p\in(6,\infty]$, then
\eq
\| \J_{3}(\psi_0)\|_{\s^\infty_x}\lesssim C(\|\psi_0(x)\|_{\mathcal{H}^1_x},\|\psi_0(x)\|_{\s^p_x}).
\eeq
\end{lemma}
\begin{proof}According to the standard computation for $tT$ potential, we have
\eq
\J_{3}(\psi_0)=\frac{1}{(2\pi)^{\frac{3}{2}}}\int_0^1ds \int d^3 \xi \beta(|\xi|\leq \frac{1}{s^{\frac{1}{2}+\frac{\epsilon}{2}}})\hat{V}(\xi,s) e^{i(x\cdot \xi+s\xi^2)} \psi_0(x+2s\xi)
\eeq
where
\eq
V(x,s):=\mathcal{N}({|\psi(s)|}).
\eeq
Control the $\s^\infty_x$ norm of $\J_{3}(\psi_0)$ directly 
\begin{multline}
\| \J_{3}(\psi_0)\|_{\s^\infty_x}\lesssim \sup\limits_{x\in \mathbb{R}^3}\int_0^1ds \int d^3\xi\beta(|\xi|\leq \frac{1}{s^{\frac{1}{2}+\frac{\epsilon}{2}}})|\hat{V}(\xi,s)| |\psi_0(x+2s\xi)|\\
(\text{H\"older's inequality})\lesssim  \sup\limits_{x\in \mathbb{R}^3}\int_0^1ds \| \beta(|\xi|\leq \frac{1}{s^{\frac{1}{2}+\frac{\epsilon}{2}}})\|_{\s^q_\xi} \|\hat{V}(\xi,s) \|_{\s^2_\xi} \|\psi_0(x+2s\xi)\|_{\s^p_\xi}\\
\lesssim \int_0^1ds \frac{1}{s^{\frac{3}{2q}+\frac{3\epsilon}{2q}}} C(\|\psi(s)\|_{\mathcal{H}^1_x})\|\psi_0(x)\|_{\s^p_x}\times \frac{1}{s^{3/p}}\\
\lesssim_{\epsilon, p} C(\|\psi_0(x)\|_{\mathcal{H}^1_x}) \|\psi_0(x)\|_{\s^p_x}
\end{multline}
where 
\eq
\frac{1}{p}+\frac{1}{q}=\frac{1}{2}
\eeq
and we use that
\begin{align}
\frac{3}{2q}+\frac{3\epsilon}{2q}+\frac{3}{p}=&\frac{3}{2}-\frac{3}{2q}+\frac{3\epsilon}{2q}\\
=&\frac{3}{2}-\frac{3}{2q}(1-\epsilon)\\
<&1
\end{align}
if we choose $\epsilon >0$ small enough such that
\eq
\frac{3}{q}(1-\epsilon)>1\label{epsilon}
\eeq 
and this can be achieved since $q<3$ due to $p>6$.

\end{proof}
According to Lemma \ref{J3ItT}, we have
\eq
\| \J_{3}(\psi_0)\|_{\s^\infty_x}\lesssim C(\|\psi_0(x)\|_{\mathcal{H}^1_x},\|\psi_0(x)\|_{\s^{p_0}_x}).\label{J3psi0}
\eeq
For $\J_4(\psi_0)$, we need following lemma:
\begin{lemma}\label{Lem1}If $\psi_0\in \mathcal{H}^1_x$ and $\mathcal{N}$ satisfies \eqref{condition}, then in \eqref{NLS}, for any $\epsilon_1\in (0,1)$, $\psi_1(s) \in \s^\infty_x+\mathcal{F}\s^{1+\epsilon_1}_x$ and its $ \s^\infty_x+\mathcal{F}\s^{1+\epsilon_1}_x$ norm is uniformly in $s\in \mathbb{R}$. To be precise,
\eq
\sup\limits_{s\in \mathbb{R}}\| \psi_1(s)  \|_{\s^\infty_x+\mathcal{F}\s^{1+\epsilon_1}_x}\lesssim_{\epsilon_1} C(\|\psi_0(x)\|_{\mathcal{H}^1_x}),
\eeq
that is,
\eq
\sup\limits_{s\in \mathbb{R}}\| e^{isH_0}\psi(s) -\psi_0 \|_{\s^\infty_x+\mathcal{F}\s^{1+\epsilon_1}_x}\lesssim_{\epsilon_1} C(\|\psi_0(x)\|_{\mathcal{H}^1_x}).
\eeq
\end{lemma}
\begin{proof}
Choose $\psi_0\in \mathcal{H}_x^1$. Due to the assumptions on $\mathcal{N}$, $\psi(t)\in \mathcal{H}^1_x$ uniformly in $t\in \mathbb{R}$. It is sufficient to look at
\eq
\J_{11}(\psi_0)(s):=\int_0^{\min\{1,s\} }du e^{iuH_0}\mathcal{N}({|\psi(u)|})\psi(u)
\eeq
since for $s\geq 1$, due to \eqref{condition}, 
\begin{multline}
\| \int_1^s due^{iuH_0}\mathcal{N}({|\psi(u)|})\psi(u) \|_{\s^\infty_x}\lesssim \int_1^s du\frac{1}{u^{3/2}}\| \mathcal{N}(|\psi(u)|)\|_{\s^2_x}\| \psi(u)\|_{\s^2_x} \\
 \int_1^\infty du \frac{1}{u^{3/2}} C(\| \psi(u)\|_{\mathcal{H}^1_x})\lesssim C(\|\psi_0(x)\|_{\mathcal{H}^1_x}).
\end{multline}
Break $\J_{11}(\psi_0)$ into two pieces
\begin{multline}
\J_{11}(\psi_0)(s)=\int_0^{\min\{1,s\}}du \beta(|P|> \frac{1}{u^{\frac{1}{2}+\frac{\epsilon_1}{2}}})e^{iuH_0}\mathcal{N}({|\psi(u)|})\psi(u)+\\
\int_0^{\min\{1,s\}}du \beta(|P|\leq \frac{1}{u^{\frac{1}{2}+\frac{\epsilon_1}{2}}})e^{iuH_0}\mathcal{N}({|\psi(u)|})\psi(u)\\
=:\J_{1,11}(\psi_0)(s)+\J_{1,12}(\psi_0)(s).
\end{multline}
For $\J_{1,11}(\psi_0)(s)$, we break $\beta(|P|> \frac{1}{u^{\frac{1}{2}+\frac{\epsilon_1}{2}}})$ into $3$ pieces
\eq
\beta(|P|> \frac{1}{u^{\frac{1}{2}+\frac{\epsilon_1}{2}}})=\sum\limits_{l=1}^3 \beta_l(|P|> \frac{1}{u^{\frac{1}{2}+\frac{\epsilon_1}{2}}}),\label{idd}
\eeq
where for $\beta_l$, see \eqref{2021bl}.

The $\s^\infty$ estimate for $\J_{1,11}(\psi_0)$ follows from, according to \eqref{condition},
\begin{multline}
\|\J_{1,11}(\psi_0)\|_{\s^\infty_x}\lesssim \sum\limits_{l=1}^3\int_0^1du\| \frac{1}{P_l}\beta_l(|P|> \frac{1}{u^{\frac{1}{2}+\frac{\epsilon_1}{2}}})e^{iuH_0}P_l[\mathcal{N}({|\psi(u)|})\psi(u)]\|_{\s^\infty_x}\\
\lesssim \sum\limits_{l=1}^3\int_0^1du u^{\frac{1}{2}+\frac{\epsilon_1}{2}} \frac{1}{u^{3/2}}\| P_l[\mathcal{N}({|\psi(u)|})\psi(u)]\|_{\s^1_x}\\
\lesssim_{\epsilon_1} C(\|\psi_0(x)\|_{\mathcal{H}^1_x})\label{2021J111}
\end{multline}
where we use
\eq
\|  \frac{1}{P_l}\beta(P_l> \frac{1}{100u^{\frac{1}{2}+\frac{\epsilon_1}{2}}})\|_{\s^\infty_x\to \s^\infty_x}\lesssim u^{\frac{1}{2}+\frac{\epsilon_1}{2}},
\eeq
and according to \eqref{H10}
\begin{multline}
\| P_l[\mathcal{N}({|\psi(u)|})\psi(u)]\|_{\s^1_x}\lesssim \| \mathcal{N}({|\psi(u)|}) \|_{\s^2_x}\times \| P_l [\psi(u)] \|_{\s^2_x}+\| P_l[\mathcal{N}'({|\psi(u)|})]\|_{\s^{6/5}_x} \times \|\psi(u) \|_{\s^6_x}\\
\lesssim C(\|\psi_0(x)\|_{\mathcal{H}^1_x}).\label{H1}
\end{multline}
For $ \J_{1,12}(\psi_0)$, compute its Fourier transform
\eq
\mathscr{F}[ \J_{1,12}(\psi_0) ](\xi)=\int_0^{\min\{1,s\}}du \beta(|\xi|\leq \frac{1}{u^{\frac{1}{2}+\frac{\epsilon_1}{2}}})e^{iu\xi^2}\hat{\phi}(\xi,u)
\eeq
with 
\eq
\phi(x,u):=\mathcal{N}(|\psi(u)|)\psi(u).
\eeq
Then
\begin{multline}
| \mathscr{F}[ \J_{1,12}(\psi_0) ](\xi)|\lesssim \sum\limits_{l=1}^3\int_0^1du \beta(|\xi|\leq \frac{1}{u^{\frac{1}{2}+\frac{\epsilon_1}{2}}})\beta_l(\xi)\times \frac{1}{|\xi|} |\xi_l \hat{\phi}(\xi,u) |\\
\lesssim \int_0^1du \beta(|\xi|\leq \frac{1}{u^{\frac{1}{2}+\frac{\epsilon_1}{2}}})\times \frac{1}{|\xi|}C(\|\psi_0(x)\|_{\mathcal{H}^1_x})\lesssim \frac{1}{|\xi|^{1+\frac{2}{1+\epsilon_1}}}C(\|\psi_0(x)\|_{\mathcal{H}^1_x})\in \s^1_\xi+\s^{1+\epsilon_1}_\xi
\end{multline}
where we use \eqref{H1} and 
\eq
|\xi_l \hat{\phi}(\xi,u) |\lesssim \| P_l[\mathcal{N}({|\psi(u)|})\psi(u)]\|_{\s^1_x}.\label{NLSxi}
\eeq
Thus, $\J_{1,12}(\psi_0)\in \s^\infty_x+\mathcal{F}\s^{1+\epsilon_1}_x$ and finish the proof.
\end{proof}
\begin{remark}\label{Dec.17.1}Here if in addition, $\psi_0\in \s^p_x$ for some $p\in [1,\frac{6}{5})$, then based on Lemma \ref{Lem1}, we have $\psi(t)\in \s^{p'}_x$($p'>6$ since $p<\frac{6}{5}$), which implies that in \eqref{NLSxi} $\xi_l \hat{\phi}(\xi,u)\in \s^{q}_{\xi}$ with $1/q+5/6+\frac{1}{p'}=1$. If we choose $\epsilon$ wisely, we are able to get $\mathscr{F}[ \J_{1,12}(\psi_0) ](\xi)\in \s^1_\xi$ and have $\psi(t)-e^{-itH_0}\psi_0\in \s^\infty_x$. For detailed statement, see Lemma \ref{NLStinfty}.
\end{remark}

\begin{corollary}\label{CLem1}If $\psi_0\in \mathcal{H}^1_x$ and $\mathcal{N}$ satisfies \eqref{condition}, then in \eqref{NLS}, for any $\epsilon_1\in (0,1)$, $\psi_1(s) \in \s^\infty_x+\mathcal{F}\s^{1+\epsilon_1}_x$ and its $ \s^\infty_x+\mathcal{F}\s^{1+\epsilon_1}_x$ norm is uniformly bounded in $s\in \mathbb{R}$. To be precise,
\eq
\sup\limits_{s\in \mathbb{R}}\| \psi(s)-e^{-isH_0}\psi_0  \|_{\s^\infty_x+\mathcal{F}\s^{1+\epsilon_1}_x}\lesssim_{\epsilon_1} C(\|\psi_0(x)\|_{\mathcal{H}^1_x}).
\eeq
\end{corollary}
According to Lemma \ref{Lem1}, by interpolation inequality, we have $\psi_1(x,t)\in \s^p_x$ for any $p\in [2,\infty)$ uniformly in $t$ and we get the $ItT$ potential method for $\J_4(\psi_0)$:
\begin{lemma}[$ItT$ for NLS-2]\label{J4ItT}If $\psi_0\in \mathcal{H}^1_x$, then
\eq
\| \J_{4}(\psi_0)\|_{\s^\infty_x}\lesssim C(\|\psi_0(x)\|_{\mathcal{H}^1_x}).
\eeq
\end{lemma}
\begin{proof}Similarly, we have
\begin{multline}
\| \J_{4}(\psi_0)\|_{\s^\infty_x}\lesssim \int_0^1ds \int d^3\xi\beta(|\xi|\leq \frac{1}{s^{\frac{1}{2}+\frac{\epsilon}{2}}})|\hat{V}(\xi,s)| |\psi_1(x+2s\xi,s)|\\
(\text{H\"older's inequality})\lesssim \int_0^1ds \| \beta(|\xi|\leq \frac{1}{s^{\frac{1}{2}+\frac{\epsilon}{2}}})\|_{\s^2_\xi\cap \s^{2+\epsilon_2}_\xi} \|\hat{V}(\xi,s) \|_{\s^2_\xi} \|\psi_1(x+2s\xi,s)\|_{\s^{\frac{1+\epsilon_1}{\epsilon}}_\xi+\s^\infty_\xi}\\
(\text{Lemma \ref{Lem1}})\lesssim \int_0^1ds \frac{1}{s^{\frac{3}{4}+\frac{3\epsilon}{4}}} \| \psi(s)\|_{\s^6_x}^3 C(\|\psi_0(x)\|_{\mathcal{H}^1_x})\times \frac{1}{s^{3\epsilon_1/(1+\epsilon_1)}}\\
(\text{Choosing $\epsilon, \epsilon_1$ sufficiently small})\lesssim\int_0^1ds \frac{1}{s^{7/8}} C(\|\psi_0(x)\|_{\mathcal{H}^1_x})\\
\lesssim C(\|\psi_0(x)\|_{\mathcal{H}^1_x}),\label{epsilon2}
\end{multline}
where
\eq
\frac{1}{2+\epsilon_2}+\frac{\epsilon_1}{1+\epsilon_1}=\frac{1}{2},
\eeq
$\epsilon_1\in (0,\frac{1}{4})$ and we also use that
\eq
\dfrac{1}{s^{3(\frac{1}{2}+\frac{\epsilon}{2})\times \frac{1}{2+\epsilon_2}}}\leq \frac{1}{s^{\frac{9}{20}+\frac{9\epsilon}{20}}}
\eeq
since
\eq
\frac{1}{2-\epsilon_2}=\frac{1}{2}-\frac{\epsilon_1}{1+\epsilon_1}>\frac{1}{2}-\frac{1/4}{5/4}=\frac{3}{10}.
\eeq
\end{proof}
According to \eqref{nJ1}, \eqref{nJ2}, \eqref{J3psi0} and Lemma \ref{J4ItT}, we get
\eq
\|( \Omega_\pm^*-1)\psi_0(x)\|_{\s^\infty_x}\lesssim C(\| \psi_0(x)\|_{\mathcal{H}^1_x\cap \s^p_x}).
\eeq
The $\s^\infty_x$ boundedness for $e^{itH_0}U(t,0)-1$ with $t\in [-\infty, \infty)$ follows in the same argument. Since for $t\in \mathbb{R}$, $e^{itH_0}U(t,0)-1: \mathcal{H}^1_x\to \s^2_x$, is bounded, by using interpolation inequality, we get 
\eq
\|(e^{itH_0}U(t,0)-1)\psi_0(x) \|_{\s^p_x}\leq C(\| \psi_0\|_{\mathcal{H}^1_x})
\eeq
for $p\in [2,\infty], t\in \mathbb{R}, \psi_0(x)\in \mathcal{H}^1_x\cap \s^{p_0}_x$. Now we come to $\s^p$ estimate of $\Omega_\pm^*$ for $p>6$ with additional assumption $\psi_0\in \s^p_x\cap \s^1_x$. Due to Lemma \ref{Lem1}, we have $\psi(t)\in \s^\infty_x+\mathcal{F}\s_x^{1+\epsilon}$ for $|t|\geq 1$, any $\epsilon>0$ if $\psi_0\in \s^1_x$. Then 
\begin{multline}
\| \int_1^\infty ds e^{isH_0} \mathcal{N}(|\psi(s)|)\psi(s) \|_{\s^p_x}\lesssim \int_1^\infty ds s^{-3(\frac{1}{2}-\frac{1}{p})}\| \mathcal{N}(|\psi(s)|) \|_{\s^2_x}\| \psi(s)\|_{\s^q_x}\\
\lesssim C(p, \|\psi_0(x)\|_{\mathcal{H}^1_x})\int_1^\infty ds s^{-3(\frac{1}{2}-\frac{1}{p})}\\
(\text{ use }p>6 \text{ and } \psi(s)\in \s^q \text{ due to interpolation})\lesssim C(p, \|\psi_0(x)\|_{\mathcal{H}^1_x})
\end{multline}
where $q$ satisfies 
\eq
\frac{1}{q}+\frac{1}{2}=\frac{1}{p'}.
\eeq
Thus, 
\eq
\|(\Omega_+^*-e^{iH_0}U(1,0))\psi_0(x)\|_{\s^p_x}\lesssim_p C(\|\psi_0\|_{\mathcal{H}^1_x})
\eeq
which implies that for $p\in (6,\infty]$(Recall that this time we have $\psi_0\in \s^p_x$),
\eq
\|\Omega_+^*\psi_0(x)\|_{\s^p_x}\lesssim_p C(\|\psi_0\|_{\mathcal{H}^1_x}).
\eeq
Similarly, we have the same result for $\Omega_-^*$ by using the a similar argument and finish the proof of Theorem \ref{LinftyNLS}.
\end{proof}
\begin{proof}[Proof of Theorem \ref{LinftyNLSp}]It follows directly from Lemma \ref{Lem1} since in Lemma \ref{Lem1}, we have $(e^{itH_0}U(t,0)-1)\psi_0\in \s^\infty+\mathcal{F}\s^{1+\epsilon}_x$ for any $\epsilon \in (0,1)$.
\end{proof}
We also have similar result for $U(t,0)-e^{-itH_0}$:
\begin{lemma}\label{NLStinfty}If $\psi_0\in\mathcal{H}^1_x$ and $\mathcal{N}$ satisfies \eqref{condition}, then for any $\epsilon\in (0,1)$, 
\eq
\sup\limits_{|t|\geq 1}\| \psi(t)-e^{-itH_0}\psi_0\|_{\s^\infty_x+\mathcal{F}\s^{1+\epsilon}_x}\leq C(\sup\limits_{t\in \mathbb{R}}\|\psi(t)\|_{\mathcal{H}^1_x },  \epsilon). \label{2021-1}
\eeq  
Furthermore, if $\psi_0\in \s^p_x\cap \mathcal{H}^1_x$ for some $p\in [1,\frac{6}{5})$ and 
\eq
\sup\limits_{t\in \mathbb{R}}\|\psi(t)\|_{\mathcal{H}^1_x} \lesssim 1,
\eeq
then
\eq
\sup\limits_{|t|\geq 1}\| \psi(t)-e^{-itH_0}\psi_0\|_{\s^\infty_x}\leq C(\sup\limits_{t\in \mathbb{R}}\|\psi(t)\|_{\mathcal{H}^1_x }, \|\psi_0\|_{\s^p_x}, p'-6). \label{2021-2}
\eeq 
\end{lemma}
\begin{proof}[Proof of Lemma \ref{NLStinfty}]For \eqref{2021-1}, it follows by using a similar argument as what we did in Lemma \ref{Lem1}. For \eqref{2021-2}, by using Duhamel's formula, write $\psi(t)$ 
\eq
\psi(t)=e^{-itH_0}\psi_0+(-i)\int_0^tds e^{-i(t-s)H_0}\mathcal{N}(|\psi(s)|)\psi(s).
\eeq
For $\s^\infty_x$ estimate, it is sufficient to estimate
\eq
\psi_2(t):=(-i)\int_{t-\frac{1}{2}}^tds e^{-i(t-s)H_0}\mathcal{N}(|\psi(s)|)\psi(s).
\eeq
Since $\psi_0\in \s^p_x$ implies $e^{-itH_0}\psi_0\in \s^{p'}_x$ for $p'>6, t\geq \frac{1}{2}$, by using a similar argument as what we did in the proof of Theorem \ref{LinftyNLS} and due to Remark \ref{Dec.17.1}, we get \eqref{2021-2}.

\end{proof}
\subsubsection{Typical examples and remarks on advanced cancelation lemma}
\begin{example}[$\s^\infty$ boundedness(Cubic NLS)]\label{EX3}
When
\eq
\mathcal{N}(|\psi(t)|)= |\psi(t)|^3,
\eeq
\eq
i\partial_t\psi(t)=H_0\psi(t)+\mathcal{N}(|\psi(t)|)\psi(t), \quad \psi(0)=\psi_0\in \s^1_x\cap \mathcal{H}^1_x\label{cubic},
\eeq
satisfies \eqref{condition}. Then 
\eq
\sup\limits_{|t|\geq 1}\|\psi(t) \|_{\s^\infty_x}\lesssim 1.\label{conex3}
\eeq
\end{example}
\begin{proof}
When 
\eq
\mathcal{N}(|\psi(t)|)= |\psi(t)|^3,
\eeq
it is the defocusing case and if $\psi_0\in \mathbb{H}^1_x$, we have a global solution $\psi(t)$ with a uniform $\mathcal{H}^1_x$ norm. We also have
\eq
\|\mathcal{N}(|\psi(t)|) \|_{\s^2_x}=\| |\psi(t)|^3\|_{\s^2_x}=\|\psi(t)\|_{\s^6_x}^3\lesssim \|\psi(t)\|_{\mathcal{H}^1_x}^3
\eeq
and
\eq
\|\mathcal{N}'(|\psi(t)|) \|_{\s^3_x}=3\| |\psi(t)|^2\|_{\s^3_x}=3\|\psi(t)\|_{\s^6_x}^2\lesssim \|\psi(t)\|_{\mathcal{H}^1_x}^2.
\eeq
So \eqref{condition} is satisfied and we have \eqref{conex3}.
\end{proof}
\begin{example}[$\s^\infty$ boundedness of mixed power nonlinearity]\label{EX4}
When
\eq
\mathcal{N}(|\psi(t)|)=- |\psi(t)|^2+|\psi(t)|^3,
\eeq
if $\psi(t)\in \mathcal{H}^1_x$, uniformly in $t$, then
\eq
i\partial_t\psi(t)=H_0\psi(t)+\mathcal{N}(|\psi(t)|)\psi(t), \quad \psi(0)=\psi_0\in \s^1_x\cap \mathcal{H}^1_x\label{cubic'},
\eeq
satisfies \eqref{condition}. Then 
\eq
\sup\limits_{|t|\geq 1}\|\psi(t) \|_{\s^\infty_x}\lesssim 1.\label{conex4}
\eeq
\end{example}
\begin{proof}
When 
\eq
\mathcal{N}(|\psi(t)|)=- |\psi(t)|^2+|\psi(t)|^3,
\eeq
according to Lemma \ref{TVZ}, we have
\eq
\| \psi(t)\|_{\mathcal{H}^1_x}\leq \|\psi_0\|_{\s^2_x}+\sup\limits_{s\in [t-1,t+1]} \|\nabla \psi(t)\|_{\s^2_x}\lesssim C(\|\psi_0\|_{\mathcal{H}^1_x}).
\eeq
We also have
\eq
\|\mathcal{N}(|\psi(t)|) \|_{\s^2_x}\lesssim\| |\psi(t)|^2\|_{\s^2_x}+\| |\psi(t)|^3\|_{\s^2_x}=\|\psi(t)\|_{\s^6_x}^3+\|\psi(t)\|_{\s^4_x}^2\lesssim C(\|\psi(t)\|_{\mathcal{H}^1_x})
\eeq
and
\eq
\|\mathcal{N}'(|\psi(t)|) \|_{\s^3_x}\lesssim 2\| |\psi(t)|\|_{\s^3_x}+3\| |\psi(t)|^2\|_{\s^3_x}\lesssim C(\|\psi(t)\|_{\mathcal{H}^1_x}).
\eeq
So \eqref{condition} is satisfied, and we have \eqref{conex4}.
\end{proof}

\section{Intertwining property}
In the time-independent case, there exists an intertwining between $f(H)$ and $f(H_0)$ with $f$ measurable
\begin{equation}
f(H)P_c=\Omega_+f(H_0)\Omega_+^* \label{A.19.1}
\end{equation}
where $P_c$ denotes the projection on the continuous spectrum of $H,$ and this projection comes from the fact that $\Omega_+$ is unitary from $L^2\to Ran(\Omega_+),$ with the range of $\Omega_+$ equal to the continuous spectrum of $H$.

 When it comes to the time-dependent case, \eqref{A.19.1} fails in most situation in that $U(t+s,t)$ will not generally have a nice limit as $t\to\infty$, see \cite{RS1979}. In this section, we will introduce an adapted type of intertwining property based on the time-dependent wave operators $\Omega_T$ (For $\Omega_T$, see \eqref{OmegaT}.)
\begin{equation}
U(T,0)=\Omega_T e^{-iTH_0}\Omega^*_+, \text{ on }\mathcal{R}(\Omega_+)\label{A.23.1}
\end{equation}
where $U(t,0)$ denotes the solution operator of a Schr\"odinger equation with a Hamiltonian $H(t)$, $\mathcal{R}(\Omega)$ is the range of $\Omega_+$, a subspace equipped with the $\s^2$ norm. It follows from
\begin{equation}
U(T,0)=U( T,T+s)U( T+s,0)=U(T,T+s)e^{-isH_0}e^{-iTH_0}e^{i(s+T)H_0}U( T+s,0), \text{ on }\s^2,
\end{equation}
\begin{equation}
\Omega_T=s\text{-}\lim\limits_{s\to \infty} U(T,T+s)e^{-isH_0}, \text{ on }\s^2
\end{equation}
and
\begin{equation}
\Omega_+=s\text{-}\lim\limits_{s\to \infty} U( 0,s)e^{-isH_0}, \text{ on }\s^2.
\end{equation}
Based on Corollary \ref{cor3} and Theorem \ref{main3}, we have
\begin{equation}
\|\Omega_T e^{-iTH_0}\beta(|P|>M)\Omega^*_+\|_{ \s^p\to \s^{p^\prime} }\lesssim \frac{1}{T^{3(\frac{1}{p}-\frac{1}{2})}}, \text{ in dimension }3\label{A.27.1.1}
\end{equation}
with $\frac{1}{p}+\frac{1}{p^\prime}=1, 1\leq p\leq 2$. The decay estimates follow if we make a  low-frequency assumption:
\begin{lemma}If
\begin{equation}
\|\Omega_T e^{-iTH_0}\beta(|P|\leq M)\Omega^*_+\|_{ \s^p\to \s^{p^\prime} }\lesssim \frac{1}{T^{3(\frac{1}{p}-\frac{1}{2})}}\label{A.27.1.2}
\end{equation}
for $1\leq p\leq 2$, some sufficiently large $M$, and if $V(x,t)$ satisfies the condition in Theorem \ref{main3}, then $U(T,0)$ satisfies decay estimates on $\mathcal{R}(\Omega_+)\cap L^p_x$ for all $T>0$.
\end{lemma}
\begin{proof} Based on Corollary \ref{cor3} and Theorem \ref{main3}, we have \eqref{A.27.1.1}. Then combining \eqref{A.27.1.1} with assumption \eqref{A.27.1.2}, we get
\begin{equation}
\|\Omega_T e^{-iTH_0}\Omega^*_+\|_{ \s^p\to \s^{p^\prime} }\lesssim \frac{1}{T^{3(\frac{1}{p}-\frac{1}{2})}}.
\end{equation}
Based on \eqref{A.23.1}, we get
\begin{equation}
\sup\limits_{T\geq 0}T^{3/2}\| U(T,0)\|_{ \mathcal{R}(\Omega_+)\cap \s^1_x\to \s^\infty_x}\lesssim 1.
\end{equation}
Later by interpolation, we get $\s^p$ decay estimates on $\mathcal{R}(\Omega_+)\cap \s^p_x$.

\end{proof}
 More information about intertwining property will be discussed in our following paper.

\footnote{A. Soffer is supported in part by NSF grant DMS-1600749 and by NSFC11671163 }

\end{document}